\documentclass{article}
\usepackage{graphicx} % Required for inserting images
\usepackage{hyperref}  % or \usepackage{url}
\usepackage{lineno}
\usepackage{url}

\usepackage{amsmath, amsthm, amssymb,amsfonts}
\usepackage[linesnumbered,ruled,vlined]{algorithm2e}

\usepackage{mathrsfs}
\usepackage{scalerel}
\usepackage[dvipsnames]{xcolor}
\usepackage[margin=1in]{geometry}
\usepackage[noblocks]{authblk}
\usepackage{ulem}
\theoremstyle{plain}

\theoremstyle{plainNoItalics}

\SetCommentSty{mycommfont}
\usepackage{xifthen}
\usepackage{bbm, bm}
\usepackage[dvipsnames]{xcolor}
\usepackage{bm} % Add this in your preamble
\definecolor{darkgreen}{rgb}{0.0, 0.5, 0}
\definecolor{darkred}{rgb}{0.7, 0.0, 0.1}
\usepackage{soul} %% to cross word
\newcommand{\mathbfit}[1]{\bm{\mathit{#1}}}

\newcommand\HHREV[1]{\textcolor{black}{#1}}
\newcommand\HHREVONE[1]{\textcolor{black}{#1}}
\newcommand\HHONE[1]{\textcolor{black}{#1}}
\newcommand\HHTWO[1]{\textcolor{black}{#1}}

\title{Sylvester-Preconditioned Adaptive-Rank Implicit Time Integrators for Advection-Diffusion Equations with \HHONE{Variable} Coefficients}

\author[a]{Hamad El Kahza}
\author[a]{Jing-Mei Qiu}
\author[b]{Luis Chac\'on}
\author[b]{William Taitano}

\affil[a]{\small{Department of Mathematical Sciences, University of Delaware, Newark, DE 19716}}
\affil[b]{\small{Theoretical Division, Los Alamos National Laboratory, Los Alamos, NM 87545}}

\newcommand{\bF}{{\bf F}}
\newcommand{\bB}{{\bf B}}

\newcommand{\bS}{{\bf S}}
\newcommand{\bR}{{\bf R}}

\begin{document}

\maketitle
\begin{abstract}
We consider the adaptive-rank integration of multi-dimensional time-dependent advection-diffusion partial differential equations (PDEs) with variable coefficients. We employ a standard finite-difference method for spatial discretization coupled with high-order diagonally implicit Runge-Kutta temporal schemes. The discrete equation is a generalized Sylvester equation (GSE), which we solve with a projection-based adaptive-rank algorithm structured around two key strategies: (i) constructing dimension-wise subspaces using a novel atypical extended Krylov strategy, and (ii) efficiently solving the basis coefficient matrix with a preconditioned GMRES solver. The low-rank decomposition is performed in 2D using SVD and with high-order SVD (HOSVD) in 3D to represent the tensor in a compressed Tucker format.
\HHREV{For $d$-dimensional problems (here, $d = 2$ or $3$), the computational complexity and memory storage of the approach are found numerically to scale as $\mathcal{O}(N {r^2}) + \mathcal{O}({r^{d+1}})$  and $\mathcal{O}(N r)+ \mathcal{O}(r^{d})$, respectively,
with $N$ the one-dimensional resolution and $r$} the \HHONE{maximal} rank during the Krylov iteration (which we find to be largely independent of $N$ on our numerical examples). We present numerical examples that illustrate the advertised properties of the algorithm.
\end{abstract}

\section{Introduction}
% \lc{Most of this section should be in green color...}
We propose an adaptive-rank numerical scheme for approximating solutions of multi-dimensional \HHTWO{time-dependent} variable-coefficient advection-diffusion equations. {The goal is to reduce the computational complexity of high-dimensional PDE solvers by {dynamically} exploiting the low-rank structure of the {time-dependent} solution when it exists, building upon existing literature on low-rank PDE solvers \cite{grasedyck2013literature, einkemmer2024review} and model order reduction (MOR) for solving high-dimensional linear and nonlinear systems \cite{kressner2010krylov, grasedyck2013literature, benner2013low, simoncini2016computational}.} 

{\color{black}
Two primary strategies in the literature exploit evolving low-rank structures for time-dependent PDEs: the dynamical low-rank (DLR) approximation~\cite{koch2007dynamical, lubich2014projector,
ceruti2022unconventional, dektor2021dynamic} and the step-and-truncate (SAT) method, which has both explicit~\cite{kormann2015semi, dektor2021rank, guo2022low, guo2023local, guo2022conservative} and implicit~\cite{rodgers2023implicit, nakao2023reduced, el2024krylov, meng2024preconditioning, nakao2025reduced} variants. Recently, \cite{einkemmer2024review} provided an overview of these two approaches with a focus on time-dependent kinetic simulations. 
{Our study employs the SAT approach, which is known for preserving the high-order spatial and temporal accuracy of traditional mesh-based solvers for time-dependent problems. Specifically, to capture the solution dynamics accurately, we apply a high-order diagonally implicit Runge-Kutta (DIRK) time discretization to the matrix differential equation arising from the semi-discretization of advection-diffusion equations on a tensor-product of 1D spatial grids. The backward Euler time discretization of a 2D constant-coefficient advection-diffusion equation yields a Sylvester equation (SE) of the form:
\begin{equation}
A \mathbf{X}_1 + \mathbf{X}_1 B^\top = \mathbf{X}_0,
\label{eq: Se}
\end{equation}
where $\mathbf{X}_0, \mathbf{X}_1 \in \mathbb{R}^{N_1 \times N_2}$ represent the solution snapshots at times $t_0$ and $t_1$, respectively, and $A$ and $B$ are coefficient matrices derived from the discretizations of advection and diffusion differential operators in two dimensions. Here $N_1$ and $N_2$ denote the number of grid points along each coordinate direction of the tensor product grid. Generalization to problems with variable advection and/or diffusion coefficients leads to the so-called generalized Sylvester equation (GSE) in matrix and Kronecker forms, respectively:
\begin{equation}
    \sum_{i=1}^{l} A_i \mathbf{X}_1 {B}^{\top}_i =\mathbf{X}_0 \quad \iff \quad \underbrace{\left(\sum_{i=1}^{l} {B}_i\otimes A_i\right)}_{\mathcal{A}}\text{vec}(\mathbf{X}_1)=\text{vec}(\mathbf{X_0}),
    \label{generalized_eq}
\end{equation}
where $A_i$ and $B_i$ are coefficient matrices from discretization of the PDE operators, leveraging the separability of variable advection-diffusion coefficients. Higher-order DIRK schemes lead to SE and GSE in forms similar to \eqref{eq: Se}-\eqref{generalized_eq}, but with the right-hand side (RHS) term, $\mathbf{X}_0$, constructed from the previous DIRK stages. 
To exploit low-rank structure in matrices (for 2D) and tensors (for 3D and beyond), such RHS term can be efficiently compressed using SVD rounding in 2D, Tucker decomposition in 3D, and hierarchical Tucker or tensor train decomposition for problems beyond 3D \cite{oseledets2011tensor, kressner2012htucker}. Ultimately, a successful development of efficient and robust solvers for the GSE \eqref{generalized_eq} is of critical importance to enable adaptive-rank integrators that accurately capture the transient states of time-dependent problems.}

{Two broad families of algorithms are commonly employed to tackle GSEs (and SEs as a special case): \emph{projection-type methods} \cite{benner2013low,hao2021sherman,simoncini2023analysis,powell2017efficient,saad1989numerical, kressner2010krylov, simoncini2016computational,shank2016efficient} and \emph{matrix-oriented classical Krylov methods} \cite{simoncini2023analysis,palitta2016matrix,benner2013low,palitta2021convergence,palitta2025subspace,kressner2011param}. Projection-type schemes first construct low-dimensional subspaces that capture the dominant features of the solution. The GSE is then projected {dimension-by-dimension} via a Galerkin condition, yielding a much smaller GSE for the expansion coefficients that compose the bases together. Matrix-oriented Krylov iterations, in contrast, construct an approximation subspace by applying powers of the operator
\(\mathcal{A}\) in~\eqref{generalized_eq} to the vectorized RHS, but exploit the Kronecker structure so that the action of \(\mathcal{A}\) is performed at the {matrix level}, and the full Kroneckerized operator is never formed explicitly. Inner products and other Krylov routines are likewise executed at the matrix-level, with low-rank truncation applied throughout to control memory. Although an appropriately chosen projection space can make a projection-type method formally equivalent to a matrix-oriented Krylov scheme, the projection-based framework offers additional freedom to design {nonstandard} subspaces that may accelerate convergence \cite{simoncini2023analysis}. On the downside, projection methods are less amenable to precondition: their efficiency hinges on the efficacy of the chosen subspaces and on the accuracy with which the reduced problem is solved. Classical Krylov methods, in contrast, benefit from well-established preconditioning strategies but offer no flexibility to alter the search space, as the basis is typically generated by an Arnoldi or Lanczos-type process applied directly to \(\mathcal{A}\). In summary, each approach has complementary strengths and weaknesses, and the choice between them should be guided by the structure of the particular GSE under investigation. Below we survey representative work in both directions.} 

In the context of solving SEs, projection-type methods based on Krylov subspace techniques have proven highly effective, particularly by exploiting the low-rank structure of PDE solutions~\cite{kressner2010krylov, simoncini2007new, rodgers2023implicit, breiten2014fast,el2024krylov}. These methods typically construct dimension-wise solution bases from Krylov or Krylov-like subspaces and apply a Galerkin projection onto the reduced subspace spanned by the outer product of dimension-wise Krylov subspaces~\cite{saad1989numerical, kressner2010krylov, simoncini2016computational}. To accelerate convergence, extended Krylov subspace methods (first proposed in~\cite{druskin1998extended} for the approximation of matrix functions)  that incorporate $A$ and $A^{-1}$, as well as $B$ and $B^{-1}$, in each direction, respectively, were successfully applied to matrix equations in~\cite{simoncini2007new, el2024krylov} within a projection-based scheme termed Krylov-Plus-Inverted-Krylov (KPIK). 

Applying projection-type methods to GSEs is much more challenging. The major difficulty lies in preventing the rapid growth of the solution approximation spaces due to the (possibly large) multi-term structure of the equation, which otherwise results in large reduced systems (of size $r^d$, where $r$ is the dimension of the per-direction Krylov subspace and $d$ is the number of dimensions).
% \QQ{of size $r^d$. Here and below, $r$ is the dimension of the per-direction Krylov subspace and $d$ is the number of dimensions. (we introduce $r$ and $d$ multiple times later...)} \lc{Yes... we should remove}
Such large reduced system demands efficient, {tailored} solution strategies. %\QQ{How about "Major difficulties lie in preventing the rapid growth of the solution approximation spaces due to the (possibly large) multi-term structure of the equation, which otherwise result in large reduced systems (of size $r^d$, where $r$ is the dimension of the per-direction Krylov subspace and $d$ is the number of dimensions); such reduced system demands efficient solver."}
Generalized Lyapunov equations of the form ($A\mathbf{X}_1+\mathbf{X}_1 A +\sum_{i=1}^lM_iX_1M_i=\mathbf{X}_0$), a special case of the GSE matrix equation \eqref{generalized_eq}, have received significant attention ~\cite{benner2013low, hao2021sherman} {exemplifying such challenges}.  In~\cite{benner2013low}, projection spaces are built from extended Krylov subspaces from the operator $A$ and standard Krylov subspaces from the remaining operators ($M_1,\cdots M_l$), but do not address how to solve the reduced system efficiently.
% \QQ{but ``how to solve the reduced system efficiently" was not addressed}. \lc{... and this} 
In~\cite{hao2021sherman}, for the case $l=1$, projection spaces are built from a rational Krylov from $A$ and inverted Krylov from $A+M_1$. The authors then use the Sherman-Woodbury formula to express the projected equation as a series of Lyapunov solves, relying on the availability of a low-rank decomposition of the Galerkin projected generalized operators in the system. The resulting complexity for the reduced system is $\mathcal{O}(r^{d+1})$ under the assumption of fast decaying singular values of the operators, but deteriorates to $\mathcal{O}(r^{3d})$ otherwise.
% In~\cite{powell2017efficient}, a low-rank (one-sided) projection-based approach is applied to problems with disparate dimension sizes $(N_2 \ll N_1)$, using block rational Krylov from all operators in the $N_1$ direction. 

% In~\cite{powell2017efficient}, solving the small projected system with either matrix-based Conjugate Gradient or a Kronecker formulation, with efficiency owing to the  disparate dimensions assumption.
%In summary, projection-based methods for GSEs must carefully tailor the construction of approximation spaces to the specific structure of the equation, with overall efficiency depending critically on both the quality of these subspaces and the ability to solve the reduced system effectively. 

Low-rank matrix-oriented Krylov methods
%The second approach to solving the Kroneckerized form of the GSEs involves low-rank matrix-oriented formulations of classical Krylov methods, 
such as Conjugate Gradient (CG) or GMRES~\cite{benner2013low,kressner2011param,palitta2025subspace,stoll2015low} provide an alternative to projection-type solvers. {In~\cite{benner2013low}, several low-rank matrix-oriented Krylov methods, including CG and BiCGSTAB, are surveyed, and a comparative numerical study against projection-based methods is presented. The study shows that, across a few benchmark problems, projection methods (for the specific approximations spaces utilized) underperforms preconditioned BiCGSTAB Krylov method in terms of both convergence rate and memory usage. In contrast, Reference ~\cite{simoncini2023analysis} shows that projection-based methods can be computationally superior to the classical conjugate gradient method when an appropriate solution subspace is chosen, and that classical Krylov methods suffer from performance degradation due to loss of orthogonality from the low-rank truncation of the subspaces. Reference~\cite{palitta2021convergence} confirms the loss of orthogonality issues in low-rank Krylov subspace methods, and proposes techniques to mitigate it. They further find that preconditioning is key to avoid substantial subspace growth in matrix-oriented low-rank Krylov methods.}

Preconditioning is indeed key to low-rank Krylov solvers for GSEs. Reference~\cite{palitta2021convergence} proposes an interesting preconditioner for GSEs stemming from variable advection--diffusion equations, which will be of particular relevance to the present study. By averaging the PDE coefficients---a strategy we term here the \emph{average-coefficient-Sylvester} (ACS) approximation---they obtain a constant-coefficient preconditioner of an SE form. The GSE system is then solved using matrix-oriented GMRES (but without low-rank truncation), where the operator application exploits the Kronecker product structure of the GSE, and the preconditioner is applied using KPIK~\cite{simoncini2007new}. In another related work on anisotropic diffusion problems~\cite{meng2024preconditioning}, a basis-update-and-Galerkin (BUG) DLR integrator is employed as a preconditioner for a low-rank matrix-oriented GMRES solver for the associated GSE. The BUG preconditioner evolves the basis and the coefficient matrix to provide effective initial guesses for both the basis vectors and coefficient updates in the GMRES iteration, demonstrating significant potential in reducing the GMRES iteration count and the maximal Krylov rank. Therefore, the performance of (low-rank) matrix-oriented Krylov methods critically depends on the availability of a cheap and suitable preconditioner and an efficient truncation mechanism to limit the iteration count and control subspace growth.

In this study,  we build upon the advances above for SE and GSE solvers within the projection-based framework and develop an adaptive-rank implicit integrator for 2D and 3D time-dependent advection-diffusion equations with variable coefficients. The main innovation of the algorithm lies in the specific treatment of both the basis construction and coefficient-finding steps within a standard projection-based Krylov framework, in which we combine solution elements proposed elsewhere to produce a novel algorithm. 
{Notably, from the projection-based framework, we borrow the idea of enriched atypical Krylov subspaces to ensure an efficient representation of the solution bases. From classical matrix-oriented Krylov methods, we adopt the idea of truncating the subspaces along the iterations to keep it as compact as possible, and the idea of preconditioning the GSE with a simpler equation that is spectrally close. We assemble these ingredients in the following way:} For the basis-building step (which we will refer to as {\em the outer iteration}), we introduce a new Krylov subspace construction strategy that combines (i) inverted Krylov subspaces for each differentiation coefficient matrix, (ii) standard Krylov subspaces for diagonal matrices arising from separability of advection-diffusion coefficients, and (iii) extended Krylov subspaces for averaged-coefficient matrices formed from the ACS approximation. We control excessive growth of the rank of the so-constructed Krylov {subspaces} bases by a targeted truncation procedure using SVD-based QR and modified Gram-Schmidt orthogonalization. In the coefficient-finding step (which we refer to as {\em the inner iteration}), we follow a standard Galerkin projection  and solve the associated system with GMRES, preconditioned by the SE resulting from the Galerkin projection of the ACS approximation. This preconditioner is expected to perform optimally because it is compact equivalent to the original system~\cite{axelsson2007mesh,axelsson2009equivalent} and therefore spectrally close. The key distinction of this approach from others in the literature is applying the {ACS SE} preconditioner to the reduced-system of the projection-based framework with $\mathcal{O}(r^{d+1})$, rather than to the full-scale system (as was done for instance in the matrix-oriented GMRES framework proposed in \cite{palitta2016matrix}). 

These algorithmic innovations enable the method to achieve  linear scaling with the one-dimensional resolution $N$, both for computational complexity [as $\mathcal{O}(N r^2) + \mathcal{O}(r^{d+1})$] and for memory storage [as $\mathcal{O}(N r )+\mathcal{O}( r^{d})$], and allow us to perform 3D simulations of advection-diffusion problems with equivalent resolutions of $\mathcal{O}\left(10^{12}\right)$ degrees of freedom on a standard laptop, orders of magnitude faster than an optimally performing multigrid method scaling as $N^d \log(N)$ \cite{jones1997parallel}. Additionally, as in \cite{el2024krylov}, the adaptive low-rank integrator is amenable for use with DIRK methods for high-order temporal accuracy (up to third-order is considered in this study). All these advertised properties will be demonstrated numerically in the test examples.
}
% \HH{Reference on stationary iterations for solvingch GSE is missing. Find appropriate space to include. }\QQ{insert reference for multigrid scaling.}\lc{Done}
% \QQ{do we want to merge this paragraph with the previous one? The line on "$(3\times10^4)^3$" could be in response, not necessarily be here?}

% \HH{For the group: we have deleted reference on fixed-point iterations from previous submission, though not tightly connected our work, should we reinclude refence? Further references would include Riemanian optimization and greedy low-rank algorithm, both are new and old works by Kressner, respectively.} \WT{My preference is to keep it out, inlight of the already long introduction, we may want to focus on critically relevant references to our work.} \lc{This is a published approach to solve GSEs, and we should not ignore it} \QQ{we could find a line earlier in the literature review, and add those references without commenting on specifics?}

The rest of this paper is organized as follows. In Section~\ref{Section_2}, we consider the case of rank-one advection and diffusion coefficients and derive the GSE that results from a finite-difference discretization in space and backward Euler in time. In Section~\ref{Section_3}, we introduce an adaptive-rank solution strategy for 2D problems, including the coefficient-averaged approximation operator, the construction of {extended Krylov} subspaces (which we abbreviate as xKrylov in the sequel), the preconditioning technique in solving the reduced system, and the low-rank computation of the residual norm. The extension of the technique to high-order temporal discretizations with multi-rank coefficients is provided in Appendix \ref{Appendix_DIRK}, and for 3D problems in Appendix~\ref{Appendix_3D}. 
In Section~\ref{section_4}, numerical results are presented to demonstrate the efficacy of the algorithm. Finally, we conclude and provide an outlook for future research in Section~\ref{Section_5}. 

\section{Problem Description and Generalized Sylvester Equation}\label{Section_2}

% , where each unknown solution matrix is simultaneously operated on across both dimensions

{For simplicity, we focus the discussion in this section on 2D, with the extension to 3D elaborated in Appendix~\ref{Appendix_3D}.}
We establish the following notation conventions. Lowercase Greek letters, such as \(\phi(x, y)\) or \(\sigma(x, y)\), denote continuous scalar fields. Uppercase bold italic Greek letters, such as \(\mathbfit{\Sigma}(x, y)\) or \(\mathbfit{\Phi}(x, y)\), are used for vector-valued and tensor-valued functions. Uppercase Greek letters, like \(\Phi\) or \(\Sigma\), represent discrete matrices arising from the discretization of continuous scalar functions or differential operators on a grid. Boldface uppercase letters, such as \(\mathbf{F}\) or \(\mathbf{S}\), represent unknown matrix solutions or matrix right-hand sides in matrix equations. The identity matrix of size \( n \) is denoted by \( {I}_{n} \).

We aim to solve the \HHTWO{time-dependent} advection-diffusion equation:
% \HH{Is there a specific reason we are generically using $v$ below? \&\& Need to add discussion regarding boundary conditions.}\lc{Not that I know of; we should remove the $v$ subscript} 
%
\begin{equation}
    \frac{\partial f}{\partial t} = \nabla_v \cdot ( \mathbfit{\Phi} \cdot \nabla f) - \nabla \cdot (\mathbfit{\Sigma} f),
\label{eq:adv-diff}
\end{equation}
\HHONE{where $ (x,y) \in \Omega = [a_1, b_1] \times [a_2, b_2]$} and \(\mathbfit{\Phi} \in \mathbb{R}^{2 \times 2}\) is an isotropic diffusion tensor:
\begin{equation*}
\mathbfit{\Phi}(x,y) = \begin{bmatrix}
\phi^{x}(x,y) & 0 \\
0 & \phi^{y}(x,y)
\end{bmatrix},
\end{equation*}
and \(\mathbfit{\Sigma}\in \mathbb{R}^{2}\) is the advection vector, given by:
\begin{equation*}
\mathbfit{\Sigma}(x,y) = \begin{bmatrix}
\sigma^x(x,y) \\
\sigma^y(x,y)
\end{bmatrix}.
% \label{eq: adv_coeffs}
\end{equation*}
We assume that the diffusion and advection coefficients are separable with finite rank. Specifically, the diagonal diffusion coefficients, \(\phi^{x}(x,y)\) and \(\phi^{y}(x,y)\), are given by:
\begin{equation}
\phi^{x}(x,y) = \sum_{i=1}^{\ell_x} \phi^{1,x}_i(x) \phi^{2,x}_i(y), \quad \phi^{y}(x,y) = \sum_{i=1}^{\ell_y} \phi^{1,y}_i(x) \phi^{2,y}_i(y).
\label{eq: diff_coeffs}
\end{equation}
Similarly, for the advection coefficients we have:
\begin{equation}
\sigma^x(x,y) = \sum_{i=1}^{k_x} \sigma^{1,x}_i(x) \sigma^{2,x}_i(y), \quad \sigma^y(x,y) = \sum_{i=1}^{k_y} \sigma^{1,y}_i(x) \sigma^{2,y}_i(y).
\label{eq: adv_coeffs}
\end{equation}
Substituting these coefficients in their separable forms into Equation \eqref{eq:adv-diff}, we obtain:
\begin{equation}
    \begin{aligned}
        \frac{\partial f}{\partial t} &= \sum_{i=1}^{\ell_x} \frac{\partial}{\partial x} \left( \phi^{1,x}_i(x) \frac{\partial f(x,y)}{\partial x} \right) \phi^{2,x}_i(y)+ \sum_{i=1}^{\ell_y} \phi^{1,y}_i(x) \frac{\partial}{\partial y} \left( \phi^{2,y}_i(y) \frac{\partial f(x,y)}{\partial y} \right) \\
        &- \sum_{i=1}^{k_x} \frac{\partial}{\partial x} \left( \sigma^{1,x}_i(x) f(x,y) \right) \sigma^{2,x}_i(y) - \sum_{i=1}^{k_y} \sigma^{1,y}_i(x) \frac{\partial}{\partial y} \left( \sigma^{2,y}_i(y) f(x,y) \right),
    \end{aligned}
\label{eq:adv-diff_sub-multirank}
\end{equation}
where, in each term, the factor independent of the differentiation variable has been pulled out of the derivative.
%
%\subsection{Generalized Sylvester Equation from Implicit Discretization of Advection-Diffusion\label{Gen_sylv_sec}}
To simplify the {presentation}, 
we assume for now {rank-one advection and diffusion coefficients}:
\[
\phi^{x} = \phi^{1,x}(x) \phi^{2,x}(y), \quad \phi^{y} = \phi^{1,y}(x) \phi^{2,y}(y),%^\top, 
\]
\[
\sigma^x = \sigma^{1,x}(x) \sigma^{2,x}(y), \quad \sigma^y = \sigma^{1,y}(x) \sigma^{2,y}(y).
\]
This will be generalized to multiple ranks below. Equation \eqref{eq:adv-diff_sub-multirank} then simplifies to: 
% \lc{I would remove ``diffusion'' and ``advection'' in the equation below}
\begin{equation}
\small
    \begin{aligned}
        \frac{\partial f}{\partial t} &={ \frac{\partial}{\partial x} \left( \phi^{1,x}(x) \frac{\partial f(x,y)}{\partial x} \right) \phi^{2,x}(y) 
        + \phi^{1,y}(x) \frac{\partial}{\partial y} \left( \phi^{2,y}(y) \frac{\partial f(x,y)}{\partial y} \right)} \\
        &{-}{\frac{\partial}{\partial x} \left( \sigma^{1,x}(x) f(x,y) \right) \sigma^{2,x}(y) 
        {-} \sigma^{1,y}(x) \frac{\partial}{\partial y} \left( \sigma^{2,y}(y) f(x,y) \right)}.
    \end{aligned}
\label{eq:adv-diff_sub}
\end{equation}
We discretize this equation in a two-dimensional tensor-product spatial grid with ${N}_1$ and ${N}_2$ grid points in the \(x\) and \(y\) directions, respectively. The grid points in the \(x\) direction are denoted as $x_i$ for $i = 0, \dots, N_1 - 1$, and similarly for the \(y\) direction. We consider the point-wise multiplication of the advection and diffusion coefficients in the discrete by setting: 
\begin{align*}
\Phi_1 &= \text{diag}(\phi^{1,y}(x_i)), 
\Sigma_1 = \text{diag}(\sigma^{1,y}(x_i)),
\quad i = 0, \cdots, N_1 -1 , \\
\Phi_2 &= \text{diag}(\phi^{2,x}(y_j)), \Sigma_2 = \text{diag}(\sigma^{2,x}(y_j)),\quad j= 0, \cdots, N_2-1. 
\end{align*}
{We approximate the spatial derivative with a second-order finite-difference scheme using a three-point stencil, and obtain} the following matrix differential equation for $\mathbf{F}$, for which its $(i, j)$ entry $F_{i,j}$ corresponds to the approximation of $f(x_i, y_j)$:
% \lc{Also remove here}
\begin{equation}
    \frac{\partial \mathbf{F}}{\partial t} =     {D_{xx}^{\Phi}\mathbf{F} \Phi_2^\top + \Phi_1 \mathbf{F} {D_{yy}^{\Phi}}^\top} 
    - {D_x^{\Sigma} \mathbf{F} \Sigma_2^\top - \Sigma_1 \mathbf{F} {D_y^{\Sigma}}^\top}.
    \label{MOL_eq}
\end{equation}
Here, \(D_{xx}^{\Phi}, D_{yy}^{\Phi}, D_x^{\Sigma}, D_y^{\Sigma}\) represent compositions of difference operators with pointwise multiplication operators. See Appendices \ref{sec:discret_diff} and \ref{sec:discret_adv} for more details.
For notational simplicity, let:
\begin{equation}
T_1 = D_{xx}^{\Phi}, \quad T_2 = D_{yy}^{\Phi}, \quad T_3 = -D_x^{\Sigma} , \quad T_4 = -D_y^{\Sigma}.
\label{operator_T}
\end{equation}
Equation \eqref{MOL_eq} can be written compactly as:
\begin{equation}
    \frac{\partial \mathbf{F}}{\partial t} = \mathscr{L}(\mathbf{F}),
    \label{MOL_operator_L}
\end{equation}
where the right-hand side induces the following mapping:
\begin{equation}
{\mathscr{L}} : \mathbf{F} \mapsto
T_1 \mathbf{F} \Phi_2^\top + \Phi_1 \mathbf{F} T_2^\top + T_3 \mathbf{F} \Sigma_2^\top + \Sigma_1 \mathbf{F} T_4^\top. 
\label{operator_L}
\end{equation}
We consider first a backward Euler implicit treatment of the temporal component of the matrix differential equation \eqref{MOL_operator_L}, which gives the following discrete GSE:
\begin{equation}
    \mathbf{F}_1 - \Delta t  \mathscr{L}(\bF_1) = \mathbf{F}_0, 
    \label{Generalized_Sylv}
\end{equation}
and associated residual equation:
\begin{equation}
\mathbf{R_{\mathscr{L}}} = \mathbf{F}_1 - \Delta t \, \mathscr{L}(\mathbf{F}_1) - \mathbf{F}_0.
%= \mathbf{0}.
    \label{Residual_L}
\end{equation}
%\sout{An iterative solver updates the solution \(\mathbf{F}_1\) until the residual is sufficiently small under a suitable norm.}
{Iterative solvers update the solution \(\mathbf{F}_1\) until the residual is sufficiently small under a suitable norm. In this work, we opt for Galerkin-projection-based  iterative solvers, where the bases spanning \(\mathbf{F}_1\) are updated by augmenting Krylov-like subspaces, as shown next.}

\section{{Krylov-Based Adaptive-Rank Solver for Generalized Sylvester Equation
}}\label{Section_3}
%{LC: this paragraph is repetitive}In what follows, we refer to equations in the form of Equation \eqref{eq: Se} as the SE, where the unknown solution matrices are operated on one dimension at a time. Additionally, we use the term GSE to describe equations in the form of Equation \eqref{generalized_eq}. 
%Numerous strategies have been proposed in the literature to solve such equations, most notably the recent work by Palitta and Simoncini in \cite{palitta2016matrix}. {Inspired by their preconditioning technique \cite{palitta2016matrix}, we further exploit the low-rank structure of the solution and design a projection-based algorithm with an optimal scaling of $\mathcal{O}(N r^2)$, i.e. a linear scaling with respect to the resolution $N$ per spatial dimension, and quadratic scaling with respect to the rank $r$.} 

We describe next the projection-based framework for solving the GSE \eqref{Generalized_Sylv} \cite{benner2013low,hao2021sherman,simoncini2023analysis,powell2017efficient}.
%We propose an adaptive-rank algorithm for the GSE \eqref{Generalized_Sylv}.
We assume a low-rank factorization of the initial condition, \(\mathbf{F}_0 \in \mathbb{R}^{N_1 \times N_2}\), which evolves to an updated solution \(\mathbf{F}_1 \in \mathbb{R}^{N_1 \times N_2}\) at time \(t^{(1)} = t^{(0)} + \Delta t\) also in the low-rank format, i.e.:%\lc{Should these $S_0$ and $S_1$ be boldface?}\HH{According to description in introduction: 'Boldface uppercase letters, such as \(\mathbf{F}\) or \(\mathbf{S}\), represent unknown matrix solutions or matrix right-hand sides in matrix equations'. Since this is not a matrix equations, does the rule apply here ?}
\begin{equation}
    \mathbf{F}_0 = U_0 S_0 V_0^\top \xrightarrow{\text{evolve in time}} \mathbf{F}_1 = U_1 S_1 V_1^\top,
\end{equation}
where \(U_0 \in \mathbb{R}^{N_1 \times r_0}\), \(V_0 \in \mathbb{R}^{N_2 \times r_0}\), \(U_1 \in \mathbb{R}^{N_1 \times r_1}\), and \(V_1 \in \mathbb{R}^{N_2 \times r_1}\) are orthonormal bases for their respective dimensions, and \(S_0 \in \mathbb{R}^{r_0 \times r_0}\) and \(S_1 \in \mathbb{R}^{r_1 \times r_1}\) are diagonal matrices with singular values at the corresponding times. \HHTWO{ Projection methods} involve constructing bases \(U_1\) and \(V_1\) and then projecting the original equation onto them to determine the coefficients of \(S_1\). 
%A key aspect of the proposed adaptive-rank approach is that we never form the full solutions explicitly. Instead, we seek bases in a dimension-by-dimension manner and then use them to perform a Galerkin projection.
% \sout{, resulting in a reduced system for the matrix of coefficients \(S_1\). An SVD truncation procedure is then applied to update the orthonormal basis sets and realize an efficient low-rank representation of the updated solution}. 
%Key questions for solving the GSE that we address in this study include:
%\begin{itemize}
%    \item How do we choose {suitable} subspaces effectively to build \(U_1\) and \(V_1\) from \(U_0\) and \(V_0\)?
%    \item How can we solve the reduced equation for \(S_1\) efficiently?
%    \item How can we check the residual of the GSE without forming the matrices explicitly?
%\end{itemize}
Efficiently solving the GSE requires choosing suitable solution subspaces effectively to build \(U_1\) and \(V_1\) from \(U_0\) and \(V_0\), and solving the reduced equation for \(S_1\) efficiently. 
To accomplish both, in the spirit of the work in \cite{palitta2016matrix} we leverage an approximate operator $\Tilde{\mathscr{L}}$ (which is of the Sylvester form) to the original operator ${\mathscr{L}}$ in \eqref{operator_L}. 
% \sout{Our main innovation leverages an approximate operator $\Tilde{\mathscr{L}}$ (which is of the Sylvester form) to the original operator ${\mathscr{L}}$ in \eqref{operator_L}. 
% \sout{$\Tilde{\mathscr{L}}$ approximates the diagonal matrices \(\Phi_i\) and \(\Sigma_i\) representing pointwise multiplication of diffusion and advection coefficients, respectively, by scaled identities.}
% This approximation not only allows us to construct dimension-wise Krylov bases $U_1$ and $V_1$ effectively and efficiently, but also enables effective preconditioning for solving the projected reduced system for \(S_1\), in the spirit of \cite{palitta2016matrix}.} %Then we project Equation \eqref{Generalized_Sylv} onto the Krylov bases, resulting in a reduced system for the coefficient matrix $S_1$. 
In the following, 
% \sout{subsections, we introduce the main steps of our proposed algorithm.} 
% as illustrated in Figure~\ref{flow-chart}, highlighting the intuition and innovations behind our design. \lc{This figure is only referred to once; do we really need it?} 
we describe the construction of the approximated operator $\tilde{\mathscr{L}}$ in Section \ref{ssec: L_tilde}; we discuss the dimension-wise basis construction in Section \ref{basis_cons};  we present the reduced equation for the matrix of coefficients $S_1$ and introduce our particular implementation of the ACS preconditioner in Section \ref{preconditioning}, and outline the %\sout{solution truncation mechanism and the} \WT{haven't mentioned anything about truncation upto this point} 
residual computation in Section \ref{residual}. {Finally, Section \ref{multi_rank_BE} extends the algorithm to the case of multi-rank coefficients.}

% \QQ{revise to shorten...}
%. \WT{Seems we've removed the flowchart so these sentences should be modified accordingly.}\HH{Updated accordingly.}

% \begin{figure}[t]
%     \centering    \includegraphics[width=0.9\textwidth]{Figure1.png}
%     \caption{Flow-chart of the proposed implicit adaptive-rank algorithm for a GSE. 
%     }
%     \label{flow-chart}
% \end{figure}

% \subsection{{Averaged-Coefficient} Approximate Operator $\tilde{\mathscr{L}}$}
% \label{ssec: L_tilde}
\subsection{{Averaged-Coefficient} Approximate Operator $\tilde{\mathscr{L}}$}
\label{ssec: L_tilde}
We incorporate ideas from \cite{palitta2016matrix} and \cite{el2024krylov} to construct an approximate operator $\Tilde{\mathscr{L}}$ and efficiently solve the time-dependent inhomogeneous advection-diffusion equation using a projection-based adaptive-rank integration framework. {The approximated operator $\tilde{\mathscr{L}}$ serves two purposes: (i) contribute to constructing the solution bases $U_1$ and $V_1$ from xKrylov subspaces, and (ii) preconditioning the solution of the coefficient matrix $S_1$.} 
% \sout{We first construct a basis informed by both the original operator $\mathscr{L}$ and the approximate operator $\Tilde{\mathscr{L}}$ \sout{(defined below)}; then we derive a corresponding projected (reduced) system that we solve iteratively using GMRES with $\Tilde{\mathscr{L}}$ as preconditioner}.
 
We seek an approximation to ${\mathscr{L}}$ that captures its spectral content and is easier to invert. As proposed in \cite{palitta2016matrix}, the {advection} diagonal operators can be approximated by scaling an identity matrix, where the scaling factor is the average of the diagonal entries. \HHTWO{In this work, we approximate not only the advection coefficients by a suitable constant multiple of the identity,
% \sout{, as performed in \cite{palitta2016matrix}}
but extend this strategy to the scalar diffusion coefficients by averaging them across the computational domain.}
Let \(\Phi_i \approx {\alpha}_i I\) and \(\Sigma_i \approx \gamma_i I\), where \(\alpha_i\) and \(\gamma_i\) are the averages of the diagonal entries of \(\Phi_i\) and \(\Sigma_i\), respectively. This approximation yields the following mapping: 
\begin{equation}
{\Tilde{\mathscr{L}}} : \mathbf{F} \mapsto
\left(\alpha_2 T_1 + \gamma_2 T_3\right) \mathbf{F}
+ \mathbf{F} \left(\alpha_1 T_2 + \gamma_1 T_4\right)^\top.
\label{approx_operator_L}
\end{equation}
A backward Euler temporal discretization using \eqref{approx_operator_L} leads to the SE:
\begin{equation}
    (A_1+A_3) \mathbf{F}_1 + \mathbf{F}_1 (A_2+A_4)^\top = \mathbf{F}_0,
    \label{Approximated_Sylv}
\end{equation}
with:
\begin{align}
     A_1 &= \left(\frac{1}{4}I - \Delta t \, \alpha_2 T_1\right) \quad A_3=\left(\frac{1}{4}I - \Delta t \, \gamma_2 T_3\right), \quad   P_1 = A_1+A_3, \label{eq:A1-A3}\\
      A_2 &= \left(\frac{1}{4}I - \Delta t \, \alpha_1 T_2\right) \quad A_4=\left(\frac{1}{4}I - \Delta t \, \gamma_1 T_4\right), \quad P_2 = A_2 + A_4   \label{eq:A2-A4}.
\end{align}
Equation \eqref{Approximated_Sylv} is an SE approximation to the GSE \eqref{Generalized_Sylv}, which can be efficiently solved using low-rank approximations~\cite{el2024krylov, simoncini2007new, shank2013krylov, saad1989numerical}.
It can be exploited as a preconditioner for a GMRES iterative solution of the GSE \eqref{Generalized_Sylv}. 
{This was done in \cite{palitta2016matrix} with {$\mathcal{O}(N^3)$ }computational complexity.}
%\HH{Authors may wonder why not just use Simoncini's approach? I think this line puts forward a motivation for our approach.}  \lc{I don't disagree with this.}
In this work, we extend this approach by performing a Galerkin projection of both the original operator $\mathscr{L}$ and its approximation $\tilde{\mathscr{L}}$ on the low-rank bases before solving the GSE with GMRES. This reduces the problem to a lower-dimensional subspace, where preconditioning can be applied more effectively on a smaller $d$-dimensional problem with $\mathcal{O}( r^{d+1})$ complexity. 

% \lc{This next sentence is repetitive; I would remove}\HHONE{Our innovation using the approximate operator $\Tilde{\mathscr{L}}$ lies in the effective construction of dimension-wise xKrylov subspaces detailed next in Section~\ref{basis_cons}, as well as preconditioning the reduced projected system elaborated in Section~\ref{preconditioning}. }
% \QQ{I suggest to only keep (ii) part.}\HH{Done.}
%but also for which are generated from multiple operators present in ${\mathscr{L}}$, $\tilde{\mathscr{L}}$ and their inverses. 
% This construction is detailed next in Section~\ref{basis_cons}; (ii) Rather than preconditioning the operator \(\mathscr{L}\) directly at the matrix level, we precondition its Galerkin projection onto the carefully selected xKrylov subspaces, as will be shown in Section \ref{preconditioning}.

\subsection{Basis Construction}\label{basis_cons}
% \QQ{In the algorithm below, (1) introduce algorithm 1 first; (2) change the notations a bit: the superscript for $V$ for the $m$ Krylov iteration; and subscript for the operators? $V^{m}$, vs. $V_{(l)}$. and introduce the convection of notation we set up... (3) why the last line of the algorithm... (4)organize the procedure as an algorithm, as in Algorithm 1}\HH{Done.}

The construction of the the xKrylov subspaces per dimension at the \(m\)-th outer iteration is performed as:
\begin{equation}
\label{eq:K_N1}
\kappa_m(M_1, \cdots, M_l, U) = \text{\HHTWO{Range}}\Big\{U, \underbrace{M_1 U, \cdots, M_l U}_{\text{1st augmentation}}, \ldots, \underbrace{M_1^{m+1} U, \cdots, M_l^{m+1} U}_{\text{mth augmentation}}\Big\},
\end{equation}
where \( M_1, \cdots, M_l \) are operators of interest {and their inverses} arising in the GSE, {and $l$ is the number of operators used for each basis augmentation step.} 
% Depending on the structure of the matrix equation~\eqref{generalized_eq}, \( l \) may be large, and the generated subspace vectors at each augmentation may exhibit linear dependencies. 
%The xKrylov basis construction involves several operators from the original GSE \eqref{Generalized_Sylv} and its approximation \eqref{Approximated_Sylv}. 
%Following \cite{el2024krylov}, we construct an xKrylov subspace in the $x$-direction 
For the $M_l$ operators, we include \(P_1\) and its inverse from the SE approximation to the GSE \eqref{Approximated_Sylv}, the individual operators \(A_1\) and \(A_3\) and their inverses, as they may contain information not captured by the inverse of \(P_1\), and the diagonal operators \(\Phi_1\), \(\Sigma_1\), \(\Phi_2\), and \(\Sigma_2\), which are averaged to obtain $\tilde{\mathscr{L}}$.  Using these operators, the generalized Krylov {subspace}
%{LC: repetitive}Additionally, we use the inverse operators from the approximated operators \(A_{1:4}\) in \eqref{eq:A1-A3}-\eqref{eq:A2-A4} to form an inverted Krylov {subspace}. Finally, we construct dimension-wise extended Krylov subspaces for the approximated operator \(\tilde{\mathscr{L}}\) using \(P_1\) and \(P_2\). 
%{LC: not needed}The reasoning behind these choices of operators for subspace construction is as follows. 
% \sout{This carefully chosen set of operators results in a compact representation of the solution matrix, enriching the basis while keeping its dimension small. Combining these, we define the xKrylov subspaces}
(for instance in the \(x\)-direction) considered is
as follows:
% \HHONE{
\[
\begin{aligned}
    \kappa_m\left(M_1, \cdots, M_6, U_0\right), \quad \text{where} \quad 
    & M_1 = P_1,
    & M_2 = P_1^{-1},  \quad 
    & M_3 = A_1^{-1},   
    & M_4 = A_3^{-1},  \quad
    & M_5 = \Phi_1,  
    & M_6 = \Sigma_1. 
\end{aligned}
\]
We have found that including all these matrices is 
% \sout{crucial for accelerating convergence}
{important to capture the range of the solution}.

Depending on the structure of the GSE~\eqref{generalized_eq}, \( l \) in \eqref{eq:K_N1} may be large, and the generated subspace vectors at each augmentation may exhibit linear dependencies. To limit subspace inflation, we employ an SVD-based truncation strategy (Algorithm~\ref{alg:basis_trunc}), such that only basis vectors corresponding to significant singular values are retained.
%Drawing inspiration from \cite{simoncini2007new},  and using the truncation procedure outlined in Algorithm \ref{alg:basis_trunc}, we provide next an algorithm for 
%
With this in mind, the construction of {the dimension-wise xKrylov subspaces} \eqref{eq:K_N1}, following \cite{simoncini2007new}, is as follows. 
\HHONE{At each iteration $m$, Algorithm~\ref{alg:basis_constr} augments the previously constructed basis $\mathcal{V}_{(m-1)}$ by applying powers of operators $M_1$ through $M_l$ to their corresponding basis blocks $V^{(1)}, \dots, V^{(l)}$. Each block matrix $V^{(i)}$ may have a rank of at most $r_0$. We apply Algorithm \ref{alg:basis_constr} to orthogonalize the {subspace and obtain a basis $\mathcal{U}_{(m)} \in \mathbb{R}^{N_1 \times r^{(m)}_x}$} and $\mathcal{V}_{(m)}\in \mathbb{R}^{N_2 \times r^{(m)}_y}$ in the \(x\) and \(y\) directions, respectively. 
{For notational simplicity, we omit the calligraphic notation and subscript \((m)\) from now on, understanding that \(U_1\in \mathbb{R}^{N_1 \times r_x} \) and \(V_1\in \mathbb{R}^{N_2 \times r_y}\) refer to the bases of the Krylov subspaces at the \(m\)th iteration.} 
% \WT{Another potential notation inconsistency. Why are we bold mathbb'ing the real numbers here?}
The truncation step in line~\ref{ind_trunc_line_1} of Algorithm~\ref{alg:basis_constr}  removes any linear dependencies {within each block} during the power iterations; line~\ref{ind_trunc_line_2} applies a global truncation across all newly generated blocks, eliminating redundancies among them. As a result, the {subspace} is kept optimal in {size}.} 

% The same procedures are applied to the \(y\)-direction {to obtain }.

%\QQ{Luis: we revised a bit to explain line 4 and 6. please take a look...}
% \QQ{the notation of  $V_{(m-1)}$ and  $V^{(l)}$ are a bit confusing...how about $\mathcal{V}_{(m-1)}$?}\HH{Updated}

\begin{algorithm}[t]
\caption{\label{alg:basis_trunc} {SVD-truncated} QR, $\mathscr{T}_{\epsilon_{\kappa}}$}
\SetAlgoNlRelativeSize{-2}
\SetNlSty{textbf}{(}{)}
\KwIn{Krylov Subspace $\kappa_m$; Tolerance $\epsilon_{\kappa}$.}
\KwOut{Truncated bases ${Q}$.}
Compute the reduced QR decomposition:
    $ \{ Q, R \} = \texttt{QR}\left(\kappa_m\right)$\label{alg1:line:basis}\;
Perform reduced SVD decomposition: $Z_1 \Tilde{R} Z_2 = \texttt{SVD}(R)$, where $\Tilde{R} = \text{diag}(\sigma_j)$\label{alg1:line:svd:basis} \;
% \WT{It should be clear from context, but we've already used $T_1$ and $T_2$ to denote discrete diffusion operators.}\HH{Fixed.}\;
Identify the last index $r$ such that $\sigma_{r} > \epsilon_{\kappa}$ \tcp*[r]{Optional: Set Maximum Allowable Rank}
Update the basis vectors: $Q \leftarrow Q Z_1(:, 1:r)$\label{alg1:line:basis:update}\;
\end{algorithm}

\begin{algorithm}[t]
\caption{\label{alg:basis_constr} Construction of xKrylov Subspace $\kappa_m^{\texttt{trunc}}$ %\lc{Why are we truncating twice in (4) and (6) below?}\HH{(1) remove redundancy within each basis block (2) remove redundancy between the blocks.} \lc{I realize that, but wouldn't the latter take care of the former?}
}
\SetAlgoNlRelativeSize{-2}
\SetNlSty{textbf}{(}{)}
\KwIn{Operators $M_1,\dots, M_l$; augmentation basis blocks $\{V^{(1)}, \dots, V^{(l)}\}$; previous Krylov basis ${\mathcal{V}}_{(m-1)}$; tolerance $\epsilon_{\kappa}$.}
\KwOut{Updated Krylov basis $\mathcal{V}_{(m)}$; basis blocks $\{V^{(1)}, \dots, V^{(l)}\}$ for the next iteration.}

\For{$i = 1,\dots,l$}{
    $V^{(i)} \leftarrow M_i V^{(i)}$\;
    $\widetilde{V}^{(i)} \leftarrow \text{orthogonalize } V^{(i)} \text{ w.r.t. } V_{(m-1)}$ (e.g., modified Gram-Schmidt)\;
    $V^{(i)} \leftarrow \mathscr{T}_{\epsilon_{\kappa}}(\widetilde{V}^{(i)})$\label{ind_trunc_line_1} \tcp*[r]{truncate to $\epsilon_\kappa$ for next iteration using Algorithm \ref{alg:basis_trunc}}
}
$V' \leftarrow [\widetilde{V}^{(1)}, \dots, \widetilde{V}^{(l)}]$\;
$V' \leftarrow \mathscr{T}_{\epsilon_{\kappa}}(V')$\label{ind_trunc_line_2} \tcp*[r]{final orthogonalization and truncation using Algorithm \ref{alg:basis_trunc}}
$\mathcal{V}_{(m)} \leftarrow [\mathcal{V}_{(m-1)}, V']$\;
\end{algorithm}

\HHONE{We comment next on the expected performance of the xKrylov iteration.
% \sout{ regarding the dependence of the subspace dimension on grid resolution for a given tolerance}.
In principle, the number of xKrylov iterations may depend on the mesh size, $N$. For instance, \cite{kressner2010krylov,beckermann2013error} theoretically show for the 2D and 3D Poisson equations that the dimension of the extended Krylov subspace (denoted by $r^{(m)}$ in the basis construction above) scales as $\mathcal{O}(\kappa^{1/4})\sim\mathcal{O}(N^{1/2})$ in 2D, {and $(\kappa^{2/3})\sim\mathcal{O}(N^{4/3})$ in 3D}, with $\kappa\sim N^2$ the condition number.
%\WT{Perhaps we can just add reference \cite{beckermann2013error} into the previous sentence and say 2D and 3D Poisson equation and remove this sentence?}\HH{Done. Thanks!} 
% \sout{While this provides a worst-case upper bound without considering solution smoothness, such scaling was not observed in practice in the references, and is therefore likely pessimistic.} 
% \sout{\HH{Noting that, in the case of constant advection and diffusion coefficients, the operators $\mathscr{L}$ and $\Tilde{\mathscr{L}}$ coincide, our proposed methodology reduces to a standard extended Krylov projection-type method.}}
% \sout{We empirically observe, therefore, that the theoretical bound remains pessimistic in the context of variable coefficients.
% We numerically verify this claim in Section~\ref{section_4} and find that, for moderately tight tolerances, the required Krylov subspace dimensions remain independent of grid resolution for a wide range of grid sizes.}
{While this provides an upper bound of the convergence rate (without considering solution smoothness), our numerical results suggest that it may be pessimistic in practice, with Krylov subspace dimensions remaining largely independent of grid resolution for moderate convergence tolerances and mesh refinements.}
}
% Hence, if $r$ is independent of $N$, the solution possesses a low-rank structure, and the constructed approximation subspace efficiently maintains moderate $r$, the complexity scales linearly with respect to $N$ cubically with respect to $r$ in 2D, and quartically with respect to $r$ in 3D; which is numerically verified in Section~\ref{section_4}}.\QQ{suggest to remove "Hence, ..."}

% \sout{
% These constructed orthonormal bases will be used to represent the solution in the approximation subspace. The coefficients composing the linear combination of these basis vectors will be obtained by projecting the corresponding GSE \WT{onto the span of such bases} and evolving the matrix of coefficients \( S_1 \).
% }\HH{This sentence provides a good segue to the next one. Consider keeping?} \lc{We have already provided plenty of perspective by this point; this sentence is repetitive and actually does not help with the flow...}
%\QQ{We revise this subsection a bit; take a closer look.}\lc{So did I...}

\subsection{Solving for the Matrix of Coefficients, {$S_1$}: ACS-Preconditioned GMRES}\label{preconditioning}

%\WT{General comment (take it or leave it): We have a focused computational complexity section later on that also talk about specific subsolvers (e.g., Bartels-Steward). I suggest removing discussion on complexity here to compress and avoid distractions.} \lc{I think it is still useful to mention complexity expectations of the PC here; the later 
With the orthonormal bases \(U_1 \in \mathbb{R}^{N_1 \times  r_x}\) and \(V_1 \in \mathbb{R}^{N_2 \times  r_y}\) at hand, we perform a Galerkin projection to derive a reduced GSE for the matrix \(S_1 \in \mathbb{R}^{r_x \times  r_y}\). 
The reduced GSE for \(\bS_1\) arises from projecting the residual via the Galerkin condition \(U_1^\top \mathbf{R}_{\mathscr{L}} V_1 = 0\), yielding:
\begin{align}
\bS_1 - \Delta t \left( \tilde{T}_1 \bS_1 \tilde{\Phi}_2 + \tilde{\Phi}_1 \bS_1 \tilde{T}_2^\top + \tilde{T}_3 \bS_1 \tilde{\Sigma}_2 + \tilde{\Sigma}_1 \bS_1 \tilde{T}_4^\top \right) = \tilde{\bB}_1, \qquad \text{with} \qquad 
\tilde{\bB}_1 = (U_1^\top U_0) S_0 (V_1^\top V_0)^\top.
\label{Galerkin_S}
\end{align}
This matrix-based equation is equivalent to the classical Kronecker formulation:
\begin{equation}
\underbrace{\left[ I_{ r_x r_y \times  r_x r_y} - \Delta t \left( \tilde{\Phi}_2 \otimes \tilde{T}_1 + \tilde{T}_2 \otimes \tilde{\Phi}_1 + \tilde{\Sigma}_2 \otimes \tilde{T}_3 + \tilde{T}_4 \otimes \tilde{\Sigma}_1 \right) \right]}_{\mathcal{A}} \text{vec}(\bS_1) = \text{vec}(\tilde{\bB}_1).
\label{Galerkin_S_kron}
\end{equation}
The projected operators are defined as follows: 
\begin{align*}
\tilde{T}_1 &= U_1^\top T_1 U_1, \quad
\tilde{T}_3 = U_1^\top T_3 U_1, \quad
\tilde{\Phi}_1 = U_1^\top \Phi_1 U_1, \quad
\tilde{\Sigma}_1 = U_1^\top \Sigma_1 U_1, \\
\tilde{T}_2 &= V_1^\top T_2 V_1, \quad
\tilde{T}_4 = V_1^\top T_4 V_1, \quad
\tilde{\Phi}_2 = V_1^\top \Phi_2 V_1, \quad
\tilde{\Sigma}_2 = V_1^\top \Sigma_2 V_1.
\end{align*}
Solving the Kronecker formulation \( \mathcal{A} \text{vec}({\bS_1}) = \text{vec}({\tilde{\bB}_1}) \) in \eqref{Galerkin_S_kron} with a direct solver scales as \(\mathcal{O}( r^{3d})\), {with $r$ the maximum rank of all dimensions and $d$} the dimensionality, which may be too expensive even for moderate \( r\)-values. Instead, we solve this equation iteratively using preconditioned GMRES. We wrap standard GMRES around the Kronecker form of the GSE, but exploiting the matrix-equation form \eqref{Galerkin_S} for both the $\mathcal{A}$-times-vector operation and the preconditioning step as indicated below. The resulting overall complexity is $\mathcal{O}(r^{d+1})$. 
% \sout{This is possible due to the Kronecker structure of our operator, yielding the equivalence between \eqref{Galerkin_S} and \eqref{Galerkin_S_kron}.}
For the GMRES operator-vector product, given a vector {\(s_k = \text{vec}(\mathbf{S}_k)\)}, the product of the matrix \(\mathcal{A}\) times a vector $s_k$ is performed as follows: 
\begin{enumerate}
    \item Given \(s_k\), reshape to \(\mathbf{S}_k \in \mathbb{R}^{ r_x \times  r_y}\), such that \(\text{vec}(\mathbf{S}_k) = s_k\).
    \item Compute \(\hat{\mathbf{Z}}_k \leftarrow \mathbf{S}_k - \Delta t \left( \tilde{T}_1 \mathbf{S}_k \tilde{\Phi}_2 + \tilde{\Phi}_1 \mathbf{S}_k \tilde{T}_2^\top + \tilde{T}_3 \mathbf{S}_k \tilde{\Sigma}_2 + \tilde{\Sigma}_1 \mathbf{S}_k \tilde{T}_4^\top \right)\){, with complexity \(\mathcal{O}( r^{d+1})\)}.
    \item Reshape \(\hat{z}_k = \text{vec}(\hat{\mathbf{Z}}_k)\).
\end{enumerate}
For efficiency, GMRES requires preconditioning
%To bound the number of GMRES iterations with respect to the basis dimension (rank) \( r\), we require effective preconditioning such 
so that the eigenvalues of the preconditioned matrix are sufficiently clustered to deliver $r$-independent convergence rates ~\cite{campbell1996gmres}. 
Here, we employ the ACS preconditioner, with 
$\mathcal{P}= I_{r_y}\otimes \tilde{P}_1 + \tilde{P}_2\otimes I_{r_x}$ arising from the Galerkin projection of \eqref{Approximated_Sylv}{, with $\tilde{P}_1 =U_1^\top P_1 U_1$ and $\tilde{P}_2 = V_1^\top P_2 V_1$.} \HHONE{The ACS preconditioner delivers optimal, rank-independent, superlinear convergence (all observed in our numerical experiments), owing to its connection to the theory of compact-equivalent operators for diffusion-like equations.\footnote{The averaged-coefficient advection-diffusion equation is compact-equivalent to the variable-coefficient one,  which are known to be effective preconditioners for Krylov methods  \cite{axelsson2007mesh,axelsson2009equivalent}.}}
In particular, we consider left preconditioning as: 
\begin{equation}
\mathcal{P}^{-1} \mathcal{A} \text{vec}(\bS_1) = \mathcal{P}^{-1} \text{vec}(\tilde{\bB}_1). 
\end{equation}
% \sout{\HHTWO{We note that, a right preconditioner could similarly be applied.}}\sout{, but the left preconditioner demonstrated robust and consistent performance across all of our tests.} \QQ{only comment on right preconditioner in the response?} \lc{OK}\HH{Isn't this important to address reviewer's concern about left/right preconditioning?} \lc{We do in the response, but it is not really a crucial detail.}\HH{Noted.}
% \sout{For fast GMRES convergence, $\mathcal{P}$ should approximate the left-hand side of Equation \eqref{Galerkin_S_kron} well while being easier to invert than the original system. 
% Here, we consider} 
%In particular,  {let $U_1$ and $V_1$ be the orthonormal bases for the xKrylov, and let $\tilde{P}_1 =U_1^\top P_1 U_1$ and $\tilde{P}_2 =V_1^\top P_2 V_1$, then} 
%\begin{equation}
%  \mathcal{P} : \mathbf{S} \mapsto {\Tilde{P}_1} \mathbf{S} + \mathbf{S} {\Tilde{P}_2^\top}.
%  \label{eq: P_op}
%\end{equation}
To compute $\mathcal{P}^{-1} \hat{z}_k = z_k$ for an arbitrary vector $\hat{z}_k$ during GMRES,  we solve the corresponding SE: 
\[
\tilde{P}_1 \mathbf{Z}_k + \mathbf{Z}_k \tilde{P}_2^\top = \hat{\mathbf{Z}}_k,
\]
where \(z_k = \text{vec}(\mathbf{Z}_k) \in \mathbb{R}^{ r_x r_y}\). Specific steps are:
\begin{enumerate}
    \item Given \(\hat{z}_k\), reshape to \(\hat{\mathbf{Z}}_k \in \mathbb{R}^{ r_x \times  r_y}\), such that \(\text{vec}(\hat{\mathbf{Z}}_k) = \hat{z}_k\).
    \item Solve \(\mathcal{P}_1 \mathbf{Z}_k + \mathbf{Z}_k \mathcal{P}_2^\top = \hat{\mathbf{Z}}_k\) {with complexity $\mathcal{O}(r^{d+1})$, e.g. with the Bartels-Stewart algorithm~\cite{golub2013matrix}.} %\WT{We repeat this later on.}\lc{We do, but it is important to mention here to justify the paragraph below...}
    \item Reshape \(z_k = \text{vec}(\mathbf{Z}_k)\).
\end{enumerate}
%The proposed ACS preconditioner effectively bounds the number of GMRES iterations required to solve Equation \eqref{Galerkin_S} independently of $ r$, as will be shown later in this study, 
As a result, the $S$-matrix solve scales as \(O( r^{d+1})\) overall, {a great reduction in complexity compared to the} \(O(r^{3d})\) scaling of a na\"ive implementation of GMRES with the Kronecker form. 

\subsection{Truncation and Residual Computation}\label{residual}

% \sout{At the end of each xKrylov iteration, we need to compute the residual to accept or reject the solution, while taking advantage of its low rank structure.} 
% \sout{This process should take advantage of the low-rank representation of the solution while avoiding forming the residual on the full grid.} 

Once the solution \( \bF_1 = U_1 S_1 V_1^\top \) is obtained, we apply truncation to control rank growth and reduce the residual computational load. Since \( S_1 \) resulting from the GMRES solve is not necessarily diagonal, we diagonalize it with SVD \cite{guo2022conservative} [also scaling as \(\mathcal{O}(r^{d+1})\)] and discard singular values below a threshold \( \epsilon \).
The SVD truncation gives \(\bF_1 = U_{1,\epsilon} S_{1,\epsilon} V_{1,\epsilon}^\top\), with \(U_{1,\epsilon} \in \mathbb{R}^{N_1 \times r_1}\), \(V_{1,\epsilon} \in \mathbb{R}^{N_2 \times r_1}\), and \(S_{1,\epsilon} \in \mathbb{R}^{r_1 \times r_1}\), with \(r_1 \leq \min(r_x,r_y)\) (Algorithm \ref{alg:trunc}).

Following \cite{el2024krylov}, we compute the residual using the identity:
\begin{align}
\|\mathbf{R_{{\mathscr{L}}}}\| &= \|\mathbf{F}_1 - \Delta t \, \mathscr{L}(\mathbf{F}_1) - \mathbf{F}_0\|  \nonumber \\
&= \underbrace{\left\| 
R_U \begin{bmatrix}
    - \tilde{B}_1 & \mathbf{0} & \mathbf{0} & \mathbf{0} & \mathbf{0} & \mathbf{0}\\
    \mathbf{0} &  {{S}_{1,\epsilon}} & \mathbf{0} & \mathbf{0} & \mathbf{0}& \mathbf{0} \\
    \mathbf{0} & \mathbf{0} & -\Delta t \, {{S}_{1,\epsilon}} & \mathbf{0} & \mathbf{0} & \mathbf{0}\\
    \mathbf{0} & \mathbf{0} & \mathbf{0} & -\Delta t \, {{S}_{1,\epsilon}} & \mathbf{0} & \mathbf{0}\\
    \mathbf{0} & \mathbf{0} & \mathbf{0} & \mathbf{0} & -\Delta t \, {{S}_{1,\epsilon}} & \mathbf{0}\\
    \mathbf{0} & \mathbf{0} & \mathbf{0} & \mathbf{0} & \mathbf{0} & -\Delta t \, {{S}_{1,\epsilon}}\\
\end{bmatrix}
R_V^\top \right\|}_{\in \mathbb{R}^{(r_x + 5 r_1) \times (r_y + 5 r_1)}} \nonumber
\\ & 
=  {\left\| 
R_U \begin{bmatrix}
  - \tilde{B}_1 & \mathbf{0} & \mathbf{0}\\
    \mathbf{0} &   {S_{1,\epsilon}} &\mathbf{0} \\
    \mathbf{0}&\mathbf{0}&-\Delta t \, (I_{\mathcal{R}_{\text{ranks}}} \otimes {S_{1,\epsilon}}
    )     
\end{bmatrix}
R_V^\top \right\|} \leq \epsilon_{\texttt{tol}},
\label{eq:Residual}
\end{align}
where 
\(\mathcal{R}_{\text{ranks}} = {\ell_x + k_x + \ell_y + k_y}\) is the total rank number {of the advection and diffusion coefficients}, and $R_U$ and $R_V$ are the upper triangular matrices arising from the following {Q-less} ${\texttt{QR}}$  decomposition:
\begin{align}
    % \begin{cases}   {S_{1,\epsilon}}
    \{\_ \ , R_U\} &= \texttt{QR}\left(\begin{bmatrix}
      {U}_{1}&  {U}_{1,\epsilon}  & {T_1} {U}_{1,\epsilon} & \Phi_1 {U}_{1,\epsilon} & {T_3} {U}_{1,\epsilon} & \Sigma_1 {U}_{1,\epsilon}
    \end{bmatrix}\right)     \label{RURV},\\
        \{\_ \ , R_V\} &= \texttt{QR}\left(\begin{bmatrix}
      {V}_{1} &  {V}_{1,\epsilon} &\Phi_2 {V}_{1,\epsilon}  & {T_2} {V}_{1,\epsilon} & \Sigma_2 {V}_{1,\epsilon}  & {T_4} {V}_{1,\epsilon} 
    \end{bmatrix}\right).
    % \end{cases} ,
\nonumber
\end{align}
The above residual result follows from the property that $U_1 (U_1)^\top \bF_0 V_1 (V_1)^\top = \bF_0$, with \(\bF_0 = U_0 S_0 V_0^\top\) (where \(U_0 \in \mathbb{R}^{N \times r_0}_0\) and \(V_0 \in \mathbb{R}^{N \times r_0}_0\) have orthonormal columns), since \(U_1\) and \(V_1\) contain \(U_0\) and \(V_0\), respectively.

\HHREV{We note that, as discussed in \cite{palitta2025subspace} and analyzed in detail later in this study, forming, storing, and evaluating the norm of the residual tensor in~\eqref{eq:Residual} can be a very expensive operation when the number of GSE terms, $\mathcal{R}_\text{ranks}$, becomes large. This is so because the complexity of a na\"ive implementation of the residual evaluation scales as $\mathcal{O}(\mathcal{R}_\text{ranks}^{d+1})$ and the memory storage scales as $\mathcal{O}(\mathcal{R}_\text{ranks}^{d})$ (although the latter has been reduced in our implementation to $\mathcal{O}(\mathcal{R}_\text{ranks}^{2})$ by careful algorithmic choreography). While this is not a pressing issue in our numerical experiments (which feature small-to-moderate  values of $\mathcal{R}_\text{ranks}$), it could become a limiting factor in applications where $\mathcal{R}_\text{ranks}$ becomes large. In such instances, it may be necessary to use alternate strategies as proxy measures of convergence
%randomized algorithms were explored as a way to evaluate the residual, but were found insufficient to bound the error in the solution. Alternatively,
such as the difference between consecutive iterates~\cite{palitta2025subspace}, which is commensurate in cost with the rest of our algorithm. 
%We have tested the latter approach in our implementation and found it reliable (see Sec. \ref{section_4}).  This can be understood from the spectral properties of the advection-diffusion operator. In particular, writing the linear solution in residual-correction form as (using Kronecker form for simplicity):
%\[
%X_1-X_1^k = -\mathcal{A}^{-1} (\mathcal{A} X_1^k - X_0) = -\mathcal{A}^{-1} R_k,
%\]
%it follows that:
%\[
%\delta_k \leq \|\mathcal{A}^{-1}\|\, \|R_k\|.
%\]
%The Frobenius norm of $\mathcal{A}^{-1}$ is dominated by the smallest eigenvalues of $\mathcal{A}$, which for our application are of $\mathcal{O}(1)$ and independent of resolution. Therefore $\delta_k \sim \|R_k\|$ even if $\mathcal{A}$ is ill conditioned (i.e., with a large maximum eigenvalue).  
%A hybrid approach that combines $\delta_k$ and the true residual~\eqref{eq:Residual} such that the residual is evaluated only infrequently may also be employed to ensure reliable convergence monitoring while mitigating the high computational cost of the residual evaluation. 
}

\begin{algorithm}[t]
\caption{\label{alg:trunc}Truncated SVD, $\mathcal{T}_\epsilon$ 
}
\SetAlgoNlRelativeSize{-2}
\SetNlSty{textbf}{(}{)}
\KwIn{Bases $U$, $V$; matrix of coefficients $\bS$; Tolerance $\epsilon$.}
\KwOut{Truncated bases $\Tilde{U}$, $\Tilde{V}$; truncated singular values $\Tilde{S}$.}

Perform reduced SVD decomposition: $T_1 \Tilde{S} T_2 = \texttt{SVD}(S)$, where $\Tilde{S} = \text{diag}(\sigma_j)$\label{alg1:line:svd}\;

Identify the last index $r_1$ such that $\sigma_{r_1}/\|S\| >  \epsilon$\;

Update matrix of singular values: $\Tilde{S} \leftarrow \Tilde{S}(1:r_1, 1:r_1)$\;

Update left singular vectors: $\Tilde{U} \leftarrow U T_1(:, 1:r_1)$\label{alg1:line:mm1}\;

Update the right singular vectors: $\Tilde{V} \leftarrow V T_2(:, 1:r_1)$\label{alg1:line:mm2}\;

\end{algorithm}

\subsection{{Generalization to Multi-rank Coefficients}}\label{multi_rank_BE}

We generalize next the adaptive-rank approach to the multi-rank coefficient case \HHTWO{in \eqref{eq: diff_coeffs} and \eqref{eq: adv_coeffs}}, using backward Euler for now.
% \sout{High-order temporal discretizations are considered later in this study \QQ{in Appendix ...?}.}
Discretization of the advection-diffusion equation in this case yields a multi-term GSE of the form:
\begin{equation}
    \mathbf{F}_1 - \Delta t \left( 
    \underbrace{ 
    \sum_{i=1}^{\ell_x} T_{1,i} \mathbf{F}_1 {\Phi_i^{2,x}}^\top 
    + \sum_{j=1}^{\ell_y} \Phi_j^{1,y} \mathbf{F}_1 T_{2,j}^\top 
    + \sum_{k=1}^{k_x} T_{3,k} \mathbf{F}_1 {\Sigma_k^{2,x}}^\top 
    + \sum_{l=1}^{k_y} \Sigma_l^{1,y} \mathbf{F}_1 T_{4,l}^\top 
    }_{\mathscr{L}(\mathbf{F}_1)}
    \right) = \bF_0.
    \label{Generalized_Sylv_multiterm}
\end{equation}
As before, we find $\mathscr{\tilde{L}}$ by averaging the diagonal terms $(\Phi_i^{2,x}, \Phi_j^{1,y}, \Sigma_k^{2,x}, \Sigma_l^{1,y})$ from \eqref{Generalized_Sylv_multiterm} 
%as follows:
%\begin{equation}
%\begin{aligned}
%    \alpha_{x,i} &= \text{average}(\Phi_i^{2,x}), \quad
%    \alpha_{y,j} &= \text{average}(\Phi_j^{1,y}), \quad
%    \gamma_{x,k} &= \text{average}(\Sigma_k^{2,x}), \quad
%    \gamma_{y,l} &= \text{average}(\Sigma_l^{1,y}),
%\end{aligned}
%\label{averaged_coeffs}
%\end{equation}
to arrive at the approximated SE:
\begin{equation}
  {P_1} \mathbf{F}_1 + \mathbf{F}_1 {P_2^\top} = \mathbf{F}_0,
    \label{Approximated_Sylv_multi}
\end{equation}
where
\begin{align}
\small
      P_1&= \sum_{i=1}^{\ell_x}  {A_{1,i}} + \sum_{k=1}^{k_x} {A_{3,k}}= \sum_{i=1}^{\ell_x} {\left(\frac{1}{\mathcal{R}_{\text{ranks}}} I_{N_1} - \Delta t \alpha_{x,i} T_{1,i}\right)} + \sum_{k=1}^{k_x} {\left(\frac{1}{\mathcal{R}_{\text{ranks}}} I_{N_1} - \Delta t \gamma_{x,k} T_{3,k}\right)} , \label{eq:P_1_gse_mr}\\ P_2 &= \sum_{j=1}^{\ell_y} {A_{2,j}} +\sum_{l=1}^{k_y} {A_{4,l}} = \sum_{j=1}^{\ell_y} {\left(\frac{1}{\mathcal{R}_{\text{ranks}}} I_{N_2} - \Delta t \alpha_{y,j} T_{2,j}\right)}+ \sum_{l=1}^{k_y} {\left(\frac{1}{\mathcal{R}_{\text{ranks}}} I_{N_2} - \Delta t \gamma_{y,l} T_{4,l}\right)\label{eq:P_2_gse_mr}} .
\end{align}

Also as before, the bases \( U_1 \) and \( V_1 \) are constructed from orthonormalization of dimension-wise xKrylov from equations \eqref{Generalized_Sylv_multiterm} and \eqref{Approximated_Sylv_multi} as:
\begin{align}
U_1 &=\kappa^{\texttt{trunc}}_m\left(P_1, P_1^{-1}, A^{-1}_{1,1:\ell_x}, A^{-1}_{3,1:k_x}, \Phi^{1,y}_{1:\ell_y}, \Sigma^{1,y}_{1:k_y}, U_0\right) , \\
V_1 &=\kappa^{\texttt{trunc}}_m\left(P_2, P_2^{-1}, A^{-1}_{2,1:\ell_y}, A^{-1}_{4,1:k_y}, \Phi^{2,x}_{1:\ell_x}, \Sigma^{2,x}_{1:k_x}, V_0\right).
\end{align}

In the multi-rank setting, after a Galerkin projection \(U^\top \mathbf{R}_{\mathscr{{L}}} V = 0\), we obtain a GSE similar to \eqref{Galerkin_S} for the matrix of coefficients \(\bS_1\):
\begin{equation}
     \bS_1 - \Delta_t \left( \sum_{i=1}^{\ell_x} \tilde{T}_{1,i} \bS_1 \tilde{\Phi}_{i}^{2,x \top} + \sum_{j=1}^{\ell_y} \tilde{\Phi}_{j}^{1,y} \bS_1 \tilde{T}_{2,j}^\top + \sum_{k=1}^{k_x} \tilde{T}_{3,k} \bS_1 \tilde{\Sigma}_{k}^{2,x \top} + \sum_{l=1}^{k_y} \tilde{\Sigma}_{l}^{1,y} \bS_1 \tilde{T}_{4,l}^\top \right) = \tilde{\bB}_1,
     \label{original_projected}
\end{equation}
where the projected operators are analogous to those defined in the rank-one case.
We solve \eqref{original_projected} for the \(\bS_1\) matrix using preconditioned GMRES. The preconditioner is constructed
by another Galerkin projection 
%\(U^\top \mathbf{R}_{\mathscr{\tilde{L}}} V = 0\), 
to find the projected SE:
\begin{equation}
  \tilde{P}_1 \mathbf{S} + \mathbf{S} \tilde{P}_2^\top = \tilde{\bB}_1,
\end{equation}
with:
\begin{align*}
    {\tilde{P}_1} &=  \sum_{i=1}^{\ell_x} \underbrace{\left(\frac{1}{\mathcal{R}_{\text{ranks}}}I_{r_x} - \Delta_t \alpha_{x,i} \tilde{T}_{1,i}\right)}_{\tilde{A}_{1,i}} + \sum_{k=1}^{k_x} \underbrace{\left(\frac{1}{\mathcal{R}_{\text{ranks}}}I_{r_x} - \Delta_t \gamma_{x,k} \tilde{T}_{3,k}\right)}_{\tilde{A}_{3,k}}, \\
    \tilde{P}_2 &=  \sum_{j=1}^{\ell_y} \underbrace{\left(\frac{1}{\mathcal{R}_{\text{ranks}}}I_{r_y} - \Delta_t \alpha_{y,j} \tilde{T}_{2,j}\right)}_{\tilde{A}_{2,j}} + \sum_{l=1}^{k_y} \underbrace{\left(\frac{1}{\mathcal{R}_{\text{ranks}}}I_{r_y} - \Delta_t \gamma_{y,l} \tilde{T}_{4,l}\right)}_{\tilde{A}_{4,l}}.
\end{align*}
This yields the multi-rank-coefficients ACS preconditioner operator:
\begin{equation}
  \mathscr{P} : \mathbf{S} \mapsto \tilde{P}_1 \mathbf{S} + \mathbf{S} \tilde{P}_2^\top .
      \label{nearby_projected}
\end{equation}

{We summarize the pseudo-algorithm for multi-rank coefficients in Algorithm \ref{alg:BE-LR}. The tolerances $\epsilon_{\kappa}$, $\epsilon_{\texttt{GMRES}}$, $\epsilon$, $\epsilon_{\texttt{tol}}$ in the algorithm are set as follows.}
\HHONE{First, we choose the residual tolerance for {xKrylov outer iterations}, $\epsilon_{\texttt{tol}}$.
% \QQ{$\epsilon_{\kappa}$???}.\lc{Set in the next sentence...}
To prevent rapid growth of the xKrylov subspaces, 
% \sout{especially when the matrix equations contain many terms,}
we choose the Krylov basis truncation threshold as $\epsilon_{\kappa}  \leq\epsilon_{\texttt{tol}}$, typically two to three orders of magnitude smaller.
%to truncate basis vectors that are nearly linearly dependent. This allows us to maintain a compact basis while effectively treating singular values of small magnitude as an affordable error. However, over-truncation may risk removing important information, slowing down convergence. 
The GMRES-ACS tolerance for the inner reduced system solve, $\epsilon_{\texttt{GMRES}}$, is initially set to be less than $\epsilon_{\texttt{tol}}$ for the first outer-basis iteration. It is subsequently adjusted to be two to three orders of magnitude smaller than the relative residual error from the previous {xKrylov} outer iteration. This adjustment ensures that errors from the inner solve do not pollute the accuracy of the outer solve. The SVD truncation threshold $\epsilon$ is chosen to be one to three orders of magnitude smaller than the tolerance \(\epsilon_{\texttt{tol}}\). 
Since solution truncation occurs before residual evaluations, the truncation error in the solution should remain below the residual tolerance to prevent residual pollution by over-truncation.} 
% \sout{Otherwise, the computed solution may meet the residual threshold, but over-truncation could prevent this acceptance.}

%\HH{Keep this latter as a safeguard for the user?}} \lc{Should be OK to keep...}
% \sout{Given this choice of tolerances, the whole adaptive-rank process is summarized in Algorithm \ref{alg:BE-LR}.} \HH{Nowhere else is the algorithm introduced.}

\begin{algorithm}[t]
\caption[Backward Euler Adaptive-Rank Integrator for the GSE]{Backward Euler Adaptive-Rank Integrator for the GSE. 
%\footnotetext{A MATLAB implementation of this algorithm is available upon request.}
}
\label{alg:BE-LR}
\SetAlgoNlRelativeSize{-2}
\SetNlSty{textbf}{(}{)}
\tcp{This algorithm solves Eq.~\eqref{Generalized_Sylv_multiterm} for $\mathbf{F}_1 = U_1 S_1 V_1^\top$, where $\mathbf{F}_0 = U_0 S_0 V_0^\top$}
\KwIn{
    Initial condition matrices $U_0$, $V_0$, $S_0$;\\
    Operators $P_1$, $P_2$, $\{T_{1,i}\}_{i=1}^{\ell_x}$, $\{T_{2,j}\}_{j=1}^{\ell_y}$, $\{T_{3,k}\}_{k=1}^{k_x}$, $\{T_{4,l}\}_{l=1}^{k_y}$;\\
    Tolerances $\epsilon_{\kappa}$, $\epsilon_{\texttt{GMRES}}$, $\epsilon$, $\epsilon_{\texttt{tol}}$;\\
    % Maximum iterations \texttt{max\_iter};
}
\KwOut{
    Updated bases $U_1$, $V_1$;\\
    Truncated singular values $S_1$
    \BlankLine
}Compute operators $\{A_{1,i}\}_{i=1}^{\ell_x}$, $\{A_{2,j}\}_{j=1}^{\ell_y}$, $\{A_{3,k}\}_{k=1}^{k_x}$, $\{A_{4,l}\}_{l=1}^{k_y}$ according to Eqs.~\eqref{eq:P_1_gse_mr}, \eqref{eq:P_2_gse_mr}\; 
\HHONE{Compute $l=2+\ell_x+\ell_y+k_x+k_y$\;
Set \(U^{(i)}=U_0\) and \(V^{(i)}=V_0\) for \(i=1, \dots, l\)\;
Set $U_1=U_0 \quad V_1=V_0$\;}
\While{not converged}{
   \tcp*[h]{Step K1.}\\
  \HHONE{Orthogonalize and} truncate to tolerance $\epsilon_{\kappa}$\;
\HHONE{$U_1, \{U^{(1)}\cdots U^{(l)}\} \leftarrow \kappa_m^{\texttt{trunc}}
\left(\{P_1, P_1^{-1}, A^{-1}_{1,1: \ell_x}, A^{-1}_{3,1: k_x},  \Phi^{1,y}_{1: \ell_y},\Sigma^{1,y}_{1: k_y} \}; \{U^{(1)}\cdots U^{(l)}\};U_1; \epsilon_\kappa \right)$\label{alg:line:qr1}\;}
\HHONE{$V_1, \{V^{(1)}\cdots V^{(l)}\}  \leftarrow \kappa_m^{\texttt{trunc}}\left(\{P_2,P_2^{-1}, A^{-1}_{2,1: \ell_y}, A^{-1}_{4,1: k_y}, \Phi^{2,x}_{1: \ell_x},\Sigma^{2,x}_{1: k_x}\} ;\{V^{(1)}\cdots V^{(l)}\};V_1;\epsilon_\kappa\right)$\label{alg:line:qr2}\; }

  \tcp*[h]{Step K2.}\\
  Solve the reduced GSE~\eqref{original_projected} for $\bS_1$ using GMRES-ACS \eqref{nearby_projected} to tolerance $\epsilon_{\texttt {GMRES}}$\;
     \tcp*[h]{Step K3.}\\
{Truncate $\{U_{1,\epsilon},S_{1,\epsilon},V_{1,\epsilon}\}\leftarrow \mathcal{T}_{\epsilon}(\{U_1,S_1,V_1\})$\label{alg:line:trunc} to tolerance $\epsilon$} (Alg. \ref{alg:trunc})\; 
   
   \tcp*[h]{Step K4.}\\
   {Compute $\{ \_, R_U \} = \texttt{QR}\left([U_1, U_{1,\epsilon}, T_{1,1:\ell_x}U_{1,\epsilon}, \Phi_{1:\ell_y}^{1,y}U_{1,\epsilon}, T_{3,1:k_x}U_{1,\epsilon}, \Sigma_{1:k_y}^{1,y}U_{1,\epsilon}]\right)$ and\\
  $\{\_, R_V\} = \texttt{QR}\left([V_1,V_{1,\epsilon} ,\Phi_{1:\ell_x}^{2,x}V_{1,\epsilon}, T_{2,1:\ell_y}V_{1,\epsilon}, \Sigma_{1:k_x}^{2,x}V_{1,\epsilon}, T_{4,1:k_y}V_{1,\epsilon}]\right)$ }\;\label{alg:line:qr3}
  
   {Compute the residual $\|\mathbf{R}\| = \left\| R_U \text{diag}\left(- \tilde{B}_1,S_{1,\epsilon}, -\Delta t \left( I_{\mathcal{R}_{\text{ranks}}} \otimes S_{1,\epsilon}\right) \right) R_V^\top \right\|$}\; \label{alg:line:residual}
   
   \If{$\|\mathbf{R}\|/\|\mathbf{F}_0\| \geq \epsilon_{\texttt{tol}}$}{
      Reject solution and cycle to augment bases further\; 
   }   
   \Else{
    {  Set $U_1\leftarrow U_{1,\epsilon}$; $U_1\leftarrow V_{1,\epsilon}$; $S_1=S_{1,\epsilon}$\;}
         Break\;
      Exit loop\;
   }
}
\end{algorithm}

\subsection{Analysis of Computational \HHREV{and Memory} Complexity\label{complexity_analysis}}
%\lc{We need to add memory complexity analysis here.}\HH{Done. Please take a look.}
We analyze next the computational complexity \HHREV{and memory storage} of Algorithm \ref{alg:BE-LR} for backward Euler. For conciseness, we consider the same spatial resolution per dimension, i.e., \( N = N_1 = N_2 \).
Starting from an initial condition \( U_0 S_0 V_0^\top \), where \( U_0 \in \mathbb{R}^{N \times r_0} \), \( V_0 \in \mathbb{R}^{N \times r_0} \), and \( S_0 \in \mathbb{R}^{r_0 \times r_0} \), the \( m \)-th xKrylov iteration generates \( U_1 \in \mathbb{R}^{N \times  r_\kappa} \), \( V_1 \in \mathbb{R}^{N \times  r_\kappa} \), and \( S_1 \in \mathbb{R}^{ r_\kappa \times  r_\kappa} \), where \(  r_\kappa = (2+\mathcal{R}_{\text{ranks}})r^{\ast}_0 +r^{(m-1)} \). Here, $r^{(m-1)}$ denotes the basis size from the previous xKrylov iteration (with $r^{(0)}=r_0$), and $r_0^{\ast}\le r_0$ is the (possibly truncated) initial rank produced at line~\eqref{ind_trunc_line_1} of Algorithm~\ref{alg:basis_constr}. \HHREV{While not used in this work, one may impose a maximum allowable rank $r_{\max}$ in the orthogonalization step of Algorithm \ref{alg:basis_trunc} to prevent the subspace from exceeding available computational resources.} 
% \lc{Looks good} 
Using SVD-truncated QR, this dimension is then reduced to smaller ranks, \( r_x \) and \( r_y \). Let \( \tilde{r} = \max(r_x, r_y) \). The matrix \( S_1 \) from the \( m \)-th inner iteration solve, initially of size \( r_x \times r_y \), is truncated to rank \( r_1 \leq \min(r_x, r_y) \) using Algorithm \ref{alg:trunc}. To simplify the complexity analysis, let %\sout{\( r = \max(r_x, r_y, \mathcal{R}_{\text{ranks}}r_1) \)} 
\( r = \max(r_0, r^{(m-1)}, \tilde{r}, (\mathcal{R}_{\text{ranks}}+1)r_1) \), {which we consider independent of the mesh-refinement $N$ for practical purposes}. 
The computational complexity of each step is then estimated as follows:

\paragraph{Step K1}
Constructing the generalized Krylov basis as outlined in lines \eqref{alg:line:qr1} and \eqref{alg:line:qr2} of Algorithm \ref{alg:BE-LR} requires performing the operations \( A_{1,1:r_x}^{-1}U_0 \) and \( A_{3,1:k_x}^{-1}U_0 \), and \( A_{2,1:r_y}^{-1}V_0 \) and \( A_{4,1:k_y}^{-1}V_0 \), which require solving the systems \( A X = U \) and  \( A X = V \) (subscripts have been omitted for simplicity). With 1D finite-differences, \( A \) has a tridiagonal structure, which can be efficiently factorized using direct-solver techniques such as the Thomas algorithm, with \(\mathcal{O}(N)\) complexity per column of \( X \) \cite{quarteroni2006numerical}. {The same complexity holds for applying $P_1^{-1}$ and $P_2^{-1}$.}
Additionally, applying \( P_1 U_0 \), \( \Phi_{1,3} U_0 \), \( \Sigma_{1,3} U_0 \), \( P_2 V_0 \), \( \Phi_{2,4} V_0 \), and \( \Sigma_{2,4} V_0 \), can be done with \( \mathcal{O}(Nr^\ast_0) \) computational complexity for sparse diagonal and tridiagonal matrices, leading to a total complexity of \( \mathcal{O}(N(2+\mathcal{R}_{\text{ranks}})r^\ast_0) \) {to augment the xKrylov subspaces}. The use of Algorithm~\ref{alg:basis_constr} to compute the orthonormal truncated bases \( U_1 \) and \( V_1 \) incurs a complexity of %\sout{\( \mathcal{O}(N r_{\kappa}^2) \)} 
\( \mathcal{O}\left( N \, \max\left ( r_\kappa\, r^{(m-1)}, \; r_\kappa^2 \right ) \right) \sim \mathcal{O}(N r^2)\). \HHREV{The memory complexity to store the truncated bases is \( \mathcal{O}\left( N \, \left ( r^{(m-1)}+ r_\kappa  \right ) \right) \sim \mathcal{O}(N r)\)}. %\HH{If this level of detail is not needed, consider replacing previous expression with \(\mathcal{O}(N r^2)\)}
% Specifically The SVD decomposition of the small-sized upper triangular matrix in line \eqref{alg1:line:svd:basis} requires complexity \( \mathcal{O}\left(  r_\kappa^3\right)\).
%Hence, the complexities of lines \eqref{alg:line:qr1}, \eqref{alg:line:qr2} in Algorithm \ref{alg:BE-LR} scale as \sout{\( \mathcal{O}\left(N  r_\kappa^2\right) \)} \( \mathcal{O}\left(N  r^2\right) \). 
\HHONE{Generally, we expect %\sout{\(  r_\kappa \)} 
\(  r\) to be much smaller and independent of $N$ when the solution is sufficiently smooth and exhibits a low-rank structure.}
 % \HH{Should we change to: Generally, we expect \(  r_\kappa \) to be much smaller and independent of $N$ when the solution is sufficiently smooth and exhibits a low-rank structure.} \lc{Ok}
% \lc{This discussion should be placed in the introduction}

\paragraph{Step K2}
The application of the reduced GSE operator using matrix-matrix multiplication and summation in the Krylov method (GMRES) features a complexity of $\mathcal{O}\left( \tilde{r}^{d+1} \right)$. Applying the ACS preconditioner requires solving an SE of size $r_x \times r_y$, which can be achieved using the Bartels-Stewart algorithm \cite{golub2013matrix} with a complexity of \( \mathcal{O}(\tilde{r}^{d+1})\sim \mathcal{O}(r^{d+1}) \). Hence, the computation of the $S$-step has $\mathcal{O}( r^{d+1})$ ~ complexity. \HHREV{The memory complexity to store the S matrix/tensor is \( \mathcal{O}\left(\tilde{r}^d \right) \sim \mathcal{O}(r^d)\).}

\paragraph{Step K3}
The SVD decomposition of \( {\widetilde{S}} \) in line \eqref{alg1:line:svd} of Algorithm \ref{alg:trunc}  is of \( \mathcal{O}(\tilde{r}^{d+1}) \) complexity. Updating the bases in lines \eqref{alg1:line:mm1} and \eqref{alg1:line:mm2} of Algorithm \ref{alg:trunc} requires \( \mathcal{O}( N \tilde{r}{r_1})\) operations{, where ${r_1}\leq  \min(r_x, r_y)$ is the post-truncation rank of the new solution from Algorithm \ref{alg:trunc}}.
\paragraph{Step K4}
The QR decomposition to obtain $R_U$ and $R_V$ in lines (11)--\eqref{alg:line:qr3} of Algorithm~\ref{alg:BE-LR} has a complexity of $\mathcal{O}\left( N \left( \tilde{r} + (\mathcal{R}_{\text{ranks}} + 1) r_1 \right)^2 \right)  \sim \mathcal{O}(Nr^2) $.  
The residual computation in line~\eqref{alg:line:residual} requires $\mathcal{O}\left( \left( \tilde{r} + (\mathcal{R}_{\text{ranks}} + 1) r_1 \right)^{d+1} \right)\sim \mathcal{O}(\mathcal{R}_{\text{ranks}}^{d+1}r^{d+1})$ due to matrix multiplications, with an additional $\mathcal{O}\left( \left( \tilde{r} + (\mathcal{R}_{\text{ranks}}+1) r_1 \right)^d \right)\sim \mathcal{O}(\mathcal{R}_{\text{ranks}}^{d}r^{d})$ for the Frobenius norm.  
\HHREV{The memory storage requirements for the residual tensor may be quite prohibitive if stored fully [it would scale as $\mathcal{O}(\mathcal{R}_\text{ranks}^d r^d)$]. Instead, we  store only the $R_U$ and $R_V$ factor matrices and the $S$ matrix, and build the residual for the norm calculation on the fly still with computational complexity $\mathcal{O}(\mathcal{R}_{\text{ranks}}^{d+1}r^{d+1})$. Storing $R_U$ and $R_V$ requires \( \mathcal{O}\left(\left(  \tilde{r} + (\mathcal{R}_{\text{ranks}} + 1) r_1 \right)^2 \right)\sim \mathcal{O}(\mathcal{R}_{\text{ranks}}^2r^2)\), and the residual tensor blocks in Eq. \ref{eq:Residual} require
%\( \mathcal{O}\left(  \tilde{r}^d \right) + \mathcal{O}\left(  r_1^d \right) \sim 
\(\mathcal{O}(r^d)\), resulting in a total storage of $\mathcal{O}(\mathcal{R}_{\text{ranks}}^2r^2) + \mathcal{O}(r^d)$. Nevertheless, it is clear from this analysis that the residual evaluation may become prohibitively expensive when $\mathcal{R}_\text{ranks} \gg 1$. In such situations, the difference between successive solution iterates may be used as a alternate measure of convergence~\cite{palitta2025subspace}, commensurate in cost to the rest of the algorithm.}
% The latter approach is very cheap to compute and very reliable, as demonstrated numerically later in this study.
%\lc{Comment: unless we propose an alternative for the residual evaluation, the computational complexity and memory scalings below would need to be corrected to include residual terms...}

% \sout{We note that henceforth \( r \) will represent the maximum of \(r_\kappa\) and \(\mathcal{R}_{\text{ranks}}r_1\). \HHONE{In principle, the rank $r$ may depend on the grid size $N$, as theoretically justified by Kressner et al.~\cite{kressner2010krylov,beckermann2013error}. However, their numerical tests indicate that this theoretical dependence is overly pessimistic. Our numerical results and practical experience further confirm that this dependence appears only for extremely fine grids and exceedingly tight tolerances, both of which are uncommon in typical computational demands. Hence, if $r$ is independent of $N$, 
%the solution possesses a low-rank structure, and the constructed approximation subspace efficiently maintains moderate $r$, 
% the complexity scales linearly with respect to $N$, cubically with respect to $r$ in 2D, and quartically with respect to $r$ in 3D, as numerically verified in Section~\ref{section_4}. These scalings motivate keeping the overall rank $r$ as low as possible.} \HH{Potential redundancy with point made in 3.2.}}\QQ{this is the place to summarize the complexity, so it is okay to leave these here?}\lc{Yes, but much simplified... What about the following?}

{It follows that the overall complexity of the algorithm scales as $\mathcal{O}(Nr^2) +  \mathcal{O}(r^{d+1})$, \HHREV{and the memory complexity scales to leading order as $\mathcal{O}(Nr)  + \mathcal{O}(r^d)$,} with $d$ the dimensionality.  Both claims will be numerically verified in Section~\ref{section_4}.}

\section{Numerical Experiments}
\label{section_4}
%\lc{I did not find a reference to the high-order DIRK implementation in the appendix anywhere?}
We demonstrate 
% \sout{the computational complexity estimates and}
the overall efficacy of the algorithm with several 2D and 3D numerical examples. 
% \sout{The 3D examples demonstrate that the key conclusions of the 2D study above (namely, effective rank adaptivity via Krylov subspace search, linear scaling of the computational complexity with one-dimensional resolution $N$, effectiveness of ACS-preconditioner, and manageable scaling with the rank $r$ with dimension $d$ of $r^{d+1}$)  carry over seamlessly to higher dimensional settings.}
We assume homogeneous Dirichlet boundary conditions for simplicity; \HHONE{other boundary conditions (e.g., periodic, Neumann, etc.) can be easily incorporated in the formulation through the operators $T_{1:4}$ \eqref{operator_T}.} 
% \sout{However, special care must be taken if the operator has a null space. In these cases, the null space should either be explicitly included during the Krylov subspace construction or addressed through a post-processing correction, such as LoMaC \cite{guo2023local}. For additional details, see the discussion in \cite{el2024krylov}.}
All simulations are performed using MATLAB. For efficient computational handling in 3D, we employ tools from the Tensor Toolbox \cite{tensor_toolbox} and the Htucker package \cite{kressner2012htucker}. \HHREV{High-order temporal integration using DIRK is used throughout the numerical section. Details on the  the DIRK scheme implementation in our adaptive-rank algorithm are provided in Appendix \ref{Appendix_DIRK}.}
% \HH{
% For benchmarking purposes, we compare against full-rank algorithm using the Kronecker formulation inverted with MATLAB's backslash command.} \lc{What are we benchmarking? Solution quality or wall-clock time? let's be specific.}
Throughout this section, we denote the diffusion Courant number and the advection Courant-Friedrichs-Lewy (CFL) number by $\lambda_D$ and $\lambda_A$, respectively, defined as:
\begin{equation} \lambda_{D} = \frac{\Delta t \phi_{\text{max}}}{\Delta x^2}, \qquad \lambda_{A} = \frac{\Delta t \sigma_{\text{max}}}{\Delta x}, \end{equation}
where $\phi_{\text{max}}$ and $\sigma_{\text{max}}$ are the maximum of the diffusion and advection coefficients, respectively. \HHREV{The numbers $\lambda_D$ and $\lambda_A$ are appropriate measures of the maximum eigenvalue  of the associated GSE, and by extension of its condition number since the minimum eigenvalue is of $\mathcal{O}(1)$}.
\HHONE{{For all simulations presented in this section, we did not need to set a maximum number of outer xKrylov iterations because the algorithm consistently converged within a small number of iterations.}}
\subsection*{\textbf{Example 4.1:} 3D Advection-Diffusion}\label{Ex_4.2}
% \lc{Focus this section on 3D only; 2D discussion is just offuscating things}
% \sout{In the 3D setting,} 
We consider a rank-1 {variable} advection coefficient in a 3D domain $x, y, z \in [-1, 1]$:
\[
\begin{aligned}
\sigma^x(x,y,z) &= \sigma^1_1(x) \sigma^2_1(y)\sigma^3_1(z), & \sigma^1_{1}(x) &= 1 - x^2, & \sigma^2_{1}(y) &= 2y, & \sigma^3_{1}(z) &= -2z \\
\sigma^y(x,y,z) &= \sigma^1_2(x) \sigma^2_2(y)\sigma^3_2(z), & \sigma^1_{2}(x) &= -2x, & \sigma^2_{2}(y) &= 1 - y^2, & \sigma^3_2(z)&=2z \\
\sigma^z(x,y,z) &= \sigma^1_3(x) \sigma^2_3(y)\sigma^3_3(z), & \sigma^1_{2}(x) &= 4x, & \sigma^2_{2}(y) &= 2y, & \sigma^3_3(z)&=1-z^2. 
\end{aligned}
\]
\HHREV{Note that this flow is divergence-free and therefore a pure vorticity flow, which results in interesting flow dynamics such as differential rotation.} For diffusion, we specify a rank-3 {variable} coefficient as follows:
\[
\phi^{x}(x,y,z) = \phi^{y}(x,y,z) = \sum_{i=1}^{3} \phi^1_i(x) \phi^2_i(y) \phi^3_i(z),
\]
\[
\begin{aligned}
{\phi}^1_{1}(x) &= \exp\left(-\left(x - 0.3 \sin(x)\right)^2\right), & 
{\phi}^2_{1}(y) &= \exp\left(-\left(y - 0.3 \cos(y)\right)^2\right), & 
{\phi}^3_{1}(z) &= \exp\left(-\left(y - 0.3 \sin(z)\right)^2\right),\\
{\phi}^1_{2}(x) &= \exp\left(-\left(x - 0.6 \sin(\pi x)\right)^2\right), & 
{\phi}^2_{2}(y) &= \exp\left(-\left(y - 0.6 \sin(\pi y)\right)^2\right),  & 
{\phi}^3_{2}(z) &= \exp\left(-\left(y - 0.6 \sin(\pi z)\right)^2\right),\\
{\phi}^1_{3}(x) &= \exp\left(-\left(x - 0.6 \sin(2\pi x)\right)^2\right), & 
{\phi}^2_{3}(y) &= \exp\left(-\left(y - 0.6 \sin(2\pi y)\right)^2\right) & 
{\phi}^3_{3}(z) &= \exp\left(-\left(y - 0.6 \sin(2\pi z)\right)^2\right).
\end{aligned}
\]
\HHREV{These diffusion coefficients are relatively balanced and therefore mostly isotropic. We will investigate anisotropic diffusion in future work.}
%\HHREV{\textbf{Remark:} In this example the discretization matrices are chosen to behave similarly across all modes/dimensions. We study the convergence of the xKrylov subspace and the computational scaling-measured in terms of wall-clock time and memory usage with respect to one-dimensional grid refinements. The assumption of identical behavior across dimensions allows us to focus parametrization along a single one. This is purely a design choice and does not impose any inherent limitation on the adaptive-rank algorithm.}
The initial condition is a rank-2 Gaussian:
\begin{align*}
f(x, y, z, t=0) &= 0.5 \exp\left(-400 \left((x-0.3)^2 + (y-0.35)^2 + (z-0.2)^2\right)\right) \\ &+ 0.8 \exp\left(-400 \left((x-0.65)^2 + (y-0.5)^2+ (z-0.55)^2\right)\right).     
\end{align*}
{%In this example, for the 2D version of Equation \eqref{eq:adv-diff}, 

% \HH{Move this section to the Computational section. Do not put $f$. }In the 2D setting, the advection coefficients $\sigma^x(x,y)=\sigma^1_1(x)\sigma^2_1(y)$ and $\sigma^y(x,y)=\sigma^1_2(x)\sigma^2_2(y)$, and the diffusion coefficients are simplified to $\phi^x(x,y) = \phi^y(x,y) = \sum_{i=1}^{3} \phi^1_i(x) \phi^2_i(y)$, with the coefficients set similarly to their 3D counterparts.
% }
% The initial condition for the 2D case is defined similarly as a rank-2 Gaussian:
% \[
% f(x, y, t=0) = 0.5 \exp\left(-400 \left((x-0.3)^2 + (y-0.35)^2 \right)\right) + 0.8 \exp\left(-400 \left((x-0.65)^2 + (y-0.5)^2\right)\right).
% \]
% \QQ{i suggest description for both 2D (first) and 3D cases, but point out for each case, what points that we will demonstrate.}
We use a spatial grid with \( N=300 \) in all dimensions, and the time-step \( \Delta t \) is set to \( \lambda_D \Delta x^2 \), with \( \lambda_D \) ranging from $140$ to $1125$. We test the adaptive-rank algorithm using three different integrators: backward Euler (BE), DIRK2, and DIRK3, as detailed in Tables~\ref{Butcher:DIRK2} and \ref{Butcher:DIRK3} in Appendix \ref{Appendix_DIRK}. The residual tolerance is \( \epsilon_{\texttt{tol}} = 10^{-4} \) for BE, DIRK2, and DIRK3, the basis truncation threshold is \(\epsilon_{\kappa} = 10^{-6}\), the SVD truncation threshold is \(\epsilon = 10^{-6}\), and the GMRES acceptance threshold is $\epsilon_{\texttt{GMRES}}=10^{-6}$. 

We investigate numerically the following: (1) high-order temporal convergence, (2) effectiveness of the xKrylov iterations in capturing the range of solution, (3) performance of the ACS preconditioner to solve the reduced GSE, and (4) verification of the overall computational complexity.} 

\begin{figure}[t]
    \centering    \includegraphics[width=\linewidth]{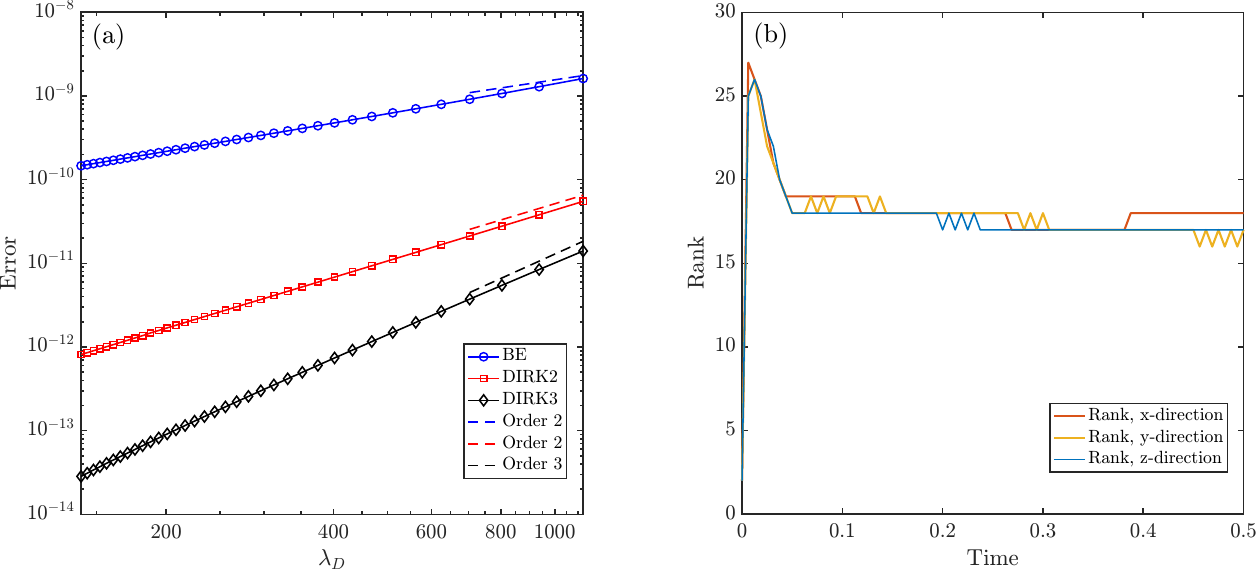}
    \caption{
    {Example 4.1. (a) Temporal convergence study of BE, DIRK2, and DIRK3 schemes. (b) Rank evolution of the solution in the \(x\), \(y\), and \(z\) directions over time for \(\lambda_D = 140\).}}
    \label{fig:convergence_comparison}
\end{figure}

\textbf{High-Order Temporal Convergence:} Figure~\ref{fig:convergence_comparison}(a) illustrates the \( L_1 \) error norm vs time stepping size (characterized by $\lambda_D$) by showcasing the performance of the adaptive-rank integrators for BE, DIRK2, and DIRK3. The expected asymptotic order of convergence of the temporal error for the adaptive-rank integrators, depicted by lines with markers, is observed. Additionally, the algorithm adeptly captures
the rank evolution of the solution in the $x$, $y$, and $z$ directions, as illustrated in Figure~\ref{fig:convergence_comparison}(b).

\textbf{Effectiveness of the xKrylov Iterations:} 
% \QQ{since we are investigating the effectiveness of xKrylov, shall we include a figure such as Fig 7 for swirling deformation?}\HH{I am concerned that we will be conveying too many messages in this example. The inclusion of the xKrylov history will not deliver any new insights to what is already delivered in Figure 8. } \lc{I agree that we should limit \# of messages...}
\HHONE{
% \sout{We conduct a subspace convergence analysis in the 3D setting.}
We evaluate the xKrylov subspace dimension required to achieve tolerances of \(\epsilon_{\texttt{tol}}= 10^{-4}, 10^{-5}, 10^{-6}\). We perform a single time-step integration with a fixed step size \(\Delta t=10^{-3}\) with grid refinements in $N$.} \HHREV{%We base our analysis on the assumption that the stiffest component in the time-dependent PDE is the diffusion operator. Consequently, the primary convergence behavior of the xKrylov subspace is governed by resolving the eigenvalues of this stiff operator. For the case of smooth variable diffusion coefficients, these coefficients can be viewed as perturbations of the eigenvalues of the Laplacian operator. Hence, we take the Laplacian as the reference operator to study the convergence of the xKrylov approximation spaces. 
The convergence bound for the extended Krylov method depends on the condition number of the advection-diffusion operator, $\kappa$. Assuming that the diffusion operator dominates the condition number, and that the variable coefficients perturb its eigenvalues only by factors of order unity, $\lambda_D$ provides a good estimate of the condition number (i.e., $\kappa = \mathcal{O}(\lambda_D)$ in the case of the diffusion equation), and is therefore a good figure of merit to assess the convergence of the xKrylov subspaces. As reported in \cite{kressner2010krylov}, the convergence factor for the extended Krylov subspace method applied to an equation with Kronecker product structure is $\kappa^{1/4}$ in 2D and $\kappa^{2/3}$ in 3D. 
%For the implicit schemes considered in this work, we typically have $\lambda_D \gg 1$, leading to a convergence rate of $\mathcal{O}(\lambda_D^{1/4})$ in 2D and $\mathcal{O}(\lambda_D^{2/3})$ in 3D.
}
% \QQ{where does " $\kappa \approx 1 + \frac{8\lambda_D}{d}\sim \lambda_D$" come from; and how "$\kappa^{1/4}$ in 2D" correspond to the slope in the figure, since it is a 3D problem?}\HH{From my own derivation. I added details in Appendix \ref{Appendix_xKrylov} and references above.} \HH{Move to referee response. $\kappa \sim \mathcal{O}(\lambda_D)$} \lc{Yes, that is what we had before; having the bit of detail does not hurt, given the referee's complaints, but I am open to removing it if we decide so.}

\begin{figure}[h]
    \centering \includegraphics[width=0.6\linewidth]{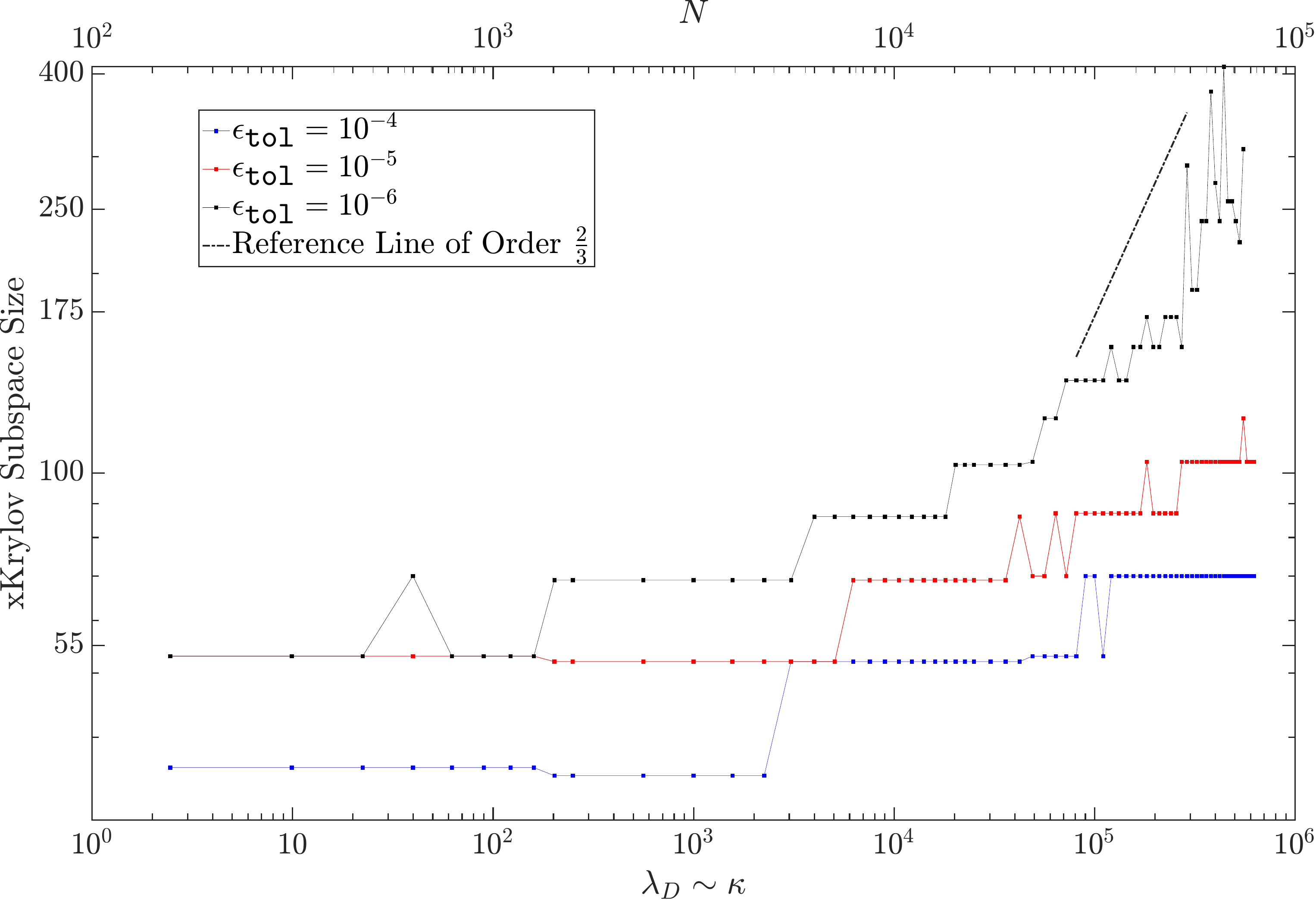}
     \caption{Example 4.1. Convergence study using the xKrylov subspace method with DIRK3 adaptive-rank integrator {for several outer basis tolerances, \(\epsilon_{tol}\)}.}
     % \HH{Should I mentioned reference line of order 1/4 wrt bottom axis, or should it be obvious to the reader?}\lc{I would put $\lambda_D^{1/4}$ in the legend instead of ``order 1/4''}\HH{Done}
    \label{fig:krylov_size_study}
\end{figure}

To isolate the effect of \(\lambda_D\) on the convergence of xKrylov, we mitigate other potential influencing factors, such as excessive truncation of the Krylov subspace, insufficient convergence of the GMRES-ACS solver, and truncation errors in the constructed solution, by setting related thresholds to \(\epsilon = \epsilon_{\texttt{GMRES}} = \epsilon_\kappa = 10^{-10}\).
Figure~\ref{fig:krylov_size_study} illustrates the xKrylov subspace dimension required for convergence as a function of the diffusion Courant number \(\lambda_D\) for a fixed \(\Delta t\) and varying the spatial 1D resolution, \(N\). Our results suggest that, for a fixed tolerance, the convergence of the Krylov subspaces remains largely insensitive to grid refinement across a broad range of \(\lambda_D\). For extremely fine 1D mesh resolutions, $N \sim 10^4$, beyond typical practical computational needs, the convergence rate of xKrylov scales as  \(\lambda_D^{2/3}\),   consistently with the theoretically expected $\kappa^{2/3}$ scaling \HHREV{in 3D}~\cite{kressner2010krylov,beckermann2013error} . 
% \QQ{Figure 2 instead?} \HH{Remove 1/4 line slope, not relevant.} \lc{Well, it is not clear we are agreeing with the 2/3 line scaling either... we are somewhere in between}%\lc{Shouldn't this be $\lambda_D^{2/3}$?} \HH{Great point. In 3D the complexity should be indeed $\lambda_D^{2/3}$. However, what we are seeing is that the complexity for 3D behaves according to the 2D rate. I belive this is due to the non-sharpness of the convergence rates. Similar result is exhibited in Kressner's paper, see Figure 7.3 (dims 5 and 10).}
% \QQ{more fundamentally, it is about $N$? Not as much about  \(\lambda_D\)?} 
% \QQ{is this 2d or 3d?}

\begin{figure}[h]
    \centering    \includegraphics[width=\linewidth]{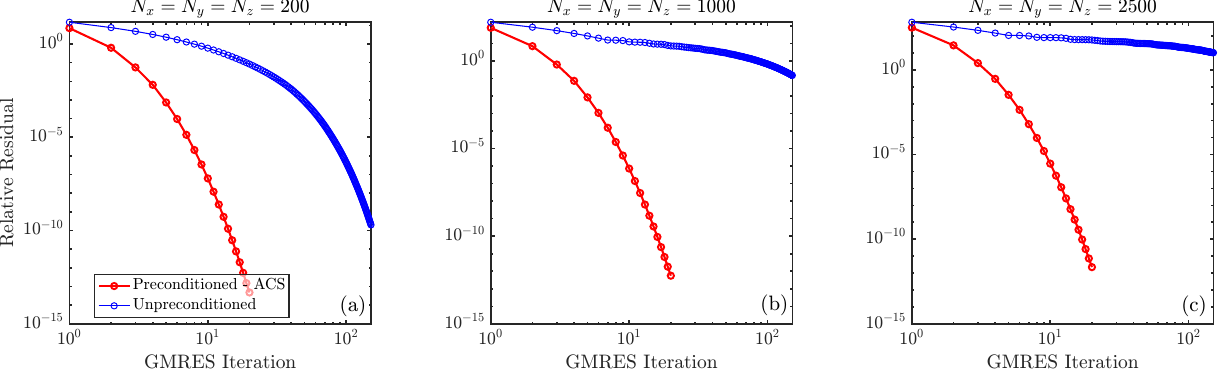}
\caption{Example 4.1. Log-log plots of GMRES residuals versus iterations, with and without preconditioning. Simulations are conducted for a single time step with \(\Delta t = 0.01\) at the third outer-basis augmentation, using grid resolutions \(N_x = N_y = N_z = 200\) in (a), \(1000\) in (b), and \(2500\) in (c).
% \lc{I would remove 2D results; do not add anything.}\HH{Done.}
}

    \label{prec_GMRES}
\end{figure}

\textbf{Performance of ACS-Preconditioner:}
\HHREV{Ill-conditioning for the projected system stems from both grid resolution (the projected system becomes more ill-conditioned as we refine the mesh) and rank growth (the projected system becomes larger in size). For effectiveness, the preconditioner must be able to deal with both sources of ill-conditioning.} We consider first the scaling of the ACS-preconditioned GMRES inner solve with respect to the 1D resolution $N$. For this, we compare performance with and without the ACS preconditioner at the first time step (\(\Delta t = 0.01\)) and the third outer-basis augmentation (i.e., keeping the subspace size constant, and therefore the rank), varying the mesh size. 
% \HHONE{We aim to demonstrate numerically that the convergence of the ACS preconditioner is superlinear and independent of both the grid resolution $N$ and the xKrylov subspace dimension $r_\kappa$, as in a compact-equivalent preconditioner~\cite{axelsson2009equivalent}.} \QQ{this statement is repeated below. maybe we could remove this statement here?}
Figure~\ref{prec_GMRES} illustrates the GMRES residual history over iterations for ACS-preconditioned and unpreconditioned cases. The results confirm that ACS-preconditioned GMRES delivers mesh-independent and super-linear convergence rates, which is not the case for unpreconditioned GMRES. 
% \QQ{"significantly accelerates the super-linear convergence"?}
% \HH{As per group discussion, GMRES gives better than linear convergence, however, the preconditioner further accelerates this superlinear convergence and makes independent of $r$ and $N$. Does this make sense, Dr. Qiu? Or should I reformulate?}
% \HHONE{This latter attributed to the conditioning of the advection and diffusion operators, which deteriorates as the grid is refined, negatively affecting the conditioning of the Galerkin-projected operators. }

\begin{figure}[h]
    \centering
    \includegraphics[width=\linewidth]{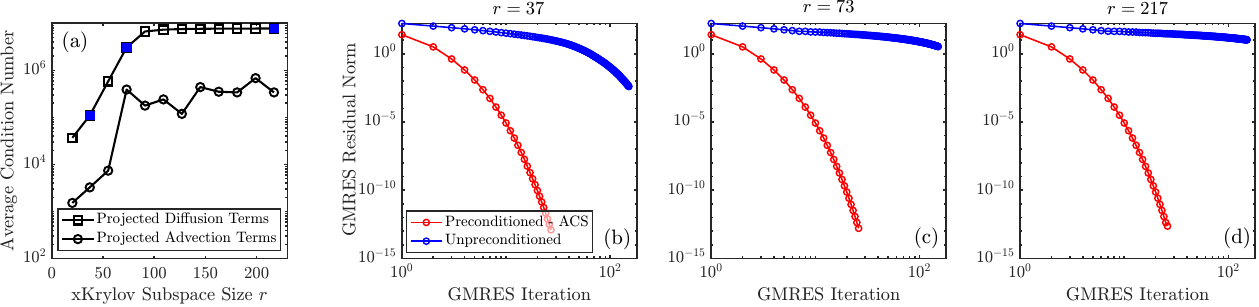}
\caption{Example 4.1. (a) Average condition number for Galerkin-projected advection and diffusion operators. (b)-(d) Log-log plots of GMRES residuals versus iterations, with and without preconditioning for various subspace sizes $r_\kappa = 37$, $r_\kappa = 73$, and $r_\kappa = 217$, respectively. These cases correspond to the highlighted blue data points in (a). Simulations are conducted for a single time step with \(\Delta t = 0.01\) and a grid of \(N_x = N_y = N_z = 1000\).
% \lc{In a perfect world these results should be 3D}\HH{Done. Perfection achieved.}
}
    \label{prec_GMRES_basis}
\end{figure}

\HHONE{
% \sout{In the following simulation, the grid size is fixed at $N = 1000$ in 3D, and a single time step of $\Delta t = 0.01$ is used. The only varying metric is the dimension of the xKrylov subspace, which we let grow without stopping criterion.}
\HHREV{We consider next the performance of the GMRES inner solve with respect to the solution rank $r$. For this test, the grid size is fixed at $N = 1000$ in 3D, and a single time step $\Delta t = 0.01$ is used.} Figure~\ref{prec_GMRES_basis} (a) illustrates the average condition number of the Galerkin-projected diffusion and advection operators as the dimension of the xKrylov subspaces increases across 12 outer iterations, confirming that the projected system gets more ill-conditioned as the rank increases.  Figures~\ref{prec_GMRES_basis} (b) - (d) show the GMRES residual histories for the ACS-preconditioned and unpreconditioned cases, corresponding to the second, fourth, and {twelfth} outer xKrylov iterations. These correspond to xKrylov subspace dimensions of $r = 37$, $r = 73$, and $r = 217$, respectively. As $r$ grows, unpreconditioned GMRES converges progressively slower, as expected. 
% \sout{due to high-frequency modes retained hindering the conditioning of the projected operators.}
In contrast, the ACS-preconditioner maintains rapid convergence rates independent of subspace size $r$.} 

\HHREV{To assess the relative cost of the preconditioner, we measure the ratio of the wall-clock timings for a single GMRES iteration with and without preconditioning for the different solution ranks considered in Fig.~\ref{prec_GMRES_basis}, yielding \(1.64\), \(2.09\), and \(1.85\), respectively. A preconditioned GMRES iteration is only about twice as expensive as an unpreconditioned one, but results in a dramatic improvement on GMRES' convergence rate.} %\lc{Improved this paragraph}%\lc{Is the point of this last sentence to assess the cost of the preconditioner in the GMRES solver?}\HH{Yes, this is in response to the reviewer's request/concern.} \lc{Ok, I'll improve; I don't think it conveys the point sufficiently clearly.}

\begin{figure}[h]
    \centering
    \includegraphics[width=\linewidth]{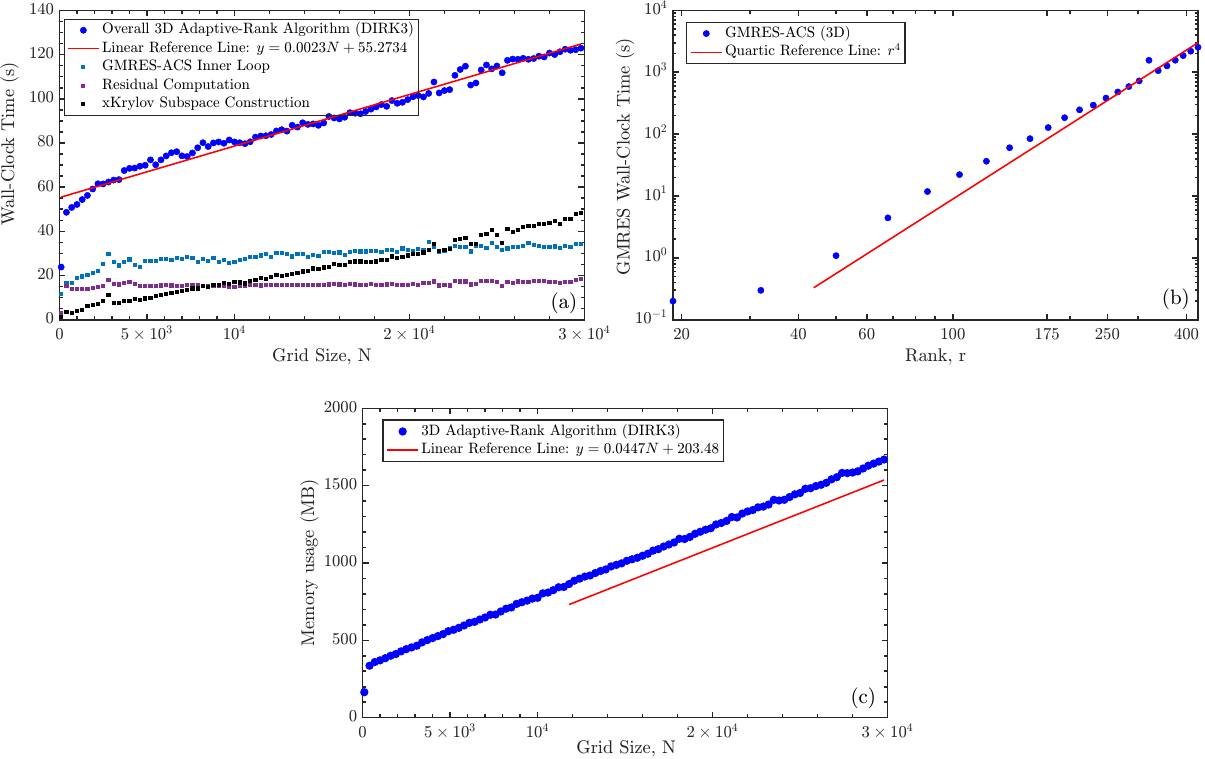}
    \caption{Example 4.1. (a) \HHREV{Profiling study depicting the computational cost of various components of the algorithm, demonstrating overall} ${\cal O}(N)$ computational complexity of DIRK3 adaptive-rank integrator in 3D. (b) Demonstration of ${\cal O}(r^{4})$ computational complexity of the GMRES-ACS solver in 3D. \HHREV{{For this test only, we allow the rank to grow indefinitely by setting the outer residual tolerance, $\epsilon_{\texttt{tol}}$, to machine precision.}
    % \HH{We should refrain from using truncation (reviewer's confusion). Can we use instead "The rank is allowed to grow even if the solution satisfies the stopping criterion."} \lc{Can we state instead that the tolerance is set to round-off?}
    (c) Demonstration of ${\cal O}(N)$ memory storage complexity of DIRK3 adaptive-rank integrator in 3D.} The execution times were recorded using MATLAB's \textbf{timeit} function. \HHREV{The memory usage was recorded using MATLAB’s \textbf{whos} command at the end of each time step, and the values reported correspond to the average memory usage across all time steps in megabytes.}
    % \HH{Remove 3D result. Revise comment on intercept.}
    % \lc{These are the only 2D results I would keep in this section.}\HH{Noted.}
    }
    \label{fig:complexity_comparison}
\end{figure}

\textbf{Verification of the Computational Complexity \HHREV{and Memory} Scalings:}
% \HH{}
We verify the $\mathcal{O}(N)$ computational complexity \HHREV{and memory storage}  claims next.
We vary $N$ from $100$ to $30000$ and the time step is fixed at \( \Delta t = 10^{-3} \). 
%Such range of mesh sizes is chosen to prevent the \( O( r^4) \) contributions of the ACS-GMRES solver from polluting the mesh-resolution scaling. 
Figure~\ref{fig:complexity_comparison}(a) displays \HHREV{a profiling study of the algorithm, depicting the computational time of various components vs. $N$ {obtained for $10$ time steps}. 
We observe that only the xKrylov component scales with mesh refinement $N$, and does so linearly. Other components do not scale with $N$, consistently with the rank of the solution being bounded.} \HHREV{Furthermore, the inner solve remains the dominant computational cost of the algorithm until, for grids of size $N \approx 20{,}000$, the xKrylov routine begins to dominate. This, in turn, highlights the importance of an efficient reduced-system solution strategy.}
% \sout{ For the same equation setup, the \(y\)-intercept, reflecting the computational cost of the reduced GSE problem, in 3D is roughly two orders of magnitude higher than in 2D.}
Figure~\ref{fig:complexity_comparison}(b) displays the computational time required for the preconditioned GMRES-ACS solver vs. rank for a single time step, using BE with \( N = 10^4 \) and \( \epsilon_{\texttt{GMRES}} = 10^{-10} \). \HHREV{For this test, 
%we did not terminate the iteration upon reaching the outer residual threshold; 
the xKrylov basis was allowed to continue to grow boundlessly, albeit with truncation applied at each augmentation step. The simulation time required for the GMRES solve was recorded at each basis augmentation.} The expected \( \mathcal{O}( r^4) \) scaling in 3D is found. {These studies verify that the complexity scaling of the approach is \( \mathcal{O}(N  r^2) +  \mathcal{O}(r^4) \) in 3D.} \HHREV{Figure~\ref{fig:complexity_comparison}(c) depicts the results of a memory scaling study demonstrating that the memory requirements also scale as $\mathcal{O}(N)$. The memory usage reported corresponds to the total size of all variables allocated within the solver function, as measured by MATLAB’s \texttt{whos} command at the end of the function call. A similar trend is expected (and observed in our simulations) in terms of memory storage: the inner tensor dominates the storage for coarser grids, while storing the xKrylov subspace dominates for finer grids.}
% \QQ{comment on inner solver dominates the computational time when $N$ is not too big; how about add profiling curves for memory?} \HH{Add line on profioing for memory.}
% , indicating the need for more efficient direct solvers for the Sylvester equation in 3D.}\QQ{what does the last sentence mean?}\HH{How about now?}

\HHREVONE{
To provide additional perspective on our results, we compare the performance of the adaptive-rank integrator for this example with a hypothetical optimal full-rank multigrid solver, which at best would scale as $\mathcal{O}(N^d \log N)$~\cite{jones1997parallel}. Assuming that the low-rank and multigrid solvers are at parity for $N=100$ (Appendix~\ref{Appendix_MG}), Table~\ref{tab:wct_ratio} reports the theoretical wall-clock-time speedup of the adaptive-rank integrator vs. the full-rank multigrid solver. The table suggests vast efficiency gains from our approach vs. optimal multigrid for this problem, scaling as 1D with grid refinement and effectively mitigating the curse of dimensionality. %\QQ{this is to compare with the computing time reported in Figure 5(a) with line $0.0023N+..$?}
}%\QQ{MG->multigrid?}\lc{Yes! Changed}
\begin{table}[t!]
\color{black}
\centering
\caption{Theoretical speedup of the adaptive-rank solver versus a hypothetical optimal multigrid solver for various 1D grid resolutions $N$ in a 3D domain.}
\label{tab:wct_ratio}
\begin{tabular}{c|c}
\hline
1D grid size $N$ & Wall-clock-time speedup in 3D \\
\hline
$1000$  & $1446.1$  \\
$4000$  & $99229.3$  \\
$10000$ & $1418193.2$  \\
\hline
\end{tabular}
\end{table}
% }
% The memory usage reported in this study corresponds to the total size of all variables allocated within the solver function. Specifically, memory was recorded using MATLAB’s \texttt{whos} command at the end of the function call. This provides an estimate of the memory footprint of all variables instantiated locally within the function, excluding overhead from the MATLAB runtime or shared libraries. The reported values reflect the sum of the \texttt{.bytes} field from \texttt{whos}, converted to megabytes. This approach allows us to isolate the algorithmic memory consumption, independently of global workspace or system-level allocations.

% \HH{I moved this sentence to the end of the study, as the computational complexity result depends on all the components discussed above. Dr. Qiu, does this address your comment below?}
% \sout{From the above study, we verify that the predominant complexity of our approach for this example is \( \mathcal{O}(N  r^2 +  r^3) \) in 2D and \( \mathcal{O}(N  r^2 +  r^4) \) in 3D, effectively reducing the complexity of full rank methods, yet maintaining their high-order convergence.}
% \QQ{assuming that iteration number is independent of $N$ and $r$...}

 %%%%%%

\subsection*{\textbf{Example 4.2:} 3D Swirling Deformation Flow}

 We consider next a 3D swirling deformation flow:
 % \sout{advection-diffusion equation}:
\begin{equation}
    f_t + \left[- 4 (1 - x^2)   y   z   f\right]_x + \left[-4 x   (1 - y^2)   z  f\right]_y + \left[8 x   y   (1-z^2)   f\right]_z = \nu \Delta f, \quad x, y, z \in [-1, 1].
    \label{swirling}
\end{equation}
% \lc{Is this flow also solenoidal?}\HH{Yes, I have just checked the code, and revised the coefficients accordingly.}\lc{Revised how? Did the coefficients change? If so, have the simulations changed?}
We set the diffusion coefficient \(\nu = 1/500\). \HHREV{The flow is solenoidal, and therefore a pure-vorticity flow.}
The initial condition is a bi-Gaussian function given by:
\begin{align*}
    f(x, y, 0) &= 0.5 \exp\left(-\frac{1}{2 \sigma^2} \left[(x - 0.3)^2 + (y - 0.35)^2 + (z - 0.3)^2\right]\right) \\ &+ 0.8 \exp\left(-\frac{1}{2 \sigma^2} \left[(x + 0.5)^2 + (y + 0.5)^2 + (z + 0.5)^2\right]\right),
\end{align*}
where \(\sigma = 0.15\). We choose \(N_x = N_y = N_z = 1000\). 
% \lc{Would be good to explain why we choose such a bizarre viscosity value}\HH{Fixed.}
The outer basis residual tolerance is \(\epsilon_{\texttt{tol}} =  10^{-5}\), the basis truncation threshold is \(\epsilon_{\kappa} = 10^{-5}\), the SVD truncation threshold is \(\epsilon = 10^{-5}\), and the GMRES acceptance threshold is $\epsilon_{\texttt{GMRES}}=10^{-10}$. 

\begin{figure}[h]
    % \centering    \includegraphics[width=0.7\linewidth]{rank_evolve.pdf}
    \centering    \includegraphics[width=\linewidth]{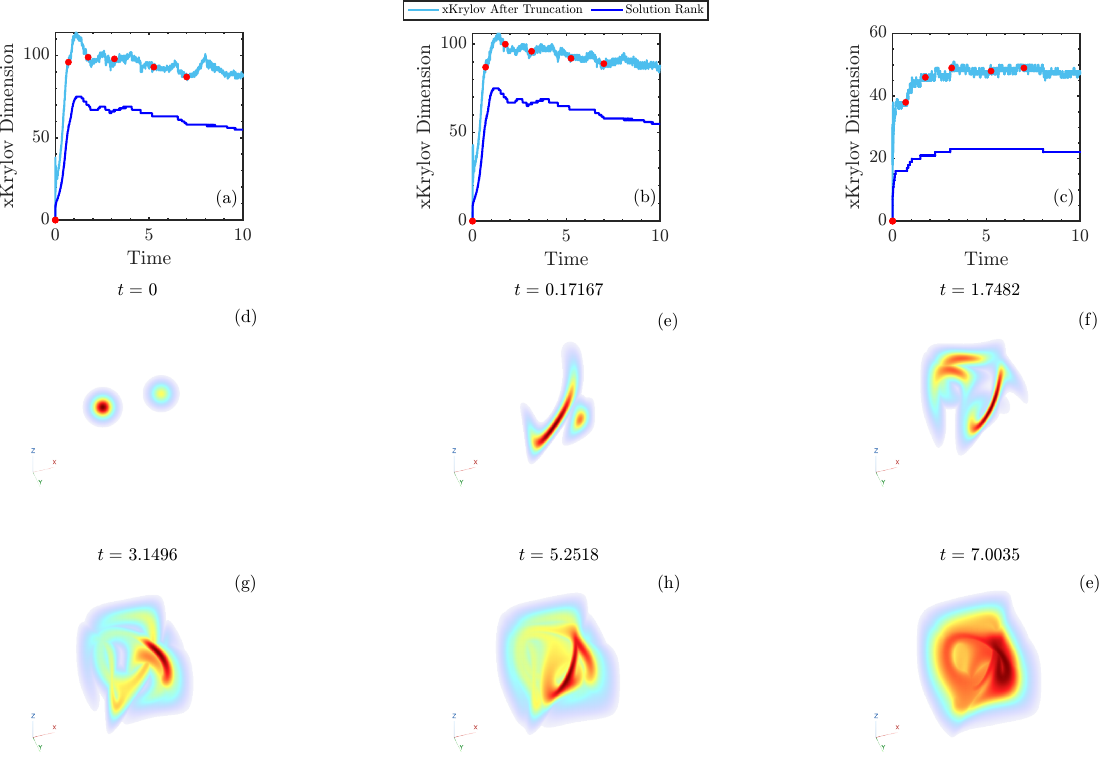}
   \caption{Example 4.2. 3D Swirling deformation problem using DIRK3 with time step corresponding to $\lambda_A=14$ on a \(1000 \times 1000 \times 1000\) grid. (a)--(c) Evolution of rank, xKrylov subspace before and after truncation, and number of xKrylov iterations over time for the $x$, $y$, and $z$ directions, respectively. (d)--(e) Solution profiles at \(t = 0\), \(t = 0.17167\), \(t = 1.7482\), \(t = 3.1496\), \(t = 5.2518\), and \(t = 7.0035\). The red dots in (a)--(c) mark these time instants and correspond to the profiles shown in (d)--(e).
   % \QQ{comment on red dots in (a).}\HH{Done.}
   The final simulation time is \(T_f = 10\).
   }\label{fig:Swirling_prob_2}
\end{figure}

\begin{figure}[h]
    % \centering    \includegraphics[width=0.7\linewidth]{rank_evolve.pdf}
    \centering    \includegraphics[width=\linewidth]{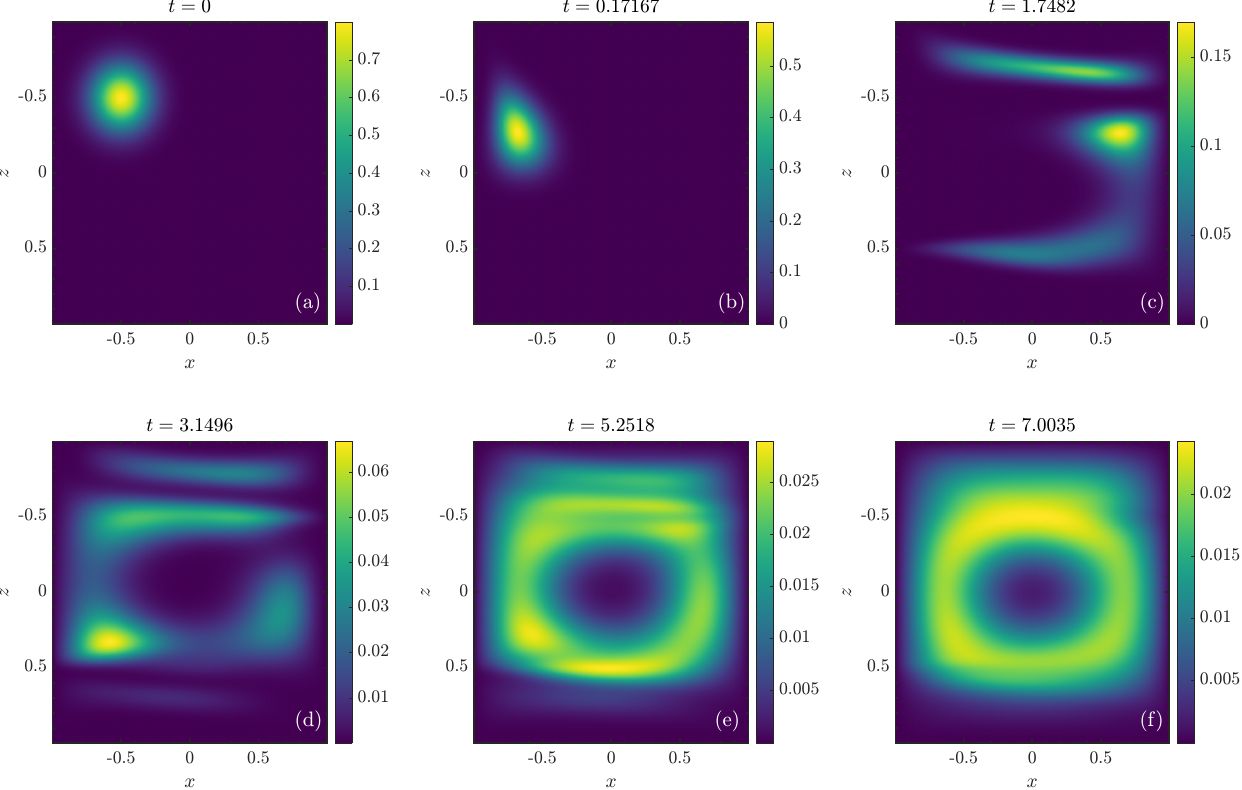}
   \caption{Example 4.2. 3D Swirling deformation problem using DIRK3 with time step corresponding to $\lambda_A=14$ on a \(1000 \times 1000 \times 1000\) grid. (a)--(f)  snapshots of the solution in the $x\text{--}z$ plane at \(y = -0.5\), corresponding to times \(t = 0\), \(t = 0.17167\), \(t = 1.7482\), \(t = 3.1496\),  \(t = 5.2518\), and \(t = 7.0035\), respectively.}\label{fig:Swirling_snapshots}
\end{figure}

This example demonstrates the algorithm's ability to adapt to dynamic rank evolution. The algorithm performs nine {xKrylov} augmentations at the initial time step, two augmentations at the third and seventh time steps, and a single augmentation otherwise. Figure~\ref{fig:Swirling_prob_2} shows rank evolution. The solution rank rises sharply in the $x$- and $y$-directions [Figures~\ref{fig:Swirling_prob_2} (a)-(b)].
% \sout{As shown in Figure~\ref{fig:Swirling_prob_2} (a)-(c),}
The solution rank in the $z$-direction evolves more gently [Figure~\ref{fig:Swirling_prob_2} (c)] and is roughly one-third of that in the $x$ and $y$ directions. 
% \QQ{It is not clear to me where the factor of $4.5$ and $3$ come from? From the figures (a)-(c)?} \lc{These factors are stated, but are not shown in any figure. We used to have the data in Fig. 6(a-c) along with what is currently shown, but we removed them because it led to confusion for the referee as to what the actual Krylov subspace sizes were. I have attempted to fix the discussion to clarify this further.} 
% \QQ{I suggest to move the above comments later, after other subfigures are discussed. Maybe we can even skip these discussions here...I also suggest to move the discussion in conclusion part here. maybe as a separate paragraph.}
% \sout{Consequently, the dimensions of the xKrylov subspaces, both before and after truncation, are approximately three times larger in the $x$- and $y$-directions compared to the $z$-direction. This demonstrates the algorithm's capability to effectively construct and adaptively truncate subspaces according to the evolving solution dynamics.}
{Figures~\ref{fig:Swirling_prob_2} (d)-(e) showcase how the solution deforms according to the prescribed advection coefficients and diffuses slowly over time.} {Figure~\ref{fig:Swirling_snapshots} shows snapshots of the solution in the \(x\text{--}z\) plane at \(y = -0.5\)}. {As the flow-driven deformation increases the complexity of the solution structure, the rank grows to capture these features; when the diffusion process starts to dominate, the rank decreases as the solution smooths.} 
%\HHREV{A snapshot at the solution time reveals that the \lc{Missing text?}}  
\HHONE{This example shows that the low-rank assumption is not intrinsic to the algorithm. Even when the rank becomes moderately large, the algorithm remains capable of  capturing the underlying dynamics effectively. Naturally, the benefit of a low-rank approximation decreases as the solution approaches a full-rank structure.}

\HHREV{From a computational perspective, the per-time-step cost breakdown is as follows: xKrylov construction and orthogonalization (10\%), ACS-GMRES inner solve (47.5\%), and residual evaluation (42.5\%). The relatively small contribution of the cost of xKrylov construction is partly a consequence of the truncation mechanism, together with the effectiveness of the chosen projection spaces, which collectively keep the approximation spaces compact. On average, in the $x$-$y$ directions, the truncation reduces the subspace dimension by a factor of about $4.5$, from roughly $450$ to $100$. In the $z$ direction, the reduction amounts to a factor of approximately $3$, from about $150$ to $50$. On the other hand, the relatively large contribution of the residual evaluation to the simulation time originates in the poor scaling of the computational complexity with the number of GSE terms (eight in this example, including the RHS) and suggests that a more economical stopping criterion (as discussed in Section~\ref{residual}) should be used.}

% }

% , \HHREV{and, on average, the truncation mechanism reduces the xKrylov subspace size by a factor of approximately $4.5$\lc{, i.e, from a dimension of about $450$ to about $100$, as reported in the figure}}.

% \HHREV{and, on average, the truncation mechanism reduces the xKrylov subspace size by a factor of approximately $3$ \lc{(from about $150$ to $50$, as reported)}}

% \QQ{Figure 7 is not commented...}\lc{Good catch! Maybe we can remove it, then? Does not really add much.}
% \HHONE{\sout{We consider next a time-dependent advection-diffusion equation with a weak diffusion coefficient. Our goal is to demonstrate that the assumption of low-rank structure is not intrinsic to our proposed algorithm. As will be shown, even when the rank moderately increases, the algorithm remains capable of effectively capturing the underlying dynamics. 
% Naturally, the benefit of using low-rank approximations decreases as the solution approaches a full-rank structure.}}\HH{This sentence is response to reviewer.}\QQ{I move this part to here...}\lc{Looks good.}

 %{LC: not particularly interesting}We observe that the Krylov subspace dimensions and the rank in the \( z \)-direction are approximately half the corresponding dimensions in the \( x \)- and \( y \)-directions, as dictated by the prescribed advection coefficients.

\subsection*{\textbf{Example 4.3:} A 2D Steady-State Balanced-Advection-Diffusion {Problem}}
%\HH{Hi Luis, just a quick note to let you know I have updated the Steady-State example figure to reflect the new version of the algorithm for completeness. I have also corrected some of the simulation details correspondingly. }
%\QQ{"Direct-to"?}
The last example prescribes diffusion and advection coefficients that exactly balance each other, resulting in a {non-trivial} steady-state solution. {This test demonstrates the ability of the algorithm to take large timesteps stably to solve directly for the steady state with expected spatial asymptotic accuracy.}

We seek an equilibrium solution satisfying the steady-state condition:
\begin{equation}
    0 = \nabla \cdot \left(\mathbfit{\Phi} \cdot \nabla f_{\text{eq}} - \mathbfit{\Sigma} f_{\text{eq}}\right),
    \label{eq:adv-diff_constructed}
\end{equation}
a solution of which is:
\begin{equation*}
    \mathbfit{\Sigma} = \mathbfit{\Phi} \cdot \nabla \log(|f_{\text{eq}}|).
\end{equation*}
This equation defines the steady-state  advective flow $\mathbfit{\Sigma}$ for given diffusion coefficient $\mathbfit{\Phi}$ and steady-state solution $f_{\text{eq}}$. We construct a steady-state solution in the domain $ 0 \leq x,y \leq 1$ to be (up to an arbitrary constant):
\begin{equation*}
    f_{\text{eq}}(x, y) = \left( x(1 - x) \right)^2 \left( y(1 - y) \right)^2,
\end{equation*}
which satisfies homogeneous Dirichlet boundary conditions. The diffusion coefficients are defined as: 
\begin{align*}
    \phi^{x}(x,y) = \phi^{y}(x,y) &= \phi^1(x) \phi^2(y), \qquad
    \phi^1(x) = x^2 (1 - x)^2, \qquad
    \phi^2(y) = y^2 (1 - y)^2.
\end{align*}
The variable advection coefficients are 
$\sigma^x(x,y) = \sigma_1^1(x) \sigma_1^2(y)$, $\sigma^y(x,y) = \sigma_2^1(x) \sigma_2^2(y),$
where:
\begin{align*}
    \sigma_1^1(x) &= 2x(1 - 3x + 2x^2), \qquad
    \sigma_1^2(y) = y^2 (1 - y)^2, \qquad
    \sigma_2^1(x) = x^2 (1 - x)^2, \qquad
    \sigma_2^2(y) = 2y(1 - 3y + 2y^2).
\end{align*}
We select the following rank-1 initial condition: 
$f(x, y,0) = \left|\sin(2\pi x) \sin(2\pi y)\right|$. 
    %, \quad 0 \leq x, y \leq 1.
%Since the solution will converge to a steady state up to a multiplicative constant, we scale the mass of the initial condition to match that of the steady-state solution. 
The outer basis residual tolerance is \(\epsilon_{\texttt{tol}} =  10^{-3}\), the basis truncation threshold is \(\epsilon_{\kappa} = 10^{-8}\), the SVD truncation threshold is \(\epsilon = 10^{-8}\), and the GMRES acceptance threshold is $\epsilon_{\texttt{GMRES}}=10^{-8}$. 
% \HH{To correct the mass of the solution, we scale it by a constant to match the steady-state mass.}
% \sout{To correct the mass of the solution, we adopt the LoMaC strategy outlined in \cite{guo2022conservative,guo2022low,guo2023local}, as implemented in \cite{el2024krylov}.} \lc{The boundary conditions are homogeneous Dirichlet, which do NOT allow a null space. So the previous sentence is very confusing...}

\begin{figure}[h]
    \centering    \includegraphics[width=0.6\linewidth]{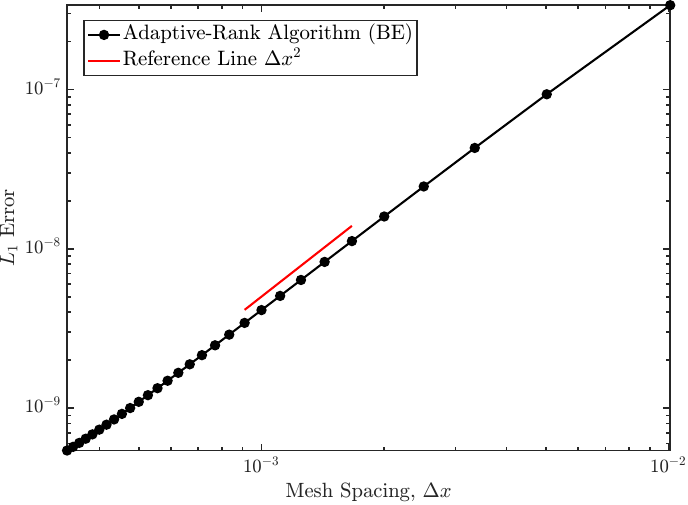}
    \caption{Example 4.3. Log-log plot of error in the numerical steady-state solution compared to the prescribed analytical steady state as a function of mesh size. The advection CFL number, \(\lambda_A\), ranges from $1,190$ to $36,076$. The plot shows second-order spatial convergence in the steady state.  }
    \label{fig:steady_state_error}
\end{figure}

Figure~\ref{fig:steady_state_error} shows the error of the numerical steady-state solution relative to the prescribed analytical steady state as a function of mesh size. The solution is obtained using ten time steps of size $\Delta t =1000$. 
% \lc{How many timesteps do we take to reach steady state?}\HH{One time step. Added detail.}
The advection CFL number \(\lambda_A\)  ranges from $1,190$ to $36,076$ as the mesh is refined. 
% \WT{It's impressive being able to take such large time-steps. If it doesn't take any wind off our sail, perhaps we should also show the number of outer xKrylov basis size and inner GMRES iteration side by side with the error plot?}
The error is obtained by comparing the numerical and analytical solutions, with the numerical solution rescaled so that its normalization constant is the same as the analytical one. The results confirm second-order spatial convergence and the ability of the algorithm to solve directly for the correct steady state. 
% \QQ{computing time and memory profiling as in the first example?} \lc{Not sure the steady-state problem will add much to the conversation; perhaps we can just comment that the CPU scaling and memory scaling is found to agree with Ex. 4.1 (if this is indeed the case...)}%\HH{Lead with this sentence.}

\section{Conclusion}
\label{Section_5}

{We propose} a computationally efficient adaptive-rank implicit algorithm for solving time-dependent advection-diffusion partial differential equations with variable coefficients. 
% \sout{We discretize the PDE using second-order centered finite differences in space and implicitly in time with up to third-order accuracy using DIRK, resulting in a generalized Sylvester matrix equation.}
% \st{As proposed elsewhere in the literature, we formulate the advection-diffusion PDE as a Generalized Sylvester Equation.
% To invert it, we follow a projection-type approach for GSEs, namely:} 
 \HHREV{We formulate the corresponding implicit system as a GSE, for which we follow a projection type approach, namely:
 % \lc{I am OK with this}
constructing dimension-wise basis from xKrylov subspaces, projecting the original matrix equation onto a {reduced system}, and
solving for {the coefficient matrix}. Our approach innovates in the way we construct the Krylov subspaces, and in the solution of the reduced system, which leverages GMRES effectively preconditioned with an effective Averaged-Coefficient Sylvester (ACS) approximation.
The approach achieves a 1D-like computational complexity (of $\mathcal{O}(N r^2) + \mathcal{O}(r^{d+1})$, with $r$ being the rank of the solution during Krylov iteration process and $d$ the dimensionality) and memory storage [of $\mathcal{O}(N r) + \mathcal{O}(r^{d})$]}, effectively mitigating the curse of dimensionality when the solution rank is bounded and making it comparable in efficiency to the constant-coefficient case \cite{el2024krylov}. 
% \HH{Although the rank is an inherent property of the solution, it has been shown theoretically \cite{kressner2010krylov} that, in principle, $r$ may depend on the grid size $N$ when the chosen approximation subspaces are extended Krylov. However, this scaling is generally pessimistic for a wide range of grid resolutions and does not account for the smoothness of the solution, which we demonstrate numerically in this work.}
\HHONE{Such complexity scaling is numerically verified with extensive tests, with the maximal rank $r$ found to be independent of the 1D grid size $N$ {for practical grid resolutions and convergence tolerances}.} Future research directions include extending the algorithm to higher dimensions and nonlinear PDEs. 

% \HHREV{As advanced in earlier studies~\cite{palitta2025subspace}, our computational complexity and memory storage analysis confirms that the residual computation can become very expensive when the number of terms in the GSE grow large even if the rank remains bounded. This motivates the use of cheaper alternatives for assessing convergence, such as monitoring the difference between consecutive solution iterates~\cite{palitta2025subspace}.}\QQ{this part of discussion could be earlier on in the numerical example; but does not have to be in the conclusion?} \lc{Since we did not do this in this study, this does belong to future work... but we can remove from here if you'd like. Not sure we need to add it in the previous section, since we already make this point. }
%{We validated our algorithm through a series of numerical experiments, demonstrating its scalability and efficiency in capturing time-dependent solutions.}
% \sout{We tested our algorithm through a series of numerical experiments, showcasing. Overall, our adaptive-rank algorithm offers a scalable and efficient solution for time-dependent advection-diffusion type PDEs with variable coefficients. It effectively balances computational cost and accuracy by efficiently constructing and truncating the Krylov basis and dynamically adjusting the solution rank. The ability to handle large meshes and time steps without sacrificing accuracy makes it a promising tool for a wide range of applications.} 
%\lc{I like these changes}

\section*{Acknowledgements}

This work was partially supported by the Multifaceted Mathematics Integrated Capability Centers (MMICCs) program
of DOE Office of Applied Scientific Computing Research (ASCR). Los
Alamos National Laboratory (LANL) is operated by Triad National Security,
LLC, for the National Nuclear Security Administration of U.S. Department
of Energy (Contract No. 89233218CNA000001). 
J.Q. and H.E. were supported by NSF grant NSF-DMS-2111253. J.Q. was also supported by Air Force Office of Scientific Research FA9550-22-1-0390, Department of Energy DE-SC0023164 and Air Force Office of Scientific Research (AFOSR) FA9550-24-1-0254 via the Multidisciplinary University Research Initiatives (MURI) Program. H.E. was also supported by MMICCs at LANL during the summer of 2024.

\appendix
\section{Discretization of the Diffusion Operator}
\label{sec:discret_diff}

We derive the matrix formulation of the discretization of the diffusion operator for the two-dimensional advection-diffusion equation with separable diffusion coefficients. We consider the discretization of the diffusion term in the $x$-direction, with other directions treated similarly:
\begin{equation}
    \frac{\partial}{\partial x} \left( \phi^{x}(x, y) \frac{\partial u(x, y)}{\partial x} \right),
\end{equation}
where the diffusion coefficient is separable:
\begin{equation}
    \phi^{x}(x, y) = \phi^{1}(x) \phi^{2}(y).
\end{equation}
We discretize this term using a 3-point stencil on a uniform grid with spacing $\Delta x$ in the $x$-direction. Let $u_{i,j} \approx u(x_i, y_j)$, with $i=0,\cdots,N_x-1$ and $j=0,\cdots,N_y-1$.
We approximate the spatial derivative using central differences. The discrete approximation at grid point $(i, j)$ is:
\begin{equation}
    \left[ \frac{\partial}{\partial x} \left( \phi^{x} \frac{\partial u}{\partial x} \right) \right]_{i,j} \approx \frac{1}{\Delta x^2} \left[ \phi^{1}_{i+\frac{1}{2}} (u_{i+1, j} - u_{i,j}) - \phi^{1}_{i-\frac{1}{2}} (u_{i,j} - u_{i-1, j}) \right]\phi^{2}_j,
    \label{eq:diffusion_discrete}
\end{equation}
with
$\phi^{1}_{i+\frac{1}{2}} \phi^{2}_j=
\phi^{1}(x_{i+\frac{1}{2}})
\phi^{2}(y_j).$
To express the discretization in matrix form, we define:
\begin{itemize}
    \item $\mathbf{F}$: matrix of solution values, $[\mathbf{F}]_{i,j} = f_{i,j}$.
    \item $\Phi^{1}_{+}$: diagonal matrix with elements $\left[ \Phi^{1}_{+} \right]_{i,i} = \phi^{1}_{i+\frac{1}{2}}$. 
    \item $\Phi^{1}_{-}$: diagonal matrix with elements $\left[ \Phi^{1}_{-} \right]_{i,i} = \phi^{1}_{i-\frac{1}{2}}$.
    \item $\Phi_{2}$: diagonal matrix with elements $\left[ \Phi_{2} \right]_{j,j} = \phi^{2}_j$.
\end{itemize}
We also define the forward and backward difference matrices in the $x$-direction for the interior nodes $i=1,\cdots,N_x-2$. The nodes $i=0$, $i=N_x-1$ rest at the boundary and enforce homogeneous Dirichlet boundary conditions (other boundary conditions can be readily considered):
\begin{align}
    D_{+} = \frac{1}{\Delta x}
    \begin{pmatrix}
    -1 & 1 & 0 & \cdots & 0 \\
    0 & -1 & 1 & \cdots & 0 \\
    \vdots & \ddots & \ddots & \ddots & \vdots \\
    0 & \cdots & 0 & -1 & 1 \\
    0 & \cdots & 0 & 0 & -1
    \end{pmatrix},
    D_{-} = \frac{1}{\Delta x}
    \begin{pmatrix}
    1 & 0 & 0 & \cdots & 0 \\
    -1 & 1 & 0 & \cdots & 0 \\
    0 & -1 & 1 & \cdots & 0 \\
    \vdots & \ddots & \ddots & \ddots & \vdots \\
    0 & \cdots & 0 & -1 & 1
    \end{pmatrix}.
\end{align}
Using these definitions, the discrete diffusion term in Eq.~\eqref{eq:diffusion_discrete} can be represented in matrix form as:
\begin{equation}
     \left( \Phi^{1}_{+} D_{+} - \Phi^{1}_{-} D_{-} \right) \mathbf{F}{\Phi_{2}^\top}.
    \label{eq:diffusion_matrix_form}
\end{equation}
We note that the operator multiplying $\mathbf{F}$ from the left is analogous to $D^{\Phi}_{xx}$ in \eqref{MOL_eq}. The discretization in the $y$-direction follows the same strategy.

\section{Discretization of the Advection Operator}
\label{sec:discret_adv}
 
We derive the matrix formulation of the discretization of the advection operator for the two-dimensional advection-diffusion equation with variable, separable advection coefficients. As before, we focus on the discretization of the advection term in the $x$-direction:
\begin{equation}
    \frac{\partial}{\partial x} \left( \sigma^x(x, y) u(x, y) \right),
\end{equation}
where the advection coefficient is separable:
\begin{equation}
    \sigma^x(x, y) = \sigma^{1}(x) \sigma^{2}(y).
\end{equation}
We discretize this term on a uniform grid with spacing $\Delta x$ in the $x$-direction and $\Delta y$ in the $y$-direction. 
To preserve linearity of the discretization without introducing too much numerical dissipation, we employ a central difference scheme to approximate the spatial derivative. This discretization is monotonicity-preserving when the following condition is satisfied:
\begin{equation}
    \Delta x < 2 \frac{\phi_{\text{max}}}{\sigma_{\text{max}}},
\end{equation}
where $\sigma_{\text{max}}$ is the maximum magnitude of the advection velocity and $\phi_{\text{max}}$ is the maximum diffusion coefficient. This condition is practical in our formulation due to the very favorable scaling of the low-rank representation with mesh size and rank.

The discrete approximation at the grid point $(i, j)$ is:
\begin{align}
     \left[\frac{\partial}{\partial x} \left( \sigma^x(x, y) u(x, y) \right)\right]_{i,j} &\approx  \frac{1}{2 \Delta x} \left[ \sigma^{1}_{i+\frac{1}{2}} (u_{i+1, j} + u_{i, j}) - \sigma^{1}_{i-\frac{1}{2}} (u_{i, j} + u_{i-1, j}) \right]\sigma^{2}_j,
\end{align}
where $\sigma^{1}_{i+\frac{1}{2}} = \sigma^{1}(x_{i+\frac{1}{2}})$ is the advection coefficient evaluated at the midpoint between $x_i$ and $x_{i+1}$, and similarly for $\sigma^{1}_{i-\frac{1}{2}}$, and $\sigma^2_j$ is evaluated at $y_j$. 
To express the discretization in matrix form, we define:
\begin{itemize}
    \item $\mathbf{F}$: matrix of solution values, $[\mathbf{F}]_{i,j} = f_{i,j}$.
    \item $\Sigma^{1}_{+}$ and $\Sigma^{1}_{-}$: diagonal matrices with elements $\left[ \Sigma^{1}_{+} \right]_{i,i} = \sigma^{1}_{i+\frac{1}{2}}$ and $\left[ \Sigma^{1}_{-} \right]_{i,i} = \sigma^{1}_{i-\frac{1}{2}}$.
    \item $\Sigma_{2}$: diagonal matrix with elements $\left[ \Sigma_{2} \right]_{j,j} = \sigma^{2}_j$.
\end{itemize}
We also define the difference matrices in the $x$-direction: 
\begin{align}
    S_{+} = \frac{1}{2}
    \begin{pmatrix}
    1 & 1 & 0 & \cdots & 0 \\
    0 & 1 & 1 & \cdots & 0 \\
    \vdots & \ddots & \ddots & \ddots & \vdots \\
    0 & \cdots & 0 & 1 & 1 \\
    0 & \cdots & 0 & 0 & 1
    \end{pmatrix}, 
    S_{-} = \frac{1}{2}
    \begin{pmatrix}
    1 & 0 & 0 & \cdots & 0 \\
    1 & 1 & 0 & \cdots & 0 \\
    0 & 1 & 1 & \cdots & 0 \\
    \vdots & \ddots & \ddots & \ddots & \vdots \\
    0 & \cdots & 0 & 1 & 1
    \end{pmatrix}.
\end{align}
Using these definitions, the discrete advection term in the $x$-direction can be represented as:
\begin{equation}
    \frac{1}{ \Delta x} \left( \Sigma^{1}_{+} S_{+} - \Sigma^{1}_{-} S_{-} \right) \mathbf{F} {\Sigma_{2}^\top}.
    \label{eq:advection_matrix_form_x}
\end{equation}
We note that the operator multiplying $\mathbf{F}$ from the left is analogous to $D^{\Sigma}_{x}$ in Equation \eqref{MOL_eq}. The discretization in the $y$-direction follows the same strategy.

\section{High-Order Temporal Discretization }\label{Appendix_DIRK}
We extend the backward Euler implicit integrator to higher-order DIRK time discretizations for the advection-diffusion operator, as compactly written in \eqref{MOL_operator_L}.
\subsection{DIRK Scheme}\label{DIRK_app_1}
DIRK methods are characterized by the Butcher tableau, Table \ref{table:DIRK}.
\begin{table}[t!]
\centering
\caption{Butcher tableau for an $s$-stage DIRK scheme. Here, $s$ represents the number of stages in DIRK, $c_i$ represents the intermediate stage, $a_{ij}$ is a lower-triangular matrix with coefficients used to approximate the solution at each stage, and $b_j$ are the quadrature weights used to update the solution in the final step.}
\label{table:DIRK}
\begin{tabular}{c|llll}
    $c_1$ & $a_{11}$ & 0 & $\hdots$ & 0 \\
    $c_2$ & $a_{21}$ & $a_{22}$ & $\hdots$ & 0 \\
    $\vdots$ & $\vdots$ & $\vdots$ & $\ddots$ & $\vdots$ \\
    $c_s$ & $a_{s1}$ & $a_{s2}$ & $\hdots$ & $a_{ss}$ \\
    \hline
    & $b_1$ & $b_2$ & $\hdots$ & $b_s$
\end{tabular}
\end{table}
We focus on stiffly accurate DIRK methods, where the coefficients satisfy \(b_i = a_{si}\) for \(i = 1, \ldots, s\), ensuring that the final solution update corresponds to the solution at the last DIRK stage. The Butcher tables for DIRK2 and DIRK3 are presented in Tables \ref{Butcher:DIRK2} and \ref{Butcher:DIRK3}, respectively.
\begin{table}[h]
\centering
\label{table:ComparisonDIRK}
\begin{minipage}{.4\linewidth}
\centering
\caption{DIRK2 Butcher table, $\gamma=1-\frac{\sqrt{2}}{2}$.}
\begin{tabular}{c|cc}
 $\gamma$ & $\gamma$ & 0 \\
$1$ &  $1-\gamma$ & $\gamma$ \\
\hline
    & $1-\gamma$    & $\gamma$
\end{tabular} 
\label{Butcher:DIRK2}
\end{minipage}%
\begin{minipage}{.7\linewidth}
\centering
\caption{DIRK3 Butcher table, $x=0.4358665215$.}
\begin{tabular}{c|ccc}
$x$ & $x$ & 0 & 0 \\
$\frac{1+x}{2}$ & $\frac{1-x}{2}$ & $x$ & 0 \\
$1$ & $-\frac{3x^2}{2}+4x-\frac{1}{4}$ & $-\frac{3x^2}{2}-5x+\frac{5}{4}$ & $x$ \\
\hline
 & $-\frac{3x^2}{2}+4x-\frac{1}{4}$ & $-\frac{3x^2}{2}-5x+\frac{5}{4}$ & $x$ 
\end{tabular}
\label{Butcher:DIRK3}
\end{minipage}
\end{table}

Similarly to the backward Euler case, the matrix-based discretization of the advection-diffusion equation using the operator \(\mathscr{L}\) defined in \eqref{MOL_operator_L} results in a GSE that must be solved at each DIRK stage. The corresponding DIRK scheme for the matrix differential equation \eqref{MOL_operator_L} from \(t^{(0)}\) to \(t^{(1)}\) is written as:
\begin{subequations}
\label{eq:RK}
\begin{align}
\begin{split}
    \mathbf{F}^{(k)} &= \mathbf{F}_0 + \Delta t \sum\limits_{\ell=1}^{k} a_{k\ell} \mathbf{Y}_{\ell}, \qquad k=1,2,\ldots,s,
    \label{eq:RKstage_a}
\end{split}
\\
\begin{split}
    \mathbf{Y}_k &= \mathscr{L}(\mathbf{F}^{(k)}; t^{(k)}), \qquad t^{(k)} = t^0 + c_k \Delta t, \qquad k=1,2,\ldots,s,
\end{split}
\\
\begin{split}
    \mathbf{F}_1 &= \mathbf{F}^{(s)} = \mathbf{F}_0 + \Delta t \sum\limits_{k=1}^{s} b_k \mathbf{Y}_k.
    \label{eq:RK_finalstep}
\end{split}
\end{align}
\end{subequations}
Within each stage $k$ of the DIRK scheme, one must solve the residual equation given by:
\begin{equation}
    \bR^{(k)}_{\mathscr{L}} = \mathbf{F}^{(k)} - a_{kk}\Delta t \mathscr{L}(\mathbf{F}^{(k)}) - \mathbf{B}^{(k)}=\mathbf{0}, \quad \text{with} \quad \mathbf{B}^{(k)} = \mathbf{F}_0 + \Delta t \sum\limits_{\ell=1}^{k-1} a_{k\ell} \mathbf{Y}_{\ell}, \quad k=1,2,\ldots,s.
    \label{disc-sylv2}
\end{equation}
Expanding the operator \(\mathscr{L}\) in the residual equation above results in the following \(k\)-th stage GSE:
\begin{equation}
    \mathbf{F}^{(k)} - a_{k,k} \Delta t  \left( T_1 \mathbf{F}^{(k)} \Phi_2{^\top} + \Phi_1 \mathbf{F}^{(k)} T_2^\top + T_3 \mathbf{F}^{(k)} \Sigma_2^\top + \Sigma_1 \mathbf{F}^{(k)} T_4^\top \right) = \mathbf{B}^{(k)}, \quad k=1,2,\ldots,s.
    \label{Generalized_Sylv_DIRK}
\end{equation}
For the DIRK update, we adopt the strategy outlined in \cite{el2024krylov}, with minor modifications to improve efficiency. In particular, we: (i) construct the basis using the time-step of the first stage and retain it for subsequent stages, (ii) evolve the projected system for the $S_1$ matrix across the DIRK stages, and (iii) truncate and evaluate the residual to determine whether to {admit the solution} or return to step (i) to further augment the basis.
Key modifications vs. \cite{el2024krylov} include truncating the solution before the residual computation, and only checking the residual at the end of the DIRK step.
The reduced matrix is solved using preconditioned GMRES, as before. 
% Further details on the discretization of the advection-diffusion equation using the DIRK scheme are provided in Appendix \ref{Appendix_DIRK}. Specifically, Appendix \ref{DIRK_app_2} describes the adaptive-rank approach for high-order integration, Appendix \ref{DIRK_app_3} extends the integration to multi-rank coefficients with DIRK schemes {and presents a pseudo-algorithm, \HH{Appendix \ref{Appendix_3D} presents the 3D formulation}}. \QQ{3D Tucker???}

% As shown in Equation \eqref{Generalized_Sylv_DIRK}, the residual equations lead to a GSE to be solved at each stage. 

\subsection{Adaptive-Rank Strategy for High-Order Temporal Discretization via DIRK}\label{DIRK_app_2}

The high-order DIRK method features a stage-by-stage generalized-Sylvester-equation solve. The strategy follows closely that of the backward Euler adaptive-rank treatment, i.e., we construct an approximate $\mathscr{\Tilde{L}}$ used for both the basis construction of $U_1$ and $V_1$ and the preconditioning of the reduced system to solve $S_1$. 

At each stage $k$ of the DIRK scheme, the corresponding DIRK formulation for the matrix differential equation from $t^{(0)}$ to $t^{(1)}$ 
requires solving an SE for the approximated operator $\mathscr{\Tilde{L}}$ by solving the $k$th-stage residual:
\begin{equation}
\mathbf{R}^{(k)}_{\tilde{\mathscr{L}}} = \mathbf{F}_1 - a_{kk}\Delta t \, \tilde{\mathscr{L}}(\mathbf{F}_1) - \mathbf{F}_0 - \mathbf{B}^{(k)} = \mathbf{0}, \quad \text{with}  \quad    \mathbf{B}^{(k)} = \mathbf{F}_0 + \Delta t \sum\limits_{\ell=1}^{k-1} a_{k\ell} \mathbf{Y}_{\ell}.
\label{disc-sylv}
\end{equation}
Expanding this residual using $\tilde{\mathscr{L}}$, defined in \eqref{approx_operator_L}, leads to the $k$th-stage SE:

\begin{equation}
\underbrace{\left[\underbrace{\left(\frac{1}{4}I - a_{kk}\Delta t \, \alpha_2 T_1\right)}_{A_1} + \underbrace{\left(\frac{1}{4}I - a_{kk}\Delta t \, \gamma_2 T_3\right)}_{A_3}\right] \mathbf{F}^{(k)}}_{P_1}
    + \mathbf{F}^{(k)} \underbrace{\left[\underbrace{\left(\frac{1}{4}I - a_{kk}\Delta t \, \alpha_1 T_2\right)^\top}_{A^\top_2} 
    + \underbrace{\left(\frac{1}{4}I - a_{kk}\Delta t \, \gamma_1 T_4\right)^\top}_{A^\top_4}\right]}_{P_2^T} = \mathbf{B}^{(k)}.
\label{Approximated_Sylv_DIRK}
\end{equation}
Note that the operators $A_{1:4}$ are fixed throughout the stages since the DIRK method uses constant diagonal Butcher tableau coefficients.

We follow \cite{el2024krylov} to perform high-order integration. The key idea is to consider the basis fixed across DIRK stages, and only evolve the matrix of coefficients ${\bS}_1$ through the stages. The procedure is as follows:

\begin{enumerate}
    \item [Step 1.] Prediction of Krylov basis functions: construct a set of orthonormal bases $U_1$ and $V_1$ from the Krylov-based low-rank implicit solver at the first DIRK stage, i.e. backward Euler with time stepping size $c_1 \Delta t$.
    
    \item [Step 2.] For $k=1:s$ (per DIRK stage)
    \begin{enumerate}
        \item Solve reduced SE for $\bS^{(k)}$ using preconditioned GMRES,
    \begin{equation}
     \bS^{(k)} - \Delta t \left( \tilde{T}_1 \bS^{(k)} \tilde{\Phi}^\top_2 + \tilde{\Phi}_1 \bS^{(k)} \tilde{T}_2^\top + \tilde{T}_3 \bS^{(k)} \tilde{\Sigma}^\top_2 + \tilde{\Sigma}_1 \bS^{(k)} \tilde{T}_4^\top \right) = \tilde{\bB}^{(k)},
        \label{eq: Sk}
    \end{equation}
    with \[
    \tilde{\bB}^{(k)} =  U_1^\top \mathbf{F}_0 V_1 +  \Delta t\sum\limits_{\ell=1}^{k-1}{a_{k\ell}\tilde{Y_{\ell}}}, 
    \]
   Here $\tilde{Y_{\ell}} = U_1^\top {Y}_{\ell} V_1 \in \mathbb{R}^{r_x \times r_y}$. The preconditioner used is given by 
   \begin{equation}
       \bS\mapsto \Tilde{P}_1\bS+\bS\tilde{P}_2,
       \label{preconditioner_DIRK}
   \end{equation}
    where $\Tilde{P}_1=U_1^TP_1U_1$ and $\Tilde{P}_2=V_1^TP_1V_1$, and is valid for all stages since the basis constructed is fixed throught the DIRK stages by assumption. Solve for $\bS^{(k)}$ and obtain the intermediate RK solutions as $U_1 S^{(k)} V_1^\top$. 
To further improve computational efficiency, from \eqref{eq: Sk}, we have,
\begin{align*}
    \tilde{Y}_{\ell} &= \frac{1}{a_{\ell \ell}\Delta t}\left( {\bS^{(\ell)}} -  {\tilde{\bB}^{(\ell)}}\right),
\end{align*}
leading to an efficient computation of $\tilde{\bB}^{(k)}$:
\begin{equation}
     \tilde{\bB}^{(k)} =  \tilde{\bB}_1 + \Delta t \sum_{\ell=1}^{k-1} \frac{a_{k\ell}}{a_{\ell \ell}} \left( {\bS^{(\ell)}} -  {\tilde{\bB}^{(\ell)}}\right),
\end{equation}
    with $\tilde{\bB}_1$ defined in \eqref{Galerkin_S}. This avoids the need to evaluate the full size $\bB^{(k)}$ in \eqref{disc-sylv}. 
    
\end{enumerate} 
\item [Step 3.]
 Perform truncation of evolved solution, i.e.  $\{U_{1,\epsilon}, S_{1,\epsilon}, V_{1,\epsilon}^\top\} \leftarrow\mathcal{T}_{\epsilon}(\{U_1, S_1, V_1^\top \})$, with $S_1:=S_1^{(s)}$.\
    \item [Step 4.] Efficiently evaluate the residual norm, 
    \[
 \left\|{\bR^{(s)}_{\mathscr{L}}} \right\|= \left\|  {R}_U
{\text{diag}\left(- \tilde{B}^{(s)}_1, S_{1,\epsilon}, -a_{ss}\Delta t \left( I_{{\mathcal{R}}_{\text{ranks}}} \otimes S_{1,\epsilon}\right) \right)}
 {R}_V^\top \right\|,
\]
where $R_U$ and $R_V$ are upper triangular matrices from a reduced Q-less \texttt{QR} decomposition as described in \eqref{RURV}, which can be done once for all DIRK stages. 
\item [Step 5.] Compare the computed residual norm $\|\bR^{(s)}_{\mathscr{L}}\|$ with the given error tolerance $\epsilon_\texttt{tol}$. If the tolerance is met, then admit the solution with $\bF_{1} = U_{1,\epsilon} S_{1,\epsilon} V_{1,\epsilon}^\top$; otherwise, we go back to Step 1 to further augment the xKrylov subspaces. 
 
\end{enumerate}

\subsection{{DIRK Generalization to Multi-Rank Advection-Diffusion Coefficients}}\label{DIRK_app_3}

We now extend the formulation to the multi-rank case, as described in \eqref{eq:adv-diff_sub-multirank}. Using a DIRK method for temporal discretization, the \(k\)-th stage GSE becomes:
\begin{equation}
    \mathbf{F}^{(k)} - a_{kk} \Delta t \left( 
    \underbrace{ 
    \sum_{i=1}^{\ell_x} T_{1,i} \mathbf{F}^{(k)} {\Phi_i^{2,x}}^\top 
    + \sum_{j=1}^{\ell_y} \Phi_j^{1,y} \mathbf{F}^{(k)} {T_{2,j} }^\top
    + \sum_{k=1}^{k_x} T_{3,k} \mathbf{F}^{(k)} {\Sigma_k^{2,x}}^\top 
    + \sum_{l=1}^{k_y} \Sigma_l^{1,y} \mathbf{F}^{(k)} {T_{4,l} }^\top
    }_{\mathscr{L}(\mathbf{F}^{(k)})}
    \right) = \mathbf{B}^{(k)}, 
    \quad k = 1, 2, \ldots, s
    \label{Generalized_Sylv_multi}
\end{equation}
where \(\mathbf{B}^{(k)}\) is as defined in \eqref{disc-sylv}. The terms \(T_{1:4,i}\) represent compositions of difference operators with diffusion and advection coefficients, similar to \eqref{operator_T}. Left subscripts 1 and 2 correspond to diffusion, while 3 and 4 correspond to advection in the \(x\) and \(y\) directions, respectively. The right subscript denotes the rank of the respective diffusion and advection coefficients. This equation generalizes \eqref{Generalized_Sylv}, which is recovered for \(\ell_x = \ell_y = k_x = k_y = 1\) and equal diffusion and advection coefficients in both the \(x\) and \(y\) directions. It also generalizes the backward Euler special case with a single stage and \(a_{1,1} = 1\), yielding \(\mathbf{B}^{(1)} = \mathbf{F}_0\), consistently with \eqref{Generalized_Sylv}.

Using approximations in \eqref{averaged_coeffs}, we arrive at the $k$th stage approximated SE:
\begin{equation}
 \underbrace{ {P_1} \mathbf{F}^{(k)} + \mathbf{F}^{(k)} {P_2^\top}}_{\mathscr{\Tilde{L}}\left(\bF^{(k)}\right)} = \mathbf{B}^{(k)},
    \label{Approximated_Sylv_multi_DIRK}
\end{equation}
where
\begin{align*}
\small
      P_1&= \sum_{i=1}^{\ell_x}  {A_{1,i}} + \sum_{k=1}^{k_x} {A_{3,k}}= \sum_{i=1}^{\ell_x} {\left(\frac{1}{\mathcal{R}_{\text{ranks}}} I_{N_1} - a_{kk} \Delta t \alpha_{x,i} T_{1,i}\right)} + \sum_{k=1}^{k_x} {\left(\frac{1}{\mathcal{R}_{\text{ranks}}} I_{N_1} - a_{kk} \Delta t \gamma_{x,k} T_{3,k}\right)} , \\ P_2 &= \sum_{j=1}^{\ell_y} {A_{2,j}} +\sum_{l=1}^{k_y} {A_{4,l}} = \sum_{j=1}^{\ell_y} {\left(\frac{1}{\mathcal{R}_{\text{ranks}}} I_{N_2} - a_{kk} \Delta t \alpha_{y,j} T_{2,j}\right)}+ \sum_{l=1}^{k_y} {\left(\frac{1}{\mathcal{R}_{\text{ranks}}} I_{N_2} - a_{kk} \Delta t \gamma_{y,l} T_{4,l}\right)} .
\end{align*}
In order to obtain a reduced equation for the $k$th stage matrix of coefficients \(\bS^{(k)}\), we perform a Galerkin projection \(U^\top \mathbf{R}_{\mathscr{{L}}} V = 0\) to obtain :
\begin{equation}
     \bS^{(k)} - a_{kk}\Delta_t \left( \sum_{i=1}^{\ell_x} \tilde{T}_{1,i} \bS^{(k)} \tilde{\Phi}_i^{2,x\top}
 + \sum_{j=1}^{\ell_y} \tilde{\Phi}_j^{1,y} \bS^{(k)} \tilde{T}_{2,j}^\top + \sum_{k=1}^{k_x} \tilde{T}_{3,k} \bS^{(k)} \tilde{\Sigma}_k^{2,x\top} + \sum_{l=1}^{k_y} \tilde{\Sigma}_l^{1,y} \bS^{(k)} \tilde{T}_{4,l}^\top \right) = \tilde{\bB}^{(k)}
     \label{eq:s_k_multi}
\end{equation}
where the projected operators are analogous to those defined in the rank-one case.

Afterwards, as before, we solve \eqref{eq:s_k_multi} for the \(\bS^{(k)}\) matrix using preconditioned GMRES. The preconditioner is constructed
by a Galerkin projection \(U^\top \mathbf{R}_{\mathscr{\tilde{L}}} V = 0\), to find the projected SE:
\begin{equation}
    \underbrace{\left[ \sum_{i=1}^{\ell_x} 
{\tilde{A}_{1,i}} + \sum_{k=1}^{k_x} 
{\tilde{A}_{3,k}} \right]}_{\tilde{P}_1} \mathbf{S}^{(k)} + \mathbf{S}^{(k)}\underbrace{\left[ \sum_{j=1}^{\ell_y} 
{\tilde{A}_{2,j}} + \sum_{l=1}^{k_y} 
{\tilde{A}_{4,l}} \right]}_{\tilde{P}_2^\top} = \tilde{\bB}^{(k)}
    \label{nearby_projected_multi}
\end{equation}
Using the operators $\Tilde{P}_1$ and $\Tilde{P}_2$, we define the ACS-preconditioner mapping:
   \begin{equation}
    \bS\mapsto \Tilde{P}_1\bS+\bS\tilde{P}_2^T,
\label{preconditioner_DIRK_multi}
   \end{equation}
Algorithm \ref{DIRK-LR} solves the time-dependent variable advection-diffusion equation using DIRK methods.

\begin{algorithm}
\caption[s-Stages Adaptive-Rank DIRK Integrator for the GSE]{s-Stages Adaptive-Rank DIRK Integrator for the GSE. %\footnotetext{A MATLAB implementation of this algorithm is available upon request.}
}
\label{DIRK-LR}
\SetAlgoNlRelativeSize{-2}
\SetNlSty{textbf}{(}{)}
\KwIn{
    Initial condition matrices $U_0$, $V_0$, $S_0$;\\
    Operators $P_1$, $P_2$, $\{T_{1,i}\}_{i=1}^{\ell_x}$, $\{T_{2,j}\}_{j=1}^{\ell_y}$, $\{T_{3,k}\}_{k=1}^{k_x}$, $\{T_{4,l}\}_{l=1}^{k_y}$;\\
    Butcher table $\{a_{ij}\}$;\\
    Time step size $\Delta t$;\\
    Tolerances $\epsilon_\kappa$, $\epsilon_{\texttt{GMRES}}$, $\epsilon$, ${\epsilon}_{\texttt{tol}}$;\\
    % Maximum iterations \texttt{max\_iter};
}
\KwOut{
    Updated bases $U_1$, $V_1$;\\
    Truncated singular values $S_1$;
}
\BlankLine
Compute operators $\{A_{1,i}\}_{i=1}^{\ell_x}$, $\{A_{2,j}\}_{j=1}^{\ell_y}$, $\{A_{3,k}\}_{k=1}^{k_x}$, $\{A_{4,l}\}_{l=1}^{k_y}$ according to Eq.~\eqref{Approximated_Sylv_multi_DIRK}\;
\HHONE{Compute $l=2+\ell_x+\ell_y+k_x+k_y$\;
Set \(U^{(i)}=U_0\) and \(V^{(i)}=V_0\) for \(i=1, \dots, l\)\;
Set $U_1=U_0 \quad V_1=V_0$\;
}
\While{not converged}{
   \tcp*[h]{Step 1.}\\
    \HHONE{Orthogonalize and} truncate to tolerance $\epsilon_{\kappa}$\;
\HHONE{$U_1, \{U^{(1)}\cdots U^{(l)}\} \leftarrow \kappa_m^{\texttt{trunc}}
\left(\{P_1, P_1^{-1}, A^{-1}_{1,1: \ell_x}, A^{-1}_{3,1: k_x},  \Phi^{1,y}_{1: \ell_y},\Sigma^{1,y}_{1: k_y} \}; \{U^{(1)}\cdots U^{(l)}\};U_1; \epsilon_\kappa \right)$\;}
\HHONE{$V_1, \{V^{(1)}\cdots V^{(l)}\}  \leftarrow \kappa_m^{\texttt{trunc}}\left(\{P_2,P_2^{-1}, A^{-1}_{2,1: \ell_y}, A^{-1}_{4,1: k_y}, \Phi^{2,x}_{1: \ell_x},\Sigma^{2,x}_{1: k_x}\} ;\{V^{(1)}\cdots V^{(l)}\};V_1;\epsilon_\kappa\right)$\; }
   \tcp*[h]{Step 2.}\\
  \For{$k = 1$ \KwTo $s$}{
    Compute $\tilde{\bB}^{(k)} = \tilde{\bB}^{(1)} + \Delta t \sum_{\ell=1}^{k-1} \frac{a_{k\ell}}{a_{\ell\ell}} \left( \bS^{(\ell)} - \tilde{\bB}^{(\ell)} \right)$\;
    Solve the reduced equation \eqref{nearby_projected_multi} for $\bS^{(k)}$ using GMRES with ACS preconditioner \eqref{preconditioner_DIRK_multi} to tolerance $\epsilon_{\texttt {GMRES}}$\;
        Compute and store $\frac{1}{a_{kk}} \left( \bS^{(k)} - \tilde{\bB}^{(k)} \right)$ and proceed to the next stage\;}
 {Truncate $\{U_{1,\epsilon},S_{1,\epsilon},V_{1,\epsilon}\}\leftarrow \mathcal{T}_{\epsilon}(\{U_1,S^{(s)},V_1\})$\label{alg5:line:trunc} to tolerance $\epsilon$}\;
   \tcp*[h]{Step 3.}\\
   {Compute $\{ \_, R_U \} = \texttt{QR}\left([U_1, U_{1,\epsilon}, T_{1,1:\ell_x}U_{1,\epsilon}, \Phi_{1:\ell_y}^{1,y}U_{1,\epsilon}, T_{3,1:k_x}U_{1,\epsilon}, \Sigma_{1:k_y}^{1,y}U_{1,\epsilon}]\right)$ and\\
  $\{\_, R_V\} = \texttt{QR}\left([V_1,V_{1,\epsilon} ,\Phi_{1:\ell_x}^{2,x}V_{1,\epsilon}, T_{2,1:\ell_y}V_{1,\epsilon}, \Sigma_{1:k_x}^{2,x}V_{1,\epsilon}, T_{4,1:k_y}V_{1,\epsilon}]\right)$ }\;
   \tcp*[h]{Step 4.}\\
  % Compute the residual $\|\mathbf{R}\| = \left\| R_U \text{diag}\left(S_1 - \tilde{B}^{(s)}, -a_{ss}\Delta t \left( I_{\mathcal{R}_{\text{ranks}}} \otimes S_1\right) \right) R_V^\top \right\|$\;
     {Compute the residual $\|\mathbf{R}\| = \left\| R_U \text{diag}\left(-\tilde{B}^{(s)} ,S_{1,\epsilon}, -\Delta t \left( I_{\mathcal{R}_{\text{ranks}}} \otimes S_{1,\epsilon}\right) \right) R_V^\top \right\|$\label{alg5:line:residual}}\;
     
  \If{$\|\mathbf{R}\|/\|\mathbf{F}_0\| \geq \epsilon_{\texttt{tol}}$}{
    Return to Step 1 to augment bases further\;
  }
  \Else{
      {  Set $U_1\leftarrow U_{1,\epsilon}$; $U_1\leftarrow V_{1,\epsilon}$; $S_1=S_{1,\epsilon}$\;}
         Break\;
      Exit loop\;
  }
}
\end{algorithm}

\section{Three-Dimensional Formulation}\label{Appendix_3D}

We outline the approach to 3D next. In 3D, we replace the SVD algorithm with a tensor Tucker decomposition, specifically tailored for the three (and higher) dimension tensor-product grids. 
\subsection{Tucker Tensor Equations}

We discretize the 3D time-dependent advection-diffusion equation on a three-dimensional tensor-product spatial grid with \(N_1\), \(N_2\), and \(N_3\) grid points in the \(x\), \(y\), and \(z\) directions, respectively.
%The grid points in the \(x\)-direction are denoted by \(x_i\) for \(i = 0, \dots, N_1 - 1\), with analogous definitions for the \(y\) and \(z\) directions. 
% Our goal is to evolve the initial condition \(\bF_0\) at time \(t_0\) to an updated solution \(\bF_1\) at time \(t_1\). 
%\HH{Should we omit the following sentences and simply reference some paper containing these identites?}\QQ{we could...Actually, we could move this notation part after the model (when we first need it)...} \lc{I would not omit these sentences, but would move them as Jingmei suggests} \HH{How about now?}
%We start by introducing notation and identities pertaining to the high-dimensional tensors which will be used throughout this note. 
{Here and in the following, we let} \(\bF \in \mathbb{R}^{N_1 \times N_2 \times N_3}\) denote a three-dimensional tensor. The \emph{mode-\(n\) product} of \(\bF\) with a matrix \(A_n \in \mathbb{R}^{I \times N_n}\), for \(n = 1, 2, 3\), is defined as:
\[
\bF \times_1 A_1 \in \mathbb{R}^{I \times N_2 \times N_3}, \, 
\bF \times_2 A_2 \in \mathbb{R}^{N_1 \times I \times N_3}, \, 
\bF \times_3 A_3 \in \mathbb{R}^{N_1 \times N_2 \times I}.
\]
%\QQ{suggest to apply $A$ on $\bF$, i.e. $A_1  \times_1 \bF$, and comment that it acts as operators on the corresponding dimension. Make analog to the PDE. this switch may lead to many changes throughout the manuscript...}\HH{Note: We decided to keep the notation as is to remain consistent with litterature.}
The resulting tensor in each case has its \(n\)-th dimension replaced by the number of rows of \(A_n\). The entries of the mode-\(n\) product are computed as:
\[
(\bF \times_n A_n)_{i_1, \ldots, i_{n-1}, j, i_{n+1}, \ldots, i_3} 
= \sum_{i_n=1}^{N_n} \bF_{i_1, \ldots, i_n, \ldots, i_3} \cdot (A_n)_{j, i_n}.
\]
For a more thorough review of tensor operations, we refer the reader to \cite{kolda2006multilinear}. The discretization of the 3D advection-diffusion equation  using finite-differences in space combined with the backward Euler method in time, yields the following GSE:
\begin{align}
   \bR_{\mathscr{L}} = \mathbf{F}_1 - \Delta t \mathscr{L}(\bF_1) - \mathbf{F}_0 = \mathbf{0},
    \label{Generalized_Sylv_multiterm_2}
\end{align}
where \(\bF_0\) denotes the initial tensor condition at time \(t_0\), and \(\bF_1\) is the sought tensor solution at time \(t_1\), and,
\begin{align}
 \mathscr{L} : \mathbf{F} \mapsto & 
        \sum_{i=1}^{\ell_x} \mathbf{F} \times_1 T_{1,i} \times_2 \Phi_i^{2,x} \times_3 \Phi_i^{3,x}
        + \sum_{i=1}^{\ell_y} \mathbf{F} \times_1 \Phi_i^{1,y} \times_2 T_{2,i} \times_3 \Phi_i^{3,y} \notag 
      + \sum_{i=1}^{\ell_z} \mathbf{F} \times_1 \Phi_i^{1,z} \times_2 \Phi_i^{2,z} \times_3 T_{3,i} \\
      + &\sum_{i=1}^{k_x} \mathbf{F} \times_1 T_{4,i} \times_2 \Sigma_i^{2,x} \times_3 \Sigma_i^{3,x} \notag 
      + \sum_{i=1}^{k_y} \mathbf{F} \times_1 \Sigma_i^{1,y} \times_2 T_{5,i} \times_3 \Sigma_i^{3,y}
        + \sum_{i=1}^{k_z} \mathbf{F} \times_1 \Sigma_i^{1,z} \times_2 \Sigma_i^{2,z} \times_3 T_{6,i}.
\end{align}
The operators \(T\) correspond to the composition of advection and diffusion difference operators with diagonal matrices, analogous to those defined in Section \ref{Section_2}.
Similarly, the operators \(\Phi\) and \(\Sigma\) are analogous to the diagonal matrices for diffusion and advection coefficients, respectively, used for pointwise multiplication, as defined in Section \ref{Section_2}. To obtain an averaged (approximate) mapping for preconditioning purposes, we define the average coefficients:
\begin{equation}
\begin{aligned}
    \alpha_i^{2,x} &= \text{avg}(\Phi_i^{2,x}), \,
    \alpha_i^{3,x} &= \text{avg}(\Phi_i^{3,x}), \,
    \alpha_i^{1,y} &= \text{avg}(\Phi_i^{1,y}), \,
    \alpha_i^{3,y} &= \text{avg}(\Phi_i^{3,y}), \,
    \alpha_i^{1,z} &= \text{avg}(\Phi_i^{1,z}), \,
    \alpha_i^{2,z} &= \text{avg}(\Phi_i^{2,z}), \\
    \gamma_i^{2,x} &= \text{avg}(\Sigma_i^{2,x}), \,
    \gamma_i^{3,x} &= \text{avg}(\Sigma_i^{3,x}), \,
    \gamma_i^{1,y} &= \text{avg}(\Sigma_i^{1,y}), \,
    \gamma_i^{3,y} &= \text{avg}(\Sigma_i^{3,y}), \,
    \gamma_i^{1,z} &= \text{avg}(\Sigma_i^{1,z}), \,
    \gamma_i^{2,z} &= \text{avg}(\Sigma_i^{2,z}).
\end{aligned}
\label{averaged_coeffs}
\end{equation}
which yields the three-dimensional approximate operator $\Tilde{\mathscr{L}}$ as:
\begin{align}
\mathscr{\Tilde{L}} : \mathbf{F} & \mapsto  
    \mathbf{F} \times_1 
    \bigg( 
         \sum_{i=1}^{\ell_x} \alpha_i^{2,x} \alpha_i^{3,x} T_{1,i} 
         + \sum_{k=1}^{k_x} \gamma_k^{2,x} \gamma_k^{3,x} T_{4,k} 
    \bigg) \notag \\
&+ \mathbf{F} \times_2 
    \bigg( 
        \sum_{j=1}^{\ell_y} \alpha_j^{1,y} \alpha_j^{3,y} T_{2,j} 
        + \sum_{l=1}^{k_y} \gamma_l^{1,y} \gamma_l^{3,y} T_{5,l} 
    \bigg) \notag \\
&+ \mathbf{F} \times_3 
    \bigg( 
        \sum_{k=1}^{\ell_z} \alpha_k^{1,z} \alpha_k^{2,z} T_{3,k} 
        + \sum_{l=1}^{k_z} \gamma_l^{1,z} \gamma_l^{2,z} T_{6,l} 
    \bigg).
\end{align}
Using this approximate operator, {and letting \({\mathcal{R}_{\text{ranks}}} = r_x + r_y + r_z + \ell_x + \ell_y + \ell_z\),} we then define the averaged three-dimensional approximate SE of the form:
\begin{equation}
    \bF_1\times_1 P_1+ \bF_1 \times_2 P_2 +\bF_1 \times_3 P_3= \bF_0,
\label{nearby_sylvester}
\end{equation}
with 
\begin{align}
\small
      P_1&=  \sum_{i=1}^{\ell_x} \underbrace{\left(\frac{1}{\mathcal{R}_{\text{ranks}}} I_{N_1} - \Delta t \alpha_i^{2,x} \alpha_i^{3,x} T_{1,i}\right)}_{{A_{1,i}}} + \sum_{k=1}^{k_x} \underbrace{\left(\frac{1}{\mathcal{R}_{\text{ranks}}} I_{N_1} - \Delta t \gamma_k^{2,x} \gamma_k^{3,x} T_{4,k}\right)}_{{A_{4,k}}}=\sum_{i=1}^{\ell_x}  {A_{1,i}} + \sum_{k=1}^{k_x} {A_{4,k}},\nonumber\\
      P_2 &= \sum_{j=1}^{\ell_y} \underbrace{\left(\frac{1}{\mathcal{R}_{\text{ranks}}} I_{N_2} - \Delta t\alpha_j^{1,y} \alpha_j^{3,y} T_{2,j}\right)}_{A_{2,i}}+ \sum_{l=1}^{k_y} \underbrace{\left(\frac{1}{\mathcal{R}_{\text{ranks}}} I_{N_2} - \Delta t\gamma_l^{1,y} \gamma_l^{3,y} T_{5,l} \right)}_{A_{5,i}}=\sum_{j=1}^{\ell_y} {A_{2,j}} +\sum_{l=1}^{k_y} {A_{5,l}}, \label{P_123}  \nonumber\\
    P_3 &= \sum_{j=1}^{\ell_z} \underbrace{\left(\frac{1}{\mathcal{R}_{\text{ranks}}} I_{N_3} - \Delta t\alpha_k^{1,z} \alpha_k^{2,z} T_{3,k} \right)}_{A_{3,i}}+ \sum_{l=1}^{k_z} \underbrace{\left(\frac{1}{\mathcal{R}_{\text{ranks}}} I_{N_2} - \Delta t\gamma_l^{1,z} \gamma_l^{2,z} T_{6,l} \right)}_{A_{6,i}}=\sum_{j=1}^{\ell_z} {A_{3,j}} +\sum_{l=1}^{k_z} {A_{6,l}}.
    \nonumber
\end{align}
These one-dimensional matrix operators, analogous to those in \eqref{Approximated_Sylv_multi}, are used (i) to construct xKrylov subspaces for the Tucker basis factors and (ii) to precondition the solution of the three-dimensional GSE arising from the Galerkin condition on the residual in Eq.~\eqref{Generalized_Sylv_multiterm_2}.

\subsection{GMRES-ACS Solution of the Reduced GSE  in 3D}
% We extend the adaptive-rank treatment in \cite{el2024krylov} and \cite{kahza2024sylvester} to the three-dimensional case. Specifically, we propose an adaptive-rank algorithm for the generalized Sylvester equation \eqref{Generalized_Sylv_multiterm}. 
We assume the existence of a low-rank Tucker factorization of the initial condition, \(\mathbf{F}_0 \in \mathbb{R}^{N_1 \times N_2 \times N_3}\), which evolves to an updated solution \(\mathbf{F}_1 \in \mathbb{R}^{N_1 \times N_2 \times N_3}\) at time \(t^{(1)} = t^{(0)} + \Delta t\), also in the low-rank Tucker format. This is represented as:
\begin{equation}
    \mathbf{F}_0 = S_0 \times_1 U_0 \times_2 V_0 \times_3 W_0 
    \xrightarrow{\text{evolve in time}} 
    \mathbf{F}_1 = S_1 \times_1 U_1 \times_2 V_1 \times_3 W_1,
\end{equation}
where \(U_0 \in \mathbb{R}^{N_1 \times r^x_0}\), \(V_0 \in \mathbb{R}^{N_2 \times r^y_0}\), \(W_0 \in \mathbb{R}^{N_3 \times r^z_0}\), \(U_1 \in \mathbb{R}^{N_1 \times r^x_1}\), \(V_1 \in \mathbb{R}^{N_2 \times r^y_1}\), and \(W_1 \in \mathbb{R}^{N_3 \times r^z_1}\) are orthonormal bases for their respective dimensions, and \(S_0 \in \mathbb{R}^{r^x_0 \times r^y_0 \times r^z_0}\) and \(S_1 \in \mathbb{R}^{r^x_1 \times r^y_1 \times r^z_1}\) are core tensors at the corresponding times. 
% As in \cite{kahza2024sylvester}, our method involves constructing the bases \(U_1\), \(V_1\), and \(W_1\) and then projecting the original equation to determine the coefficients of \(S_1\).
The bases \(U_1\), \(V_1\), and \(W_1\) are constructed from orthonormalization of dimension-wise xKrylov, as discussed in \ref{basis_cons}
% , derived from equations \eqref{Generalized_Sylv_multiterm} and \eqref{nearby_sylvester}, as follows: \QQ{need to be updated...}
% \begin{align}
% U_1 &= \texttt{orth}\left(\kappa_m\left(P_1, P_1^{-1}, A^{-1}_{1,1:\ell_x}, A^{-1}_{4,1:k_x}, \Phi^{1,y}_{1:\ell_y}, \Phi^{1,z}_{1:\ell_z}, \Sigma^{1,y}_{1:k_y}, \Sigma^{1,z}_{1:k_z}, U_0\right) ;\epsilon_\kappa\right), \\
% V_1 &= \texttt{orth}\left(\kappa_m\left(P_2, P_2^{-1}, A^{-1}_{2,1:\ell_y}, A^{-1}_{5,1:k_y}, \Phi^{2,x}_{1:\ell_x}, \Phi^{2,z}_{1:\ell_z}, \Sigma^{2,x}_{1:k_x}, \Sigma^{2,z}_{1:k_z}, V_0\right);\epsilon_\kappa\right), \\
% W_1 &= \texttt{orth}\left(\kappa_m\left(P_3, P_3^{-1}, A^{-1}_{3,1:\ell_z}, A^{-1}_{6,1:k_z}, \Phi^{3,x}_{1:\ell_x}, \Phi^{3,y}_{1:\ell_y}, \Sigma^{3,x}_{1:k_x}, \Sigma^{3,y}_{1:k_y}, W_0\right);\epsilon_\kappa\right).
% \end{align}

% \WT{Define in words what ${\cal T}_{\epsilon_{\kappa}}$ and $\kappa_{m}$.}\HH{Done.}
%
% The operators \(A\) are used to construct solution bases in the form of inverted Krylov subspaces. 
% For certain advection flows, particularly in swirling problems, the averaging may yield zero due to equal advection and diffusion in opposing directions. While the current methodology performs well, our numerical experiments suggest that alternative averaging metrics, such as the \(L_1\)-norm, can further enhance basis construction.
%
Using a Galerkin projection on the residual, \(\mathbf{R}_\mathscr{L} \times_1 U_1^\top \times_2 V_1^\top \times_3 W_1^\top\), we obtain a reduced-size GSE for the core tensor \(\bS_1\) as follows:
\begin{align}
    \mathbf{S}_1 - \Delta t \mathscr{L}(\mathbf{S}_1) \times_1 U_1^\top \times_2 V_1^\top \times_3 W_1^\top 
    = \mathbf{S}_0 \times_1 U_1^\top U_0 \times_2 V_1^\top V_0 \times_3 W_1^\top W_0.
    \label{reduced_S}
\end{align}
This projected system has dimensions \( r_x \times r_y \times r_z \), where \( r_x \), \( r_y \), and \( r_z \) denote the dimensions of the truncated xKrylov subspaces in the \( x \), \( y \), and \( z \) directions, respectively, and will be solved using the GMRES-ACS preconditioner.
%\QQ{add a comment on the reduced size of the system now?}
 % can be solved directly \cite{simoncini2020numerical, chen2020recursive}

% In the constant-coefficient case \cite{el2024krylov}, the approximated operator \( \mathscr{\Tilde{L}}(\bF) \) is equal to the original operator \( \mathscr{L}(\bF) \), and Eq. \eqref{reduced_S} can be solved using standard Sylvester solvers, eliminating the need for GMRES. 
Eq. \eqref{reduced_S} is solved iteratively using GMRES, with $\mathscr{L}$-operator application cast as a series of tensor-matrix multiplications, along with the ACS preconditioning operator given by:
\begin{align}
  \mathcal{P} : \mathbf{S} \mapsto \mathbf{S} \times_1 \tilde{P}_1 + \mathbf{S} \times_2 \tilde{P}_2 + \mathbf{S} \times_3 \tilde{P}_3,
  \label{nearby_projected_3d}
\end{align}
where \(\Tilde{P}_{i} = P_i \times_1 U_1^\top \times_2 V_1^\top \times_3 W_1^\top\) denotes the projection of operators in Eq. \eqref{nearby_sylvester}, analogous to those introduced in Section~\ref{Section_3}. {This operator is used to precondition each GMRES iteration of the inner system solve.} The solution of the {projected} SE has a computational cost of \(\mathcal{O}(r^4)\) \cite{chen2020recursive,simoncini2020numerical}. The tensor-times-matrix operations, executed using the \texttt{ttm} command in the Tensor Toolbox, incur a cost of \(\mathcal{O}(r^4)\). The total complexity in three dimensions is therefore \(\mathcal{O}(N r^2 + r^4)\).

\section{Technical Details of the Hypothetical Multigrid Comparison}
{
\color{black}
\label{Appendix_MG}
We provide here a brief explanation of the assumptions involved in the performance comparison between a hypothetical optimally scaling full-rank multigrid method and the adaptive-rank approach. 
Optimal multigrid methods exhibit a computational scaling of $\mathcal{O}(N^d \log N)$ \cite{jones1997parallel}. Let us denote the wall-clock time (WCT) of a multigrid method as: \[
WCT_{\mathrm{MG}} = C_{\mathrm{MG}} N^d \log N.
\]  
From Fig. \ref{fig:convergence_comparison}-a, the wall-clock time of the adaptive-rank method may be written as:
\[
WCT_{\mathrm{LR}} = C_\mathrm{LR} (0.0023\,N + 55.2734).
\]  
The speedup of the adaptive-rank method vs. multigrid is given by the ratio:
\begin{equation}
\frac{WCT_{\mathrm{MG}}}{WCT_{\mathrm{LR}}} = C_{\mathrm{WCT}} \frac{ N^3 \log N}{0.0023\,N + 55.2734},
\label{WCT_expr}
\end{equation}
where $C_{\mathrm{WCT}} = 1.2 \times 10^{-5}$ is found assuming that both methods achieve comparable wall-clock times for $N = 100$.
This formula is used to populate Table~\ref{tab:wct_ratio}.

}

\bibliographystyle{abbrv}
\bibliography{sample}

\end{document}